\documentclass{article}
\usepackage{amssymb}
\usepackage{latexsym}
\usepackage{amsthm}
\usepackage{enumerate}
\usepackage{epsfig}
\usepackage{color}
\usepackage{float}
\usepackage{subfigure}
\usepackage{amsmath}
\usepackage{ifthen}
\usepackage{makeidx}
\usepackage{lscape}

\newboolean{printing_margins}
\setboolean{printing_margins}{false}

\newboolean{laptop_margins}
\setboolean{laptop_margins}{false}

\ifthenelse{ \boolean{printing_margins} }
{ 
\setlength{\oddsidemargin}{-0.5in}
\setlength{\evensidemargin}{-0.5in}
\setlength{\textwidth}{7.5in}
\setlength{\topmargin}{-0.5in}
\setlength{\textheight}{9.5in}
}
{ 
\setlength{\oddsidemargin}{.25in}
\setlength{\evensidemargin}{.25in}
\setlength{\textwidth}{6in}
}

\ifthenelse{ \boolean{laptop_margins} }
{ 
\setlength{\textheight}{7in}
}{}

\def\a{\alpha}

\def\g{\gamma}

\def\R{\mathbb{R}}

\def\Z{\mathbb{Z}}

\def\d{\partial}
\def\cross{\times}

\def\orbit{G \centerdot}

\def\half{\frac{1}{2}}

\newcommand{\geom}[2]{\text{geom}(#1 \cap #2)}

\newtheorem{theorem}{Theorem}[section]
\newtheorem{lemma}[theorem]{Lemma}
\newtheorem{corollary}[theorem]{Corollary}

\newenvironment{Proof}{\noindent \textbf{Proof.}}{\qed\vspace{.5\baselineskip}}

\theoremstyle{definition}
\newtheorem{definition}[theorem]{Definition} 
\newtheorem{notation}[theorem]{Notation}
\newtheorem{remark}[theorem]{Remark}
\newtheorem{claim}[theorem]{Claim}

\newtheorem{convention}[theorem]{Convention}

\newcounter{case}

\newenvironment{case}[1]{\stepcounter{case} \addvspace{.5\baselineskip} \noindent\textbf{Case \thecase}. \textit{#1}}{\hfill\fbox{Case \thecase}}

\newtheorem{subcase}{Case}[case]
\newtheorem{sub2case}{Case}[subcase]

\newcounter{step}

\restylefloat{figure}

\begin{document}

\title{Period three actions on lens spaces}
\author{Joseph Maher\footnote{email:maher@math.okstate.edu}}
\date{\today}

\maketitle

\tableofcontents

\section{Introduction}

In this paper we prove the following theorem:

\begin{theorem} \label{theorem:main}
A free $\Z_3$ action on a lens space is standard.
\end{theorem}

An immediate corollary of this is that a free $\Z_{3^k}$ action on $S^3$ is standard. Together with earlier work of Milnor \cite{milnor}, Livesay \cite{livesay}, Myers \cite{myers}, Rubinstein \cite{rubinstein} and Thomas \cite{thomas} this implies that:

\begin{corollary}
A free action of a group of order $2^a 3^b$ on $S^3$ is standard.
\end{corollary}

This is a special case of the spherical spaceform conjecture, which
has now been resolved by Perelman's proof of Thurston's geometrization
conjecture \cite{mt}.

This paper is an extension of Maher and Rubinstein \cite{mr}, which dealt with free period-three actions on the three sphere, and is not self contained. The basic strategy is the same as the three sphere case. We first explain the strategy, and indicate the extra work we need to do to make it work for lens spaces. We then prove the extra results, but do not repeat parts that we can quote directly from \cite{mr}.

\begin{remark}
We show that the quotient is Seifert fibered, but in fact it is always a lens space. Given the list of fundamental groups of elliptic $3$-manifolds, it is an elementary but tedious exercise to check that none of these groups have an index three normal cyclic subgroup, unless the initial group was in fact cyclic.
\end{remark}

\subsection{Acknowledgements}

I would like to thank Daryl Cooper and Hyam Rubinstein for their
advice and encouragement. I would also like to thank the referee for
his exceptional diligence and much helpful advice, which has greatly
improved the paper. This paper was typeset in \LaTeX, and all the
figures were drawn in \verb|xfig|.

\subsection{Outline}

An action of $\Z_3$ on a lens space $L$ is given by a diffeomorphism $g:L \to L$, which is free, and which is period three. This means that $g$ generates a group $G \cong \Z_3$, and as $g$ is free, the quotient $L/G$ is a manifold. We will give $L$ the round metric coming from its description as a cyclic quotient of the round three-sphere. We say $g$ is linear, if it is covered by a linear map on $S^3$, namely an element of $SO(4)$. If $g$ is linear then the quotient of $L$ by $G$ is again a lens space. We say the action of $G$ is standard if $g$ is conjugate by a diffeomorphism to a linear map.

We can show that the action is standard by finding an invariant Seifert fibering of $L$. We will find an invariant Seifert fibering by studying sweepouts of $L$. A sweepout of $L$ is a family of surfaces which ``fill up'' the manifold. A simple example is the foliation of $L$ by tori with two singular leaves which are circles, coming from a genus one Heegaard splitting of $L$. Think of the leaves as parameterised by time, starting with one singular leaf at $t = 0$ and ending with the other one at $t = 1$. At a non-singular time $t$, the sweepout consists of a single torus $S_t$. We can think of the union of the leaves as a $3$-manifold, in this case a lens space, with a height function for which each level set is a torus. The map from the leaf space to $L$ is degree one, and is an embedding on each level set of the height function.

For our purposes, we require a more general definition, in which the sweepout surfaces at non-singular times may be finitely many spheres, together with at most one Heegaard torus. We say a generalised sweepout is a $3$-manifold $M$, with a height function $h$ on it, so that the level sets at regular values are at most one torus, together with a union of $2$-spheres, and a degree one map $f: M \to L$, which is an embedding on each level set of the height function. We shall think of the height function on $M$ as time. The $3$-manifold $M$ is in fact a connect sum of lens spaces and copies of $S^2 \cross S^1$.

We can look at the three images of the sweepout surfaces under the group $G$. Generically, they will intersect in double curves and triple points. In order to distinguish the three images of the surfaces under $G$, we colour them using the same colour conventions as \cite{mr}, which we now explain. The initial images of the surfaces $h^{-1}(t)$ in $L$ will be labelled red. The images of the red surfaces under $g$ will be labelled green, and the images of the red surfaces under $g^2$ will be labelled blue. If two surfaces intersect in a double curve, we will label the double curve by the complementary colour, i.e.  the colour of the surface not involved in the intersection. For example, if a red surface intersects a green surface in a double curve, we will label that double curve blue.

By general position, we can arrange that the sweepout surfaces intersect transversely for all but finitely many times, and that the non-transverse intersections all come from a finite list of possibilities, corresponding to critical points of the height function, on either $M$ itself, or on the double or triple point sets of $M$. Critical points of the height function on $M$ change the surfaces, either by increasing or decreasing the number of $2$-spheres, or by splitting the torus into a $2$-sphere, or the reverse operation. Critical points of the height function on the double set change the number of double curves, by either creating or destroying a double curve, or by saddling curves together. Critical points of the height function on the triple set change the number of triple points. In fact the number of triple points is always a multiple of six, as every triple point has three images under $G$, and there must also be an even number of triple points.

We call these non-transverse intersections {\bf moves}, and we can describe a sweepout by drawing the configurations in $L$ in between the critical times. Each neighbouring pair of pictures will differ by one of the moves described above. We look for a sweepout that is ``simple'', by defining a complexity for sweepouts, and showing how to change the sweepout to reduce complexity. We say that the {\bf complexity} of the sweepout at a generic time $t$ is the ordered pair $(n, d)$, where $n$ is the number of triple points, and $d$ is the number of inessential simple closed double curves, i.e. those double curves without triple points which bound discs in a sweepout surface. We order the pairs $(n,d)$ lexicographically. We say that the complexity of a sweepout is the maximum complexity that occurs over all generic times. We say that a sweepout is a minimax sweepout if it has minimal complexity.

The double curves form an equivariant graph in the lens space. If this graph contains an invariant knot or link, which is isotopic to a knot or non-trivial link lying on a Heegaard torus, then the quotient manifold is Seifert fibered, so the action is standard. The basic strategy is to show that a minimax sweepout has a local maximum which contains such an invariant knot or link.

We can change a sweepout by using a procedure we call a {\bf modification}. This involves taking an equivariant neighbourhood $N \subset L$, and a time interval $I$, so that the intersection of $\d N$ with the sweepout only changes by an isotopy during the time interval $I$. We can replace the sweepout surfaces inside $N$ during the time interval $I$ with some new set of sweepout surfaces with the same boundary. This does not change the degree of the sweepout map $f:M \to L$, but we need to make sure that the new sweepout surfaces are still spheres, together with at most one Heegaard torus. We may change the moves that occur in $N$ during the time interval $I$.

We can use modifications to remove innermost inessential double curves for short time intervals, by ``pinching them off'' at the beginning of the time interval $I$, and then putting them back at the end of the time interval. Similarly, we can remove bigons that contain no double curves in their interiors for a short time interval, by pushing the double curves forming the boundary across the bigon, and then putting them back at the end of the interval. These modifications are described in detail in \cite[Section 2.6]{mr}. We call the moves that we insert in the new sweepout, at the beginning and end of the time interval $I$, {\bf compound moves}. 

This means that if some local maximum occurs for which there is a disjoint inessential innermost double curve or bigon with no double curves in its interior, then we can reduce the complexity of the local maximum by removing the double curve or bigon for the duration of the local maximum. We call this {\bf undermining} the local maximum.

The main argument in the proof works by changing all the local maxima to local maxima consisting of compound local maxima, and then showing that we can either undermine the compound local maxima, or find an invariant knot or non-trivial link that lies on a Heegaard torus, thus showing that the action is standard. There are three main steps in doing this, and this is where the main differences from \cite{mr} occur. The three main steps are:

\begin{itemize}

\item[4.1] Every sweepout contains triple points, or is standard.
\item[5.1] There are disjoint bigons for non-triple point moves.
\item[3.1] Local maxima consisting of compound moves can be undermined, or else the action is standard.

\end{itemize}

In Lemma \ref{lemma:triple}, we show a sweepout contains triple points, or else the action is standard. This argument is similar in spirit to the argument in Section $4$ of \cite{mr}, and only minor changes are needed to deal with the fact that there may be sweepout surfaces that are tori as well as spheres.

In order to convert all local maxima into compound local maxima, we need to show that there are disjoint bigons for non-triple point moves, which we do in Lemma \ref{lemma:main}. This argument is more complicated than in the $3$-sphere case, and takes up the bulk of this paper. One reason for this is that on a $2$-sphere any pair of intersecting curves create at least four bigons, whereas on a torus, curves may intersect without creating any bigons. We break the proof into two main parts. First we show that a configuration with triple points contains at least three bigons. This suffices to show that there are disjoint bigons for all the non-triple point moves except for saddle moves. We then show there are disjoint bigons for saddle moves, by using the fact that there are at least three bigons in the configuration before and after the saddle move.

To show there are at least three bigons, first we show in Section \ref{section:intersection} there can't be too many parallel curves in a configuration with only two bigons. If there are few bigons, but many parallel curves, then there must be many pairs of curves which bound annuli which only contain parallel essential double arcs. We call these annuli {\bf strips}. The red Heegaard torus divides the lens space into two solid tori, which contain properly embedded subsurfaces of the blue and green surfaces. We show that if there are only two bigons, there cannot be intersecting strips on the same side of the red torus, and so this restricts the number of parallel double curves that may occur. 

If there are at most two bigons, then there may be at most one innermost inessential double curve of each colour, so the double curves of a given colour divide the torus up into one of three configurations. There may be a punctured torus component, a pair of pants component, or all components may be annuli. The blue and green double curves may divide the torus up in different ways. As there may not be too many parallel curves, there are a limited number of cases that may arise, and we go through each one in turn, and show that configurations with two bigons may not occur, usually by showing that the number of triple points in the configuration is not a multiple of six. We give a more detailed outline of this at the beginning of Section \ref{section:disjoint bigons}.

Finally, we need to show Lemma \ref{lemma:special cases}, we can either undermine local maxima with compound triple point moves at each end, or else the action is standard. This is straight forward, as there are only a few extra cases to consider that are not dealt with in Section $6$ of \cite{mr}, so in fact we deal with this lemma first out of the three lemmas, in Section \ref{section:special cases}.

\section{The main argument.}

\subsection{Preliminaries.}

Our definitions are the same as in \cite[Section 2]{mr}, except we now use a sweepout which may contain tori, as well as spheres. At any given time during the sweepout, there will be at most one torus, which will be a Heegaard torus.

\begin{definition}
A {\bf generalised sweepout} is a triple $(M, f, h)$, where

\begin{itemize}

\item $M$ is a closed, orientable $3$-manifold. 

\item The smooth map $h:M \to \R$ is a height function, such that for all but finitely many $t \in \R$, the inverse image, $h^{-1}(t)$ consists of a finite collection of $2$-spheres, union at most one torus. The map $h$ should be a Morse function away from the first and last singular sets, which should each consist of a single circle.

\item The smooth map $f:M \to L$  is degree one. 

\item The map $f|_{h^{-1}(t)}$ is an embedding on the level set $h^{-1}(t)$ for every $t$. If $h^{-1}(t)$ contains a torus, the image of this in $L$ is a Heegaard torus.

\end{itemize}

We will often write a sweepout as $(M, \phi)$, where $\phi$ denotes the map $(f \cross h):M \to L \cross \R$. We will think of $t \in \R$ as the time coordinate.
\end{definition}

\begin{remark}
The manifold $M$ is homeomorphic to a connect sum of lens spaces, and copies of $S^2 \cross S^1$. 
\end{remark}

\begin{convention}
As we are dealing with tori, we will often need to draw diagrams of tori. All red squares in the diagrams we draw will be assumed to have opposite sides identified to make tori.
\end{convention}

The following well known definitions are additional to those in \cite{mr}.

\begin{definition} {\bf Lens space.}

A {\bf lens space} $L$ is a quotient of $S^3$ by a cyclic subgroup of $SO(4)$ which acts freely. 
\end{definition}

\begin{definition}{\bf Free action.} 

We say that a {\bf free action} of the group $G \cong \Z_3$ on the lens space $L$ is generated by the diffeomorphism $g:L \to L$ if $g$ has no fixed points, and $g^3$ is the identity. Therefore $g$ generates a cyclic group of order three, $G = <\!g\!> \cong \Z_3$.
\end{definition}

We start with some free action of $Z_3$ on a lens space. The lens space $L$, and the diffeomorphism $g$ are fixed for the duration of this paper, and do not change.

\begin{definition} {\bf Standard action.}

The action of $G$ is {\bf standard} if $G$ is conjugate by a diffeomorphism of $L$ to an action by isometries.
\end{definition}

This is equivalent to the quotient space being Seifert fibered. We now give some conditions which suffice to show that the action is standard. The arguments below are well known, see \cite{myers} for example. 

\begin{lemma} \label{lemma:invariant unknot}
If there is a smooth invariant curve in $L$ which bounds an embedded disc, then the action of $G$ is standard.
\end{lemma}

\begin{Proof}
Let $N$ be an invariant regular neighbourhood for the smooth invariant curve. The complement $L - N$ has a compressing disc, as the curve bounds an embedded disc, so by the loop theorem the quotient $(L - N)/G$ also has a compressing disc. Cutting the torus along this disc creates a $2$-sphere, which bounds a $3$-ball as $L$, and hence $L/G$, is irreducible. The invariant loop cannot be contained in the $3$-ball, as the $3$-ball lifts to three disjoint $3$-balls in $L$, so the $2$-sphere bounds a $3$-ball on the other side. This means the complement is a $3$-ball union a $1$-handle, which is a solid torus, so $L/G$ is the union of two solid tori, namely a lens space. 
\end{Proof}

\begin{lemma} \label{lemma:invariant curve on torus}
If there is a smooth invariant curve in $L$ which lies on a Heegaard torus, then the action of $G$ is standard.
\end{lemma}

\begin{Proof}
By Lemma \ref{lemma:invariant unknot}, we may assume that the curve does not bound a disc in the Heegaard torus, and is not meridional on either side. Then $L - N$ is Seifert fibered, with base orbifold a disc, and at most two singular fibers. By \cite{gh}, the quotient is also a Seifert fibered space with base orbifold a disc. If the meridian of $N/G$ is not isotopic to a fiber of $(L - N)/G$ then we can extend the Seifert fibering to all of $L/G$, so the action of $G$ is standard. If the meridian of $N$ is isotopic to a fiber, then if $(L - N)/G$ has two or more singular fibers, then $\pi_1(L/G)$ is a non-trivial free product, which can't be a finite quotient of $S^3$. This means $(L - N)/G$ may have at most one singular fiber, so is a solid torus, so $L/G$ is in fact a lens space.
\end{Proof}

\begin{definition}
A link in $L$ is trivial if all components of the link bound discs which are all disjoint.
\end{definition}

\begin{lemma} \label{lemma:link}
If there is an invariant non-trivial link in $L$ which lies on a Heegaard torus, then the action of $G$ is standard. 
\end{lemma}

\begin{Proof}
Suppose no component of $L$ is an essential non-meridional slope on the Heegaard torus $H$. Then every component of $L$ bounds a disc, either in the Heegaard torus $H$, or in one of the solid torus complements of $H$. If components of $L$ bound meridional discs, then we may choose these meridional discs to be disjoint from each other, and from all the other components of $L$. If components of $L$ are inessential in the Heegaard torus $H$ then they bound discs in $H$. If these discs are nested, then we can make them disjoint by pushing them off to one side of $H$, so in fact all the components of $L$ bound disjoint discs, and so the link is trivial. So we may assume that $L$ has at least one component which is an essential non-meridional disc. By passing to an invariant sublink, we may assume that all the components of the link are parallel to some non-meridional slope on the Heegaard torus. The complement of a regular neighbourhood of the link is Seifert fibered, with base orbifold a pair of pants, and two singular fibers. By \cite{gh}, the quotient $(L - N)/G$ is also Seifert fibered, so as the fibers are not meridional, we can extend the Seifert fibering to all of $L/G$, so the action is standard.
\end{Proof}

We now  show that if there are essential double curves on the torus component of the sweepout surfaces, and no triple points, then the action is standard.

\begin{lemma} \label{lemma:parallel}
If all double curves are essential in the torus component of the sweepout surfaces, and there are no triple points, then the action of $G$ is standard.
\end{lemma}

\begin{Proof}
If the double curves are not meridians on either side of one of the Heegaard tori, then we can apply Lemma \ref{lemma:link}, so we may assume that the curves have meridional slope on one side of the Heegaard torus.

Each torus is divided up into annuli by the essential curves. So each side of the red torus contains blue and green annuli with parallel boundaries, and these annuli intersect only in double curves which are also parallel to their boundaries. This divides each solid torus up into regions which are solid tori. All of these solid tori have meridional discs that intersect each double curve exactly once, except for one on each side of the Heegaard torus. The number of times a meridional disc intersects each double curve is preserved under $G$, so the solid tori with meridional discs with intersection numbers different from one are invariant under $G$. However, these regions are regular neighbourhoods of the core of the Heegaard tori, and either one gives an invariant Heegaard splitting, so the action is standard. 
\end{Proof}

\subsection{The main argument.}

The argument from \cite[Section 3]{mr} goes through unchanged, assuming the following three lemmas for sweepouts of tori, which are the analogues of Lemmas 4.1, 5.1 and 6.1 from \cite{mr}.

The three lemmas are:

\begin{itemize}

\item[\ref{lemma:triple}:] Every sweepout contains triple points, or else the action is standard.

\item[\ref{lemma:main}:] There are disjoint bigons for non-triple point moves

\item[\ref{lemma:special cases}:] We can either undermine local maxima with compound triple point moves at each end, or else the action is standard.

\end{itemize}

We now give an extremely brief summary of the main argument from \cite[Section 3]{mr}. Lemma \ref{lemma:main} implies that we can undermine non-triple point moves. We can then use the genuine triple point move trick from \cite[3.2, step 2]{mr} to replace the top layer of local maxima with compound local maxima. The local maxima with compound double curve moves can be undermined as in \cite[3.2, step 3]{mr}. This leaves the special case local maxima, namely local maxima which have compound triple point moves at each end, and no other moves in the middle. These can be either undermined, or we can find an invariant Seifert fibering by Lemma \ref{lemma:special cases}. Finally, Lemma \ref{lemma:triple} shows that if we can undermine all special local maxima with triple points, then the action is standard, completing the proof.

\section{Special Cases} \label{section:special cases}

Recall that a special case local maximum is a local maximum consisting of two compound triple point moves. All the cases which are dealt with in \cite[Section 6]{mr} arise as before, and we do not repeat this material here. We only show how to deal with the extra cases that arise when some of the sweepout surfaces are tori. We will need the following definition from \cite{mr}.

\begin{definition}{\bf Special case modification neighbourhood.}

Suppose $(M, \phi)$ is a sweepout containing a special case local maximum. Suppose that $N_I$ is a modification neighbourhood with the following properties:

\begin{itemize}

\item $N_I$ contains the move neighbourhoods for the compound triple point moves in its interior, and is disjoint from all the other move neighbourhoods of the sweepout. 

\item $N_t$ is a tubular neighbourhood for the union of the two bigons in the compound moves, for the times $t$ in between the compound triple point moves. 

\end{itemize}

Then we say that $N_I$ is a {\bf special case modification neighbourhood} for the special case local maximum.
\end{definition}

\begin{lemma} \label{lemma:special cases}
A special case local maximum can be undermined, or else the action of $G$ is standard.
\end{lemma}

\begin{Proof}
Let $N_I$ be the special case modification neighbourhood. If the components of $N_t$ are all $3$-balls, then the argument of \cite{mr} goes through as before. We need to deal with the case that $N_t$ has components which are not $3$-balls. If the bigons have only three vertices in common, then the components are all $3$-balls. If the bigons have six vertices in common, and an edge in common, then the components of $N_t$ are still $3$-balls. The remaining cases are when the bigons have six vertices in common, and no edges in common. There are two cases, depending on whether $N_t$ is connected, or whether it has three connected components.

\setcounter{case}{0}
\begin{case}{$N_t$ is connected.}

In this case $N_t$ is an invariant solid torus in which the bigons share all six vertices in common, for times during the local maximum between the compound moves. The orbit of the green edges from the red bigons is an invariant curve. The green and blue edges lie in the red torus, and the red edges can be isotoped into the red torus across the blue bigons, so the invariant curve is isotopic to a curve that lies on a Heegaard torus. So by Lemma \ref{lemma:invariant curve on torus} the action of $G$ is standard.
\end{case}

\begin{case}{$N_t$ has three components.}

The intersection of the modification neighbourhood with the red surfaces is a subsurface of the red surfaces, whose boundary consists of simple closed curves. We first show that if any of these simple closed curves in $\d N_t \cap S_t$ bound a disc in a red surface, then we can undermine the local maximum. We have done all the cases previously in \cite[Section 6.2.2]{mr}, except the following one, in which the bigons involved in the compound triple points are coplanar, and the component of $N_t$ containing the red bigons intersects the red surface in a subsurface which has a boundary component which bounds a disc, $D$ say, in the red surface.

\begin{figure}[H]
\begin{center}
\epsfig{file=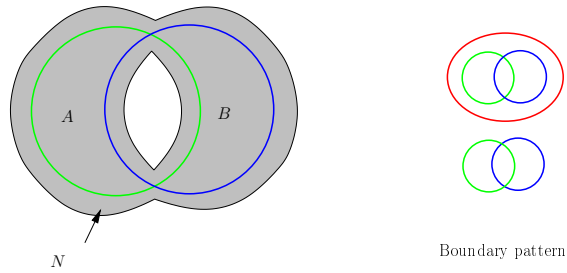, height=100pt}
\end{center}
\caption{The red bigons $A$ and $B$ are coplanar and share common vertices.}
\end{figure}

If either of the other components of $N_t$ intersect the disc $D$, then they intersect the disc $D$ in subsurfaces, which have at least one boundary component which also bounds a disc, and we have dealt with these cases already in \cite[Section 6.2.2]{mr}. There may be no triple points or simple closed curves inside the disc, as we could use them to undermine the local maximum, so $N_t$ union a regular neighbourhood of the disc is a ball with saddle reducible-boundary pattern, so we may undermine the local maximum by Lemma \cite[Lemma 6.11]{mr}

This means that all three components of $N_t$ intersect the red torus in essential annuli. If these annuli are not meridional on either side, then the action of $G$ is standard, by Lemma \ref{lemma:link}. The remaining case to consider is when the slope of the annuli is meridional on one side of the Heegaard torus.

\begin{figure}[H]
\begin{center}
\epsfig{file=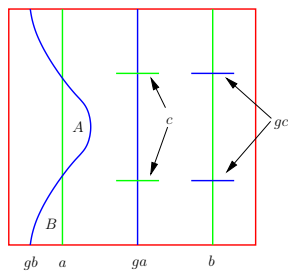, height=100pt}
\end{center}
\caption{The neighbourhood $N_t$ intersects the red torus in essential annuli.}
\label{picture532}
\end{figure}

Figure \ref{picture532} above shows the double curves that lie inside the meridional annuli. There may be no other triple points, though there may be other double curves. As the essential curves have the same slope, each pair of curves intersects in an even number of points, so the images of the triple points lie on a common curve, which we shall call $c$.

Consider the configuration after the bigon has been removed. There must be some essential curves, as the double curves that originally formed the boundary of the red bigons will both be essential. If a green double curve and its blue image are both inessential in the red torus, then we can pinch off the curve without cutting the torus into a sphere. If a curve is inessential in the red torus, but has an image which is essential, then the essential curve must bound a meridian disc which contains no double curves in its interior. So there is a meridian disc parallel to this one, whose boundary is disjoint from the double curves on the red sphere, and whose interior is disjoint from the sweepout surfaces, and which is disjoint from its images under $G$. If we surger along this disc for a time interval longer than the local maximum, we have changed the local maximum to one happening on spheres, which we can undermine. So we may assume that all of the curves are essential, and so the action is standard by Lemma \ref{lemma:parallel}. 
\end{case}

This completes the proof of Lemma \ref{lemma:special cases}.
\end{Proof}

\section{Every sweepout contains triple points.}

In this section we prove the following lemma.

\begin{lemma} \label{lemma:triple}
If there is a sweepout without triple points then the action is standard.
\end{lemma}

At the beginning of the sweepout there is a clear connected invariant region, which we will call an initial region. At the end of the sweepout, there are no clear regions. We will show that if there are no triple points, then a move which breaks up the initial region occurs in a configuration which contains curves which are essential in the torus, and in which we can remove the inessential curves without cutting the torus into a sphere. This enables us to apply Lemma \ref{lemma:parallel} to show that the action is standard.

\begin{definition}{\bf An initial region.}
\index{initial region} \index{region!initial}

We say that a region is {\bf initial} if it is  clear, connected and invariant under $G$. 
\end{definition}

We now prove Lemma \ref{lemma:triple}. \\

\begin{Proof}
At the beginning of the sweepout, there is an initial region. This is because before the first appear move, there are no sweepout surfaces, and so all of $L$ is a clear connected invariant region, and so is an initial region.

We now consider each type of move in turn, and show that if the move breaks up an initial region, then the action is standard. A single move is supported in the orbit of a $3$-ball which is disjoint from its images, so a single move cannot eliminate an initial region by shrinking it down to a point. However a single move might eliminate an initial region by disconnecting it.

\setcounter{case}{0}
\begin{case}
{Appear and vanish moves.}

If a double curve free $2$-sphere orbit appears inside an initial region, then the interior of the $2$-sphere orbit is a new region, namely a $3$-ball orbit coloured by the same colour as the sphere. The region on the outside of the $2$-sphere orbit is still connected, and so is still an initial region. Therefore appear and vanish moves may not eliminate an initial region.
\end{case}

\begin{case}
{Cut and paste moves.}

A paste move does not disconnect any region, so we need only consider cut moves. If a cut move disconnects an initial region $C$, then the cut disc $D$ is contained inside $C$. The images of $D$ also lie inside $C$, so we can find a path $\gamma$ between two images of $D$, which does not intersect any other sweepout surface, and which does not hit the third image of $D$. Some image of this path under $G$ is an arc $\g$ from $D$ to $gD$. There are two cases, depending on whether or not the curve $\orbit \g$ intersects the cut discs transversely.

\begin{figure}[H]
\begin{center}
\epsfig{file=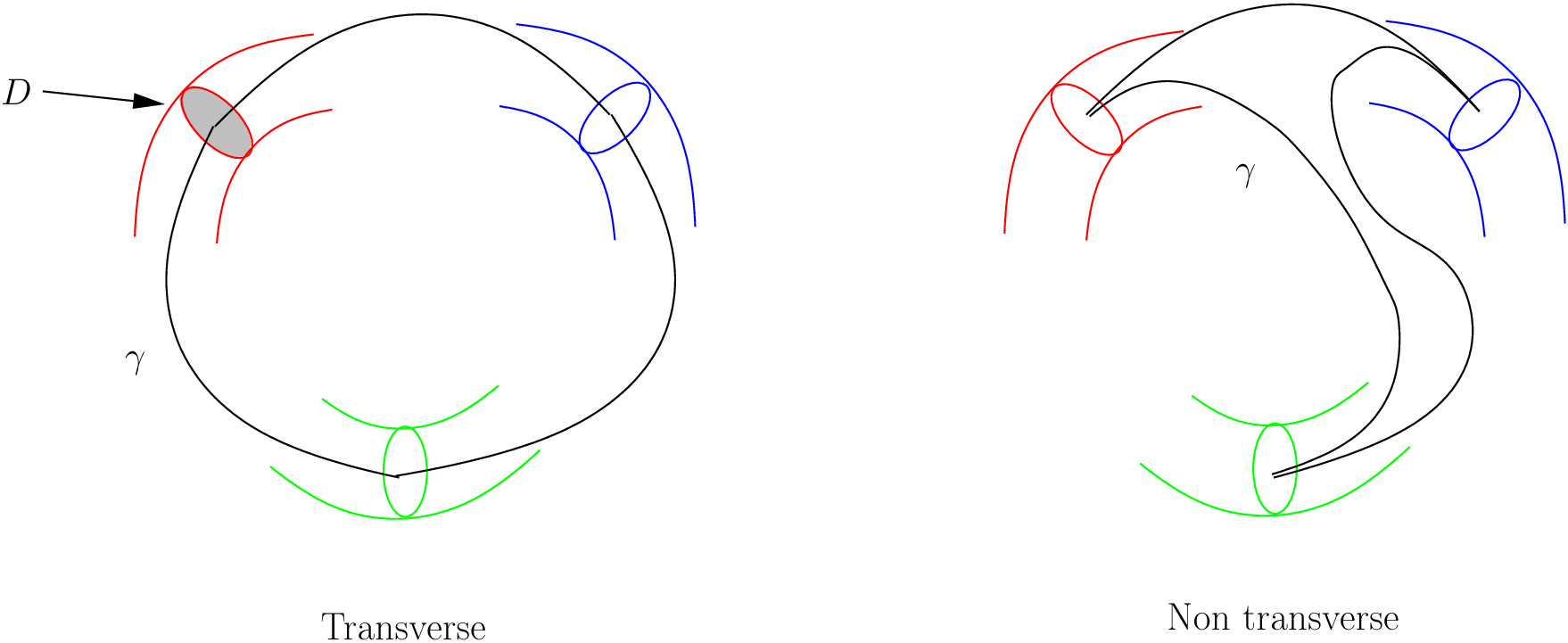, height=150pt}
\end{center}
\caption{A cut move}
\end{figure}

In the non-transverse case the region containing $\orbit \gamma$ is still an initial region after the cut move. In the transverse case, the boundary of the disc $D$ is either inessential, or essential in the red torus. If the boundary of $D$ bounds a disc $D'$ in the red sweepout surfaces, then $D \cup D'$ is a sphere. As the simple closed curve $\orbit \gamma$ lies in the clear region it is disjoint from the red surfaces, so it intersects this $2$-sphere precisely once transversely, a contradiction, as $L$ is irreducible, so this case cannot occur.

So we may assume that the boundary of $D$ is essential in the red torus. If there are any inessential double curves, then they bound discs in the sweepout surfaces disjoint from the orbits of the cut disc, so we can use cut and death moves to remove these without disconnecting the initial region. Therefore we may assume that there are no double curves that bound discs in any sweepout surface. If there are only essential curves, then the the action is standard by Lemma \ref{lemma:parallel}. If there are no double curves, then the cut move does not disconnect the initial region, as there is a path from one side of the cut disc to the other, parallel to the red torus.
\end{case}

\begin{case}
{Double curve births and deaths.}

\begin{figure}[H]
\begin{center}
\epsfig{file=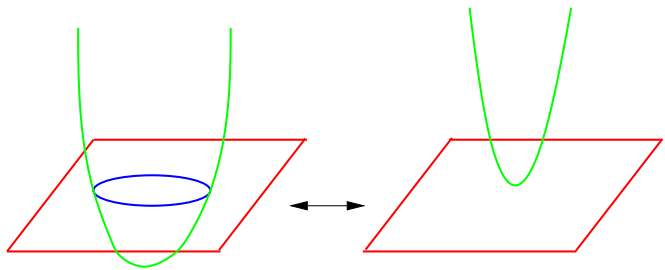, height=100pt}
\end{center}
\caption{Birth/death of a blue double curve}
\end{figure}

The region that is created or destroyed in the move cannot be invariant. The other regions that intersect the move neighbourhood are not disconnected by the move. 
\end{case}

\begin{case}
{Saddle moves.}

\begin{figure}[H]
\begin{center}
\epsfig{file=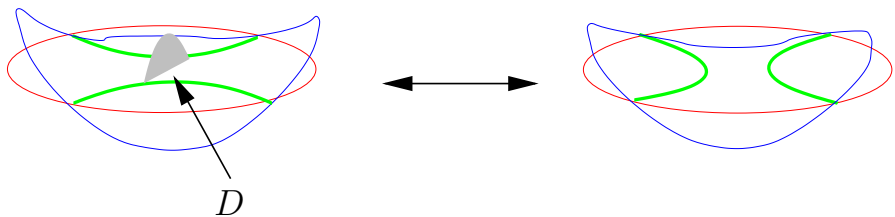, height=80pt}
\end{center}
\caption{A saddle move}
\end{figure}

The only region which could become disconnected as a result of a saddle move is the region containing the saddle disc $D$ in the diagram above. If the disc $D$ intersects two distinct green circles, then as there are no triple points, we can choose a path from one side of $D$ to the other, in the interior of the initial region, parallel to either one of the green double curves. So in this case, the saddle move does not disconnect the initial region.

Now suppose that the disc $D$ intersects a single green double curve, $\a$ say. The disc $D$ has three disjoint images under $G$. The initial region is connected, so we can find a path between two of the images of $D$ which does not hit the third image of $D$. Some image of this path under $G$ is a path from $D$ to $gD$, which only intersects the discs in its endpoints, and is disjoint from all the sweepout surfaces. Label this path $\g$. Then $\orbit \gamma$ is a simple closed curve that connects the three discs in one of two ways, illustrated in Figure \ref{picture146} below. In the non-transverse case, the region containing $\orbit \gamma$ is still an initial region after the saddle move, so we may assume we are in the transverse case. 

\begin{figure}[H]
\begin{center}
\epsfig{file=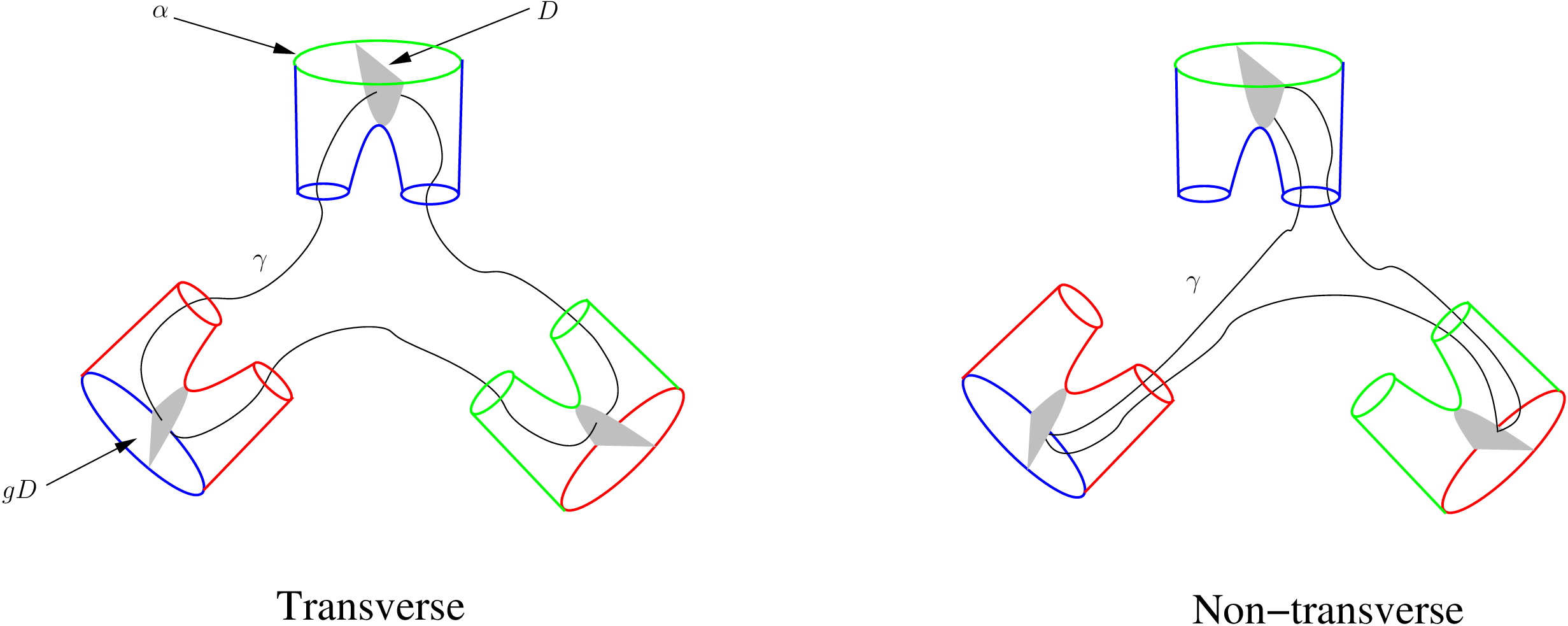, height=150pt}
\end{center}
\caption{The invariant curve $\orbit \g$}
\label{picture146}
\end{figure}

The saddle arc in the red surface meeting the green double curve $\alpha$ is either essential or inessential in the complement of $\alpha$.

First, suppose the saddle arc is inessential, so the saddle arc, union some green arc of $\alpha$, bounds a disc $D'$ in the red surface. This is illustrated below in Figure \ref{picture173}.

\begin{figure}[H]
\begin{center}
\epsfig{file=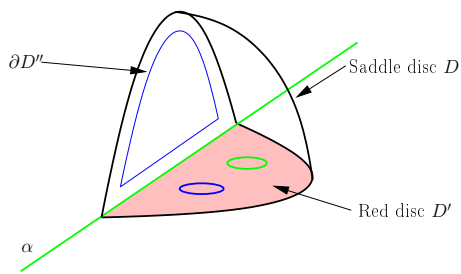, height=100pt}
\end{center}
\caption{The saddle arc is inessential.}
\label{picture173}
\end{figure}

If the saddle arc in the red surface which meets the blue double curve is contained in the red disc $D'$, then this saddle arc is also inessential, so we may choose to work with this one instead, so we may assume we have chosen an inessential saddle arc that bounds an innermost disc.

The union of $D$ and $D'$ forms a disc whose boundary lies in one of the blue surfaces. Let $D''$ be a disc parallel to $D \cup D'$. We may choose this disc so that it is disjoint from the double curves, and from $\orbit D$, so it intersects the sweepout surfaces in simple closed curves only. Cut moves do not disconnect the initial region, by Case 2, so we may remove the intersections of the sweepout surfaces with $D''$. We may need to change the path $\gamma$ if it passes through a cut disc, however we may always choose a new path from $D$ to $gD$, which we shall also call $\gamma$, as cut moves do not disconnect the initial region. The boundary of $D \cup D'$ is a simple closed curve in the blue surface which bounds a disc $D''$ disjoint from any other sweepout surface. If the curve is essential in the blue surface then do a cut move along $D''$. This turns the blue surface into a union of spheres, and does not disconnect the initial region, by the previous case. So we may assume the boundary of $D \cup D'$ bounds a blue disc with no double curves in its interior. However, the union of this blue disc with the red disc $D'$, and the saddle disc $D$, is a $2$-sphere, which the curve $\orbit \g$ intersects precisely once transversely, a contradiction.

If the saddle arc is essential, then the green double curve it
intersects must be inessential in the red torus. This is because if
the green double curve were essential, then the essential saddle arc
in the red torus would go from one side of the green essential curve
to the other, and this can only happen if there is exactly one
essential green curve on the red torus. This is a contradiction, as
the union of all the green curves bounds the subsurface of the red
torus lying inside the blue surfaces, and hence there must be an even
number of essential green curves.

If either of the saddle arcs in the red surface is inessential, we can apply the argument above, so we may assume that both saddle arcs are essential, and intersect inessential double curves. We can remove all the other inessential double curves without disconnecting the initial region. If this cuts the torus into spheres, then the saddle arcs are now inessential, and we can apply the arguments above, or else they are both still essential, and now all double curves in the torus are essential, and the action is standard by Lemma \ref{lemma:parallel}.
\end{case}

This completes the proof of Lemma \ref{lemma:triple}.
\end{Proof}

\section{Disjoint bigons for non-triple point moves.} \label{section:disjoint bigons}

The aim of this section is to prove the following lemma:

\begin{lemma} \label{lemma:main}
If the configuration contains triple points, then for any non-triple point move the configuration contains a bigon disjoint from the move.
\end{lemma}

On a sphere, any pair of intersecting curves creates at least four disjoint bigons. On a torus, a pair of essential curves may intersect without bigons. Even a pair of curves which are both inessential may intersect in only two bigons.

\begin{figure}[H]
\begin{center}
\epsfig{file=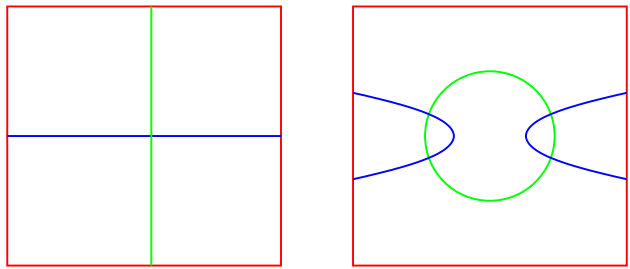, height=100pt}
\end{center}
\caption{Not many bigons}
\end{figure}

The aim of this section is to prove Lemma \ref{lemma:main}, and show that in spite of the fact that there are ``fewer'' bigons on a torus, every non-triple point move has a disjoint bigon. 

If there are any triple points on the sphere components of the sweepout surfaces, then there are disjoint bigon-orbits by \cite[Lemma 5.1]{mr}, so we may assume that all triple points lie on the torus components of the sweepout surfaces.
We prove the result in two main steps. The main technical result is to show that there must be at least three bigons. This immediately gives disjoint bigons for all of the non-triple point moves, except saddle moves, by the arguments in cases 1--3 in the proof of Lemma 5.1 from \cite{mr}. The case of saddle moves is handled by applying the fact that there are at least three bigons to the configuration before and after the saddle move, and doing some elementary combinatorics. We now give a brief overview of this section.

The red torus divides the blue and green tori into subsurfaces
properly embedded in the solid tori on either side of the red Heegaard
torus. The green surfaces in one solid torus component of the
complement of the red torus have boundary consisting of blue double
curves, and are images of the subsurfaces of the red torus divided
along the green double curves. Similarly, the blue surfaces have
boundary consisting of green double curves, and are images of the
subsurfaces of the red torus divided along the blue double curves. If
there are many double curves, but few bigons, then there must be many
parallel double curves, which bound annuli containing essential arcs
only, which we shall call {\bf strips}. The double curves with triple
points divide the surface into complementary regions. We shall call
the closures of the connected components of these regions {\bf faces}.
The boundary of a face is composed of double arcs of alternating
colours, and need not be connected.

In Section \ref{section:intersection} we assume that there are only
two bigons, and then show that there are restrictions on the number of
intersections between discs and strips of different colours on a given
side of one of the Heegaard tori. If a green and a blue subsurface
intersect on the same side of the red Heegaard torus, then the green
surface divides the solid torus bounded by the red Heegaard torus into
pieces with boundary consisting of red and green surfaces. The blue
faces are properly embedded in these pieces, with the green double
arcs in the boundary of the blue faces lying in the red surfaces, and
the red double arcs lying in the green surfaces. An Euler count
argument shows that if there are only two bigons then a face can have
at most eight arcs in its boundary, so we can construct the different
ways in which a face can be embedded in a particular piece with green
and red boundary. The faces fit together to form the original blue
surface, and we show that there are only finitely many ways in which
this can happen. In particular, we show that if there are only two
bigons, then a green strip and a blue strip may not intersect on the
same side of the red torus. This puts an upper bound on the number of
parallel double curves that may arise in a configuration with only two
bigons.

In Section \ref{section:three bigons} we show that there are at least
three bigons. If there are at most two bigons, then all simple closed
curves that bound discs of a given colour are parallel, so there are
three different ways in which the green curves divide the red torus up
into subsurfaces.  If there are no essential curves, then the
subsurfaces consist of a punctured torus, together with some
inessential annuli and a single disc. If there are both essential and
inessential curves, then the subsurfaces consist of a pair of pants,
together with some annuli, which may be either essential or
inessential, and a disc. If there are only essential curves, then all
subsurfaces are annuli. The green and blue curves may divide the red
torus up in different ways. Furthermore the results of the Section
\ref{section:intersection} mean that if there are only two bigons,
there may not be more than three parallel double curves, as parallel
double curves create strips that lie on alternating sides of the tori.
We go through each case in turn, and show that configurations with two
bigons may not arise, usually by counting the number of triple points,
which must be divisible by six.

This means that there are disjoint bigons for all non-triple point moves, except for saddle moves. Finally, in Section \ref{section:saddle} we show that there are disjoint bigons for saddle moves. If there is a configuration in which a saddle move intersects all the bigons, then before and after the saddle move there are at least three bigons, none of which is disjoint from the saddle move. This puts strong restrictions on the configurations which may occur, and we are able to eliminate them by elementary combinatorial arguments.

\subsection{Definitions and preliminary combinatorics.}

\begin{definition}{\bf Annuli types.}

Let $A$ be a properly embedded annulus in a solid torus. Then $A$ has precisely one of the following properties.

\begin{enumerate}

\item Both boundary components are inessential, and bound disjoint discs in the boundary of the solid torus.

\item Both boundary components are inessential and parallel.

\item One boundary component is inessential, and the other is essential, in fact a meridian.

\item Both boundary components are essential and parallel.

\end{enumerate}

These four cases are illustrated below.
\end{definition}

\begin{figure}[H]
\begin{center}
\epsfig{file=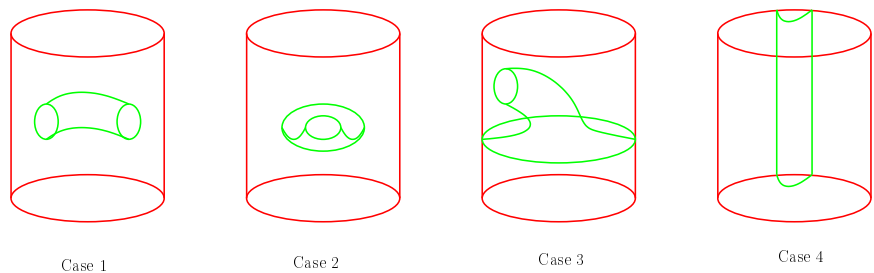, height=100pt}
\end{center}
\caption{Properly embedded annuli in a solid torus.}
\label{picture:annulus types}
\end{figure}

\begin{remark}
In the first three cases, the green annuli may be knotted inside the red torus.
\end{remark}

\begin{definition}{\bf Inner and outer solid tori.}
\label{definition:inner}

An annulus with parallel essential boundary components divides the solid torus into two solid tori. At least one of these has meridional discs that intersect each boundary curve once, call one of these the {\bf inner solid torus}, and the other solid torus the {\bf outer solid torus}.
\end{definition}

\begin{definition}{\bf Faces, squares, hexagons, etc.}

The complementary components of the double curves with triple points in the sweepout surfaces are called {\bf faces}. If a face is a disc, we say it is an $n$-gon if its closure contains $n$ triple points in its boundary. In particular, a $4$-gon is a square, a $6$-gon is a hexagon, and an $8$-gon is an octagon. A face need not have connected boundary.
\end{definition}

A face may contain an arbitrary collection of simple closed double curves without triple points. When we draw pictures of faces we will omit these simple closed curves without triple points.

\begin{definition}{\bf Strips, 1bigons and 2bigons.}

Let $A$ be an annulus bounded by a pair of double curves of the same colour in one of the sweepout tori, which contains no simple closed curves or triple points in its interior. We say that the annulus $A$ is a {\bf strip}, if all double arcs it contains are essential. We include the degenerate case in which there are no double arcs in the annulus. We say $A$ is a {\bf 1bigon} if all faces in $A$ are squares, except for precisely one bigon, and one hexagon. We say $A$ is a {\bf 2bigon} if all faces in $A$ are squares, except for precisely two bigons, and either an octagon, or a pair of hexagons.
\end{definition}

\begin{figure}[H]
\begin{center}
\epsfig{file=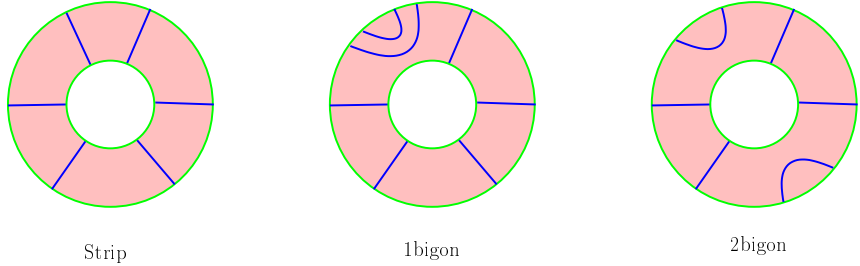, height=100pt}
\end{center}
\caption{Examples of strips, 1bigons and 2bigons.}
\end{figure}

We will make extensive use of the following facts about the configuration. The union of the green curves on the red torus represents a zero homology class in the red torus, as it bounds the subsurface of the red torus which lies on one side of the separating blue surfaces. Similarly, the blue curves represent a zero homology class in the red torus. Double curves come in multiples of three, one of each colour. In particular, in each blue-green diagram the set of green curves gets mapped to the set of blue curves, so there are the same number of blue and green curves, and furthermore, the number of green curves with $n$ triple points is the same as the number of blue curves with $n$ triple points, for each $n$.

Triple points come in multiples of six, so the total number of triple points must be divisible by six. Triple points also come with an orientation, which is preserved by $G$, so in particular, each triple point has three images under $G$, none of which may be adjacent.

\begin{lemma}\label{lemma:reflection}
Let $D$ be a green disc with blue boundary, which bounds a disc $D'$ in the red torus, so that $D \cup D'$ is a $2$-sphere. Furthermore, suppose that all double arcs in $D$ are parallel, and all double arcs in $D'$ are parallel. Then all the blue faces in the $3$-ball bounded by $D$ and $D'$ contain the same number of triple points in each boundary component.
\end{lemma}

\begin{Proof}
If we imagine the triple points on the boundary of the disc as equally spaced points on the unit circle, then the parallel arc identifications act as a reflection. The length of the orbits of the points under these maps corresponds to the number of arcs in each boundary curve in the union of the two inessential discs. If the reflections are the same, we get bigon boundaries. Otherwise, two reflections give a rotation, and so in fact all the orbits have the same length.
\end{Proof}

\begin{lemma} \label{lemma:at most octagon}
If there are two bigons, then all other faces are squares, except for either an octagon, a pair of hexagons, or an annulus with two vertices in each boundary component.
\end{lemma}

\begin{Proof}
The Euler characteristic of a torus is zero. As every vertex is four-valent, an $n$-gon contributes $1 - \frac{n}{2} + \frac{n}{4} = 1 - \frac{n}{4}$ to the Euler characteristic, so bigons contribute $+\half$, squares contribute $0$, hexagons contribute $-\half$, etc. An annular face with $n$ and $m$ triple points in each boundary component contributes $-\frac{n+m}{4}$, and a punctured torus face with $n$ vertices in its boundary contributes $-1-\frac{n}{4}$. Two bigons contribute $+1$ to the Euler characteristic, so the other faces sum to $-1$. The only ways in which this can happen are if all other faces are squares, except for either two hexagons, or a single octagon, or an annulus formed by tubing two bigons together.
\end{Proof}

If the boundary of an inner solid torus consists of a pair of strips,
then any faces which are discs in the interior of the solid torus are
bigons with essential boundary. We will use this observation
repeatedly, so we state it as a lemma.

\begin{lemma} \label{lemma:meridional bigons}
Suppose an inner solid torus has boundary which is the union of two
strips. Then a face which is a disc in the interior of the inner solid
torus is a bigon with essential boundary. 
\end{lemma}

\begin{Proof}
A choice of orientation for the disc in the interior of the inner
solid torus induces orientations in its boundary arcs in each of the
strips. The double arcs in each strip are all essential and parallel,
so the induced orientation is the same in each arc in the boundary of
the strip. This means that the geometric intersection number of the
boundary of the disc with a boundary curve of one of the strips is the same as
the algebraic intersection number. As an inner solid torus has
meridional discs that hit each boundary curve of a strip once, this
means that the intersection number must be one, and the boundary of
the face which is a disc is a bigon with essential boundary.
\end{Proof}

We now show that there must always be at least two bigons.

\begin{lemma}
A configuration with triple points contains at least two bigons. \label{lemma:at least two bigons}
\end{lemma}

\begin{Proof}
Any inessential curve with triple points creates at least two bigons, so we assume all curves with triple points are essential. Any pair of intersecting curves with non-minimal intersection creates at least two bigons, so we may assume all curves are essential, with different slopes, and have minimal intersection. So the green curves divide the red torus up into strips, as do the blue curves, and every face is a square. 

Consider an outermost green strip with blue boundary on one side of
the red torus, which divides the solid torus on one side of the red
torus into inner and outer solid tori. The boundary of the inner solid
torus consists of two strips, so by Lemma \ref{lemma:meridional
 bigons} any blue faces which are discs inside the inner solid torus
must be bigons, a contradiction.
\end{Proof}

So from now on we may assume that there are at least two bigon-orbits. Next we show that every configuration has at least three bigon-orbits. A non-triple point move which is not a saddle move may involve at most two bigon-orbits, so this deals with all cases, except for saddle moves, which we consider separately, in Section \ref{section:saddle}.

\subsection{Intersection lemmas.} \label{section:intersection}

We now prove some useful lemmas that show that if certain kinds of discs or strips intersect on the same side of the red torus, then there are bounds on the number of intersections that they may have. The results of this section are summarised in the list below.

If there are only two bigons, then:

\begin{itemize}

\item[\ref{lemma:intersecting discs}] Two inessential discs may intersect in at most two arcs.

\item[\ref{lemma:disc intersect strip}] A disc may intersect a strip only in one of the following three ways.

\begin{enumerate}

\item The strip has parallel inessential boundary components, the disc is essential, and there are at most two arcs of intersection.

\item The strip has parallel essential boundary, the disc is essential, and there are at most two arcs of intersection.

\item The strip has parallel essential boundary, the disc is inessential, and there are at most four arcs of intersection.

\end{enumerate}

\item[\ref{lemma:no strips}] Two strips may not intersect.

\end{itemize}

We now prove these lemmas.

\begin{lemma} \label{lemma:intersecting discs}
Suppose there are only two bigons, and a pair of inessential discs of different colours intersect on the same side of the red torus; then they have at most two arcs of intersection.
\end{lemma}

\begin{Proof}
Assume there are at most two bigons, and consider the green and blue discs only, ignoring all the other green and blue surfaces. As there are only two bigons,  the inessential green disc contains parallel red arcs, and its boundary also bounds a red disc in the red torus, which contains parallel green arcs. Again, as there are only two bigons, the blue disc contains parallel red arcs, so the blue faces in the blue disc which are not bigons are squares. These two discs bound a $3$-ball, which we shall call the inner $3$-ball, and we will call its complement in the solid torus bounded by the red torus, the outer solid torus. Then Lemma \ref{lemma:reflection}, the reflection lemma, implies that all of the blue faces inside the inner $3$-ball have the same boundary length. As these blue faces all lie inside the blue inessential disc, they are all either bigons or squares.

If the blue faces are bigons, then there may be at most two of them, so there may be at most two arcs of intersection. If the blue faces are squares, then the two bigons lie in the outer solid torus. In order to avoid creating red bigons outside the red disc, the boundary arcs of the blue bigons lying outside the red disc must be essential, as illustrated in Figure \ref{picture13} below. 

\begin{figure}[H]
\begin{center}
\epsfig{file=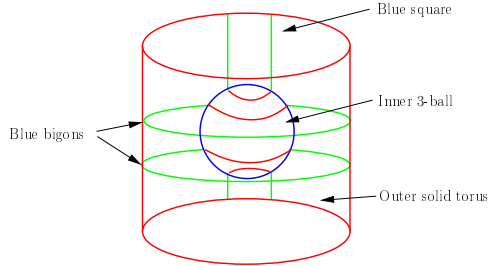, height=100pt}
\end{center}
\caption{Bigons in the outer solid torus.}
\label{picture13}
\end{figure}

The only possible blue squares lying in the outer solid torus are ones parallel to the blue square shown in Figure \ref{picture13} above. However, as there are blue squares inside the inner $3$-ball, then any blue square in the outer solid torus has both red arcs attached to the same blue square in the inner $3$-ball, forming a blue strip with two essential red arcs. Therefore there may be no blue squares in the outer solid torus, so there may be at most two arcs of intersection.
\end{Proof}

\begin{lemma} \label{lemma:disc intersect strip}
Suppose there are exactly two bigons, and a green strip and a blue disc have intersecting boundaries on the same side of the red torus. Then one of the following cases occurs.

\begin{enumerate}

\item The strip has parallel inessential boundary components, the disc is essential, and there are at most two arcs of intersection.

\item The strip has parallel essential boundary, the disc is essential, and there are at most two arcs of intersection.

\item The strip has parallel essential boundary, the disc is inessential, and there are at most four arcs of intersection.

\end{enumerate}
\end{lemma}

\begin{Proof}
Consider just the disc and the annulus, and ignore all other surfaces inside the solid torus bounded by the red torus. There are four different ways in which the boundary of the annulus may be embedded in the boundary of the solid torus, as illustrated in Figure \ref{picture:annulus types}. We consider each one in turn.

\setcounter{case}{0}

\begin{case}{The strip has inessential boundary components which bound disjoint discs.}

As there are exactly two bigons, the strip cannot have inessential boundary components that bound disjoint discs.
\end{case}

\begin{case}{The strip has parallel inessential boundary components.}

If the green strip has parallel inessential boundary components, then there are at least two bigons in the innermost disc bounded by the blue double curves, so there cannot be any bigons in the red annulus between the two blue double curves, so this annulus is a strip. The union of the green and red strips is a torus, which is contained in the solid torus bounded by the red torus. This torus has compressing discs on both sides, namely the blue faces inside the inner solid torus, and the red subdisc of the red torus bounded by the innermost blue boundary component of the green strip. The solid torus bounded by the two strips is therefore an unknotted solid torus contained in the larger solid torus bounded by the red torus, so the boundaries of the red disc and the blue compressing discs have intersection number one. This means that the blue faces cannot be squares, and so are bigons. But there may be at most two bigons, so the disc and the strip may have at most two arcs of intersection.
\end{case}

\begin{case}{The strip has exactly one essential boundary component.}

There are blue bigons with one boundary arc essential in the the green strip, and whose other boundary is therefore essential in the red pants bounded by the blue double curves, so the strip is unknotted. Therefore the strip divides the solid torus bounded by the red torus into another solid torus, whose boundary contains two green annuli, one essential and one inessential, as illustrated in Figure \ref{picture192} below.

\begin{figure}[H]
\begin{center}
\epsfig{file=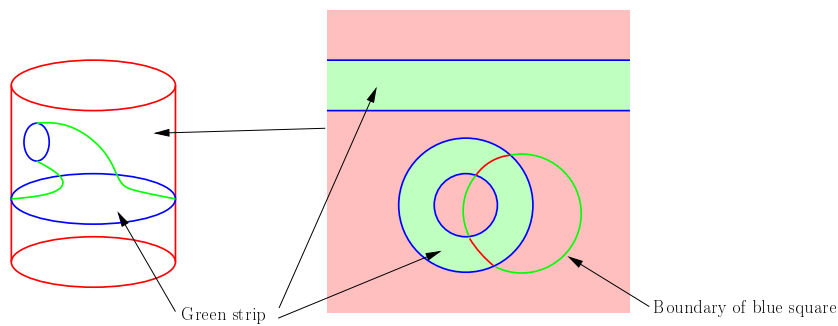, height=120pt}
\end{center}
\caption{The green strip divides the red solid torus into another solid torus.}
\label{picture192}
\end{figure}

There is a blue face with a boundary arc inside the red disc. This blue face has at least two red boundary arcs in the green strip adjacent to the red disc, so it is a square. This blue square is either a meridian disc or an inessential disc in the solid torus bounded by the red pants and the green strip. As the blue square is disjoint from the essential green strip, the square cannot be a meridian disc, so it must be an inessential disc, as illustrated in Figure \ref{picture192} above. This creates an inessential arc in the red pants, creating at least three bigons, a contradiction.
\end{case}

\begin{case}{The strip has both boundary components essential.}

The strip divides the solid torus on one side of the red torus into an inner solid torus and an outer solid torus, as in Definition \ref{definition:inner}. If there are more than two arcs, then there are squares on both sides of the strip. All squares in the inner solid torus have boundary that consists of two essential arcs in the strip, together with two inessential arcs connecting them inside the annulus in the red torus. In particular, there are two bigons in this annulus, so there may be no bigons in the red annulus in the boundary of the outer solid torus, so this red annulus is a strip. The outer solid torus therefore has boundary consisting of two strips, so the faces in the outer solid torus are either all bigons or all squares. If they are bigons, there may be at most two of them, so there are at most two arcs of intersection, so we may assume that they are all squares, and the inner solid torus contains both bigons. All the squares on the outside are parallel to the one shown in the left hand side of Figure \ref{picture14} below. The top of the red cylinder should be identified with the bottom by a half twist so that there is a single green annulus.

\begin{figure}[H]
\begin{center}
\epsfig{file=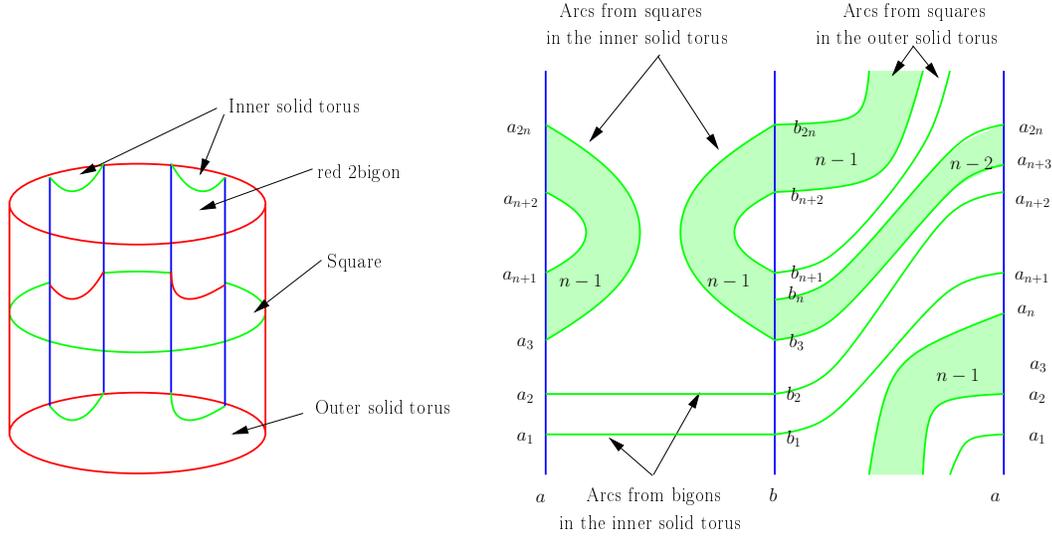, height=200pt}
\end{center}
\caption{Counting triple points in double curves.}
\label{picture14}
\end{figure}

The boundary of the inner solid torus consists of two annuli: a green
strip containing red arcs, and a red 2bigon containing green arcs. We
now explain how this forces a particular arrangement of blue faces
inside the inner solid torus. The inner solid torus contains two
bigons, and the boundary of each bigon consists of a green arc and a
red arc. The red arc lies in the green strip, and is therefore an
essential arc, so the green arc in the boundary of the bigon must be
essential in the red annulus which is a 2bigon. This means that each
bigon in the inner solid torus is a meridional disc for the inner
solid torus. If there are any blue squares in the inner solid torus,
then their red boundary arcs are essential arcs in the green strip.
The blue squares lie in the complement of the bigons, which are
meridional discs for the inner solid tori, so they have inessential
boundaries in the boundary of the inner solid torus, so as the red
arcs in the boundary of each blue square are essential, the green arcs
must be inessential. As the red annulus is a 2bigon, there are at most
two families of inessential green arcs with endpoints on each boundary
component of the red 2bigon, so in fact all blue squares must be
parallel inside the inner solid torus, and their boundaries consist of
a pair of essential red arcs, and a pair of inessential green arcs,
with one green arc meeting each blue double curve boundary component
of the red 2bigon. As the blue squares are all parallel, the two blue
bigons are also adjacent inside the inner solid torus.

The strips forming the boundaries of the inner and outer solid tori
share a common pair of blue double curves, label one of these double
curves $a$ and the other one $b$. The double curve $a$ contains a pair
of adjacent triple points which lie in the blue bigons contained in
the inner solid torus, label these triple points $a_1$ and $a_2$. Now
label the remaining triple points $a_3, \ldots, a_{2n}$ following the
circular order coming from the blue double curve $a$. We may assume
there are $2n$ triple points on each of the blue double curves as
there are an even number of triple points on each double curve. Label
the triple points in the blue double curve $b$ in a similar manner. To
be precise, label the triple points which lie in the boundaries of the
blue bigons in the inner solid torus $b_1$ and $b_2$, so that $a_1$
and $b_1$ lie in the boundary of a common blue bigon. Now label the
remaining triple points $b_3, \ldots b_{2n}$, following the circular
ordering coming from the blue double curve $b$. This means that the
green arcs contained in the boundaries of squares in the inner solid
torus connect $a_k$ to $a_{2n+3-k}$ and $b_k$ to $b_{2n+3-k}$, for $3
\leqslant k \leqslant n$. Also, the green arcs contained in the
boundaries of squares in the outer solid torus connect $a_k$ to
$b_{k+n}$, taking indices modulo $2n$. This is illustrated on the
right hand side of Figure \ref{picture14} above. Opposite sides of the
right hand diagram are identified to form the red torus, in particular
there is only one blue double curve $a$ in the red torus, as its two
images in Figure \ref{picture14} are identified. The green shaded
regions denote some number of parallel arcs; the number inside the
region indicates exactly how many.

There is a simple closed green double curve containing the eight triple points $ \{ a_1$, $b_1$, $a_{n+1}$, $a_{n+2}$, $b_2$, $a_2$, $b_{n+2}$, $b_{n+1} \}$. However, there must be exactly one double curve, so in fact this implies that this must be all the triple points, so in fact $n=2$, and there are four arcs of intersection between the disc and the strip. However, in this case the disc is inessential, so if the disc is essential, the number of arcs of intersection is at most two. 
\end{case}

This completes the proof of Lemma \ref{lemma:disc intersect strip}.
\end{Proof}

We wish to prove the following lemma:

\begin{lemma} \label{lemma:no strips}
Suppose there are only two bigons. Then a green strip and a blue strip may not intersect on the same side of the red torus.
\end{lemma}

We break the proof up into two parts. First we show Lemma \ref{lemma:no strips unless parallel}, that if there are two bigons, then a pair of intersecting strips must have parallel and essential boundaries, with exactly four arcs of intersection. We then complete the proof by showing in Lemma \ref{lemma:no strips parallel} that, if there are only two bigons, the configuration can not have parallel essential curves.

\begin{lemma} \label{lemma:no strips unless parallel}
Suppose there are only two bigons, and a green strip and a blue strip intersect on the same side of the red torus. Then the boundary components of the strips are essential, with the same slope, and the strips have exactly four arcs of intersection. In particular, this means that the outer solid tori for the strips have meridional discs that intersect each boundary curve of the annuli twice.
\end{lemma}

\begin{Proof}
The green annulus may be one of the four types shown in Figure \ref{picture:annulus types}. We deal with each case in turn.

\setcounter{case}{0}

\begin{case}{The green strip has inessential boundary components which bound disjoint discs.}

In this case, the two blue curves bound disjoint discs. Each curve contains triple points, so each disc contains at least two disjoint bigons, so there are at least four bigons. 
\end{case}

\begin{case}{The green strip has inessential boundary components which bound nested discs.}
The innermost blue curve contains two bigons, so if there are only two
bigons, all green arcs in the annulus bounded by the two blue double
curves must be essential, and so the blue double curves bound a red
strip in the red torus. Let $T$ be the torus formed from the union of
the green and red strips. Any square in the region bounded by $T$ has
intersection number two with each blue curve. We can compress the
torus $T$ along such a blue square to form a $2$-sphere, which must
bound a $3$-ball by irreducibility, so in fact the torus $T$ bounds a
solid torus, with a blue square as a meridional disc. However, the
torus $T$ also has a compressing disc on the other side, given by the
disc in the red torus with boundary the innermost blue boundary
component of the green strip. The boundaries of two compressing discs
have intersection number two, giving an $\R \mathbb{P}^3$ summand to
the solid torus bounded by the red torus, a contradiction.
\end{case}

\begin{case}{The green strip has exactly one essential boundary component.}

The green annulus divides the solid torus into a $3$-dimensional region, which we shall call $R$, which has torus boundary. The region $R$ need not be a solid torus, as the green annulus may be knotted. The boundary of $R$ consists of a red disc and a red pair of pants, and two copies of the green annulus, as illustrated in Figure \ref{picture192}. Consider the blue double curve $\g$ which is inessential in the red torus. This appears twice in $\d R$, once as the boundary of the red disc, which we shall call $\g_1$, and also as a boundary component of the red pair of pants, which we shall call $\g_2$. We now show that $\g_2$ is essential in $R$. We can make a simple closed curve in $\d R$ by taking one of the essential red arcs in the green strip, and connecting its end points by an arc in the red pants. If the blue double curve bounds a disc $D$ in $R$, then the union of $D$ and the red disc from a $2$-sphere inside the solid torus bounded by the red torus, which this curve intersects once, a contradiction.

The inessential blue curve $\g$ has triple points, so there are at least two disjoint bigons inside the red disc that it bounds. Consider a blue square in $R$ which contains a green boundary arc lying in the red disc. The two red arcs in the boundary of the square lie in the green strip, and the other green arc lies in the red pair of pants, with both endpoints in the same boundary component. If the green arc is inessential in the red pair of pants, then it creates an extra bigon, so the green arc must be essential in the red pair of pants. However, this means that the boundary of the square is now homotopic inside $R$ to the essential blue double curve, a contradiction.
\end{case}

If any boundary component of the blue annulus is inessential, then we are in one of the cases above, using the blue annulus instead of the green annulus, so we may now assume that all boundary curves are essential.

\begin{case}{The green strip has both boundary components essential.}

If the blue strip has an inessential boundary component, then we can swap the colours green and blue, and use one of the cases above, so we may assume that both strips have both boundary components essential.

The green strip divides the red solid torus into inner and outer solid tori, as in Definition \ref{definition:inner}. All the blue squares in the inner solid torus have boundary consisting of two essential red arcs in the green strip, and two inessential green arcs in the red annulus, and they are all parallel. In particular, the green curves do not cross from one boundary component to the other in the red annulus, so the boundary components of the blue strip have the same slope as the boundary components of the green strip.

If there are only two bigons, there may be no bigons in the other red annulus, so all green curves cross from one boundary component to the other, and so must be parallel. So the meridian discs of the outer solid torus intersect each essential curve twice, as illustrated in Figure \ref{picture14}.

Label one of the blue double curves $a$ and the other one $b$. We may
assume there are $2n$ triple points on each curve, as there are an
even number of triple points on each double curve. On the curve $a$,
we will label the triple points $a_1, \ldots a_{2n}$, following the
circular ordering, so that the triple points $a_1$ and $a_{2n}$ lie in
the outermost blue square in the inner solid torus. We will label the
triple points on $b$ with labels $b_1, \ldots b_{2n}$, so that $b_1$ and $b_{2n}$
lie in the same blue square in the inner solid torus as $a_1$ and
$a_{2n}$, and so that the orientations induced on the double curves by
the circular orderings are the same. All squares in the inner and
outer solid tori are parallel, so this determines the pattern of green
curves in the red torus, which is illustrated below on the left hand
side of Figure \ref{picture3}. The shaded regions indicate that there
are possibly many parallel curves, and the opposite sides of the
square are identified to form a torus.

Each green double curve contains four triple points, $\{ a_k, a_{2n+1-k}, b_{n+1-k}, b_{n+k} \}$. As there are two green double curves, this means $n=2$, so there are four arcs of intersection between the strips. This is illustrated on the right hand side of Figure \ref{picture3}.

\begin{figure}[H]
\begin{center}
\epsfig{file=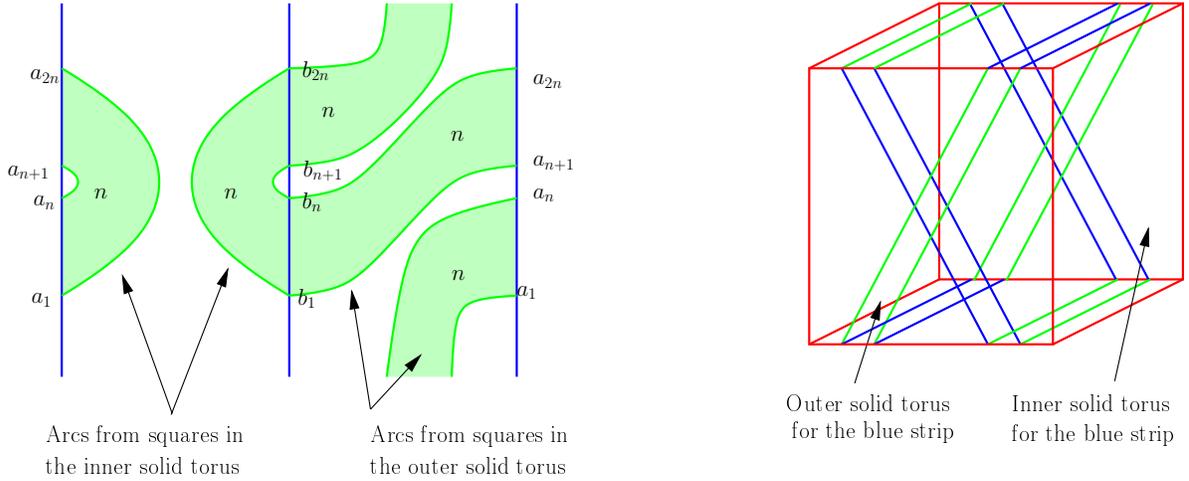, height=180pt}
\end{center}
\caption{A green strip and a blue strip which intersect.} 
\label{picture3}
\end{figure}

In the right hand side of Figure \ref{picture3} above, the solid torus is drawn as a cube with top and bottom face identified with a half twist, so that there is a single blue annulus, and a single green annulus.
\end{case}

This completes the proof of Lemma \ref{lemma:no strips unless parallel}.
\end{Proof}

We now eliminate the special case from Lemma \ref{lemma:no strips unless parallel}.

\begin{lemma} \label{lemma:no strips parallel}
If there are only two bigons, a configuration cannot have intersecting essential double curves with the same slope. In particular, this implies there are no intersecting strips.
\end{lemma}

\begin{Proof}
We consider two cases, depending on whether or not there are inessential curves.
\setcounter{case}{0}

\begin{case}{No inessential curves.}

By Lemma \ref{lemma:no strips unless parallel}, Case 4, there is an
essential green double curve and an essential blue double curve which
have geometric intersection number two. So if there are no inessential double
curves, then all essential double curves are parallel to
the ones shown below in Figure \ref{picture41}.

\begin{figure}[H]
\begin{center}
\epsfig{file=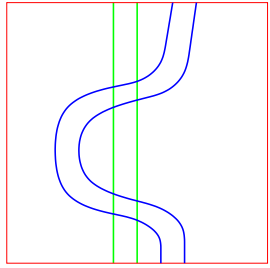, height=100pt}
\end{center}
\caption{Two bigons, essential curves with the same slope.}
\label{picture41}
\end{figure}

There are the same number, $n$ say, of blue and green curves with triple points. The number of triple points is then $2n^2$. The number of triple points is divisible by six, so $n$ is divisible by $3$. If $n$ is at least $3$, then there is a set of three consecutive parallel green double curves, which bound a pair of adjacent red strips with green boundary in the red torus, so the images of these under $g$ give rise to green strips with blue boundary, which lie on opposite sides of the red torus. Similarly, as $n$ is at least $3$, there is a set of three consecutive parallel blue double curves, which create a pair of adjacent red strips with blue boundary in the red torus, and images of these under $g^2$ give rise to a pair of blue strips with green boundary, which lie on opposite sides of the red torus. Therefore there are green and blue strips on both sides of the red torus, and so each side of the red torus contains an intersecting pair of green and blue strips, as in Lemma \ref{lemma:no strips unless parallel}, Case 4, as illustrated in Figure \ref{picture3}.

Consider a green annulus with blue boundary curves which contains
inessential red arcs. A blue boundary component of this annulus has
triple points, and is parallel to the other blue double curves, and so
intersects a blue strip with green boundary curves, and red essential
arcs. Now consider a red arc in the green annulus which bounds an
innermost disc in the green annulus. This green disc is a bigon with a
red double arc which is essential in the blue strip, so the other blue
double arc in the boundary must also be essential in the red annulus
it is contained in. The outer solid torus for the blue strip has
meridional discs that intersect each boundary component of the annuli
twice, by Lemma \ref{lemma:no strips unless parallel}, so the bigon
cannot lie in the outer solid torus, as it has essential boundary
arcs, and would therefore be a meridional disc intersecting the
boundary components of the annuli once.

The green bigon must therefore lie in the inner solid torus for the blue strip. As the red boundary arc of the bigon is essential in the blue strip, the blue boundary arc of the bigon runs from one boundary component of the red annulus to the other, as illustrated in Figure \ref{picture162} below.

\begin{figure}[H]
\begin{center}
\epsfig{file=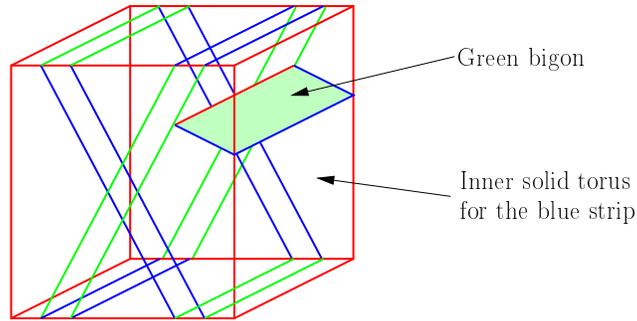, height=120pt}
\end{center}
\caption{A green bigon in the inner solid torus for the blue strip.}
\label{picture162}
\end{figure}

The boundary of the blue strip is a pair of green double curves, which divide the red torus into a pair of annuli, one lying in the inner torus for the blue strip, and the other lying in the outer solid torus. The red annulus in the outer solid torus contains essential blue arcs coming from the squares which are meridional discs for the outer solid torus. The green bigon in the inner solid torus creates essential blue arcs in the other red annulus as well. This is a contradiction, as Figure \ref{picture41} shows there may only be essential arcs in one of the red annuli.
\end{case}

\begin{case}{Inessential curves.}
First suppose there are inessential curves of only one colour.
Without loss of generality we may assume that the inessential double
curves are blue. As we are assuming there are two bigons, all
inessential blue curves are parallel in the red torus, and divide the
red torus into an innermost red disc containing parallel green arcs, a
collection of red annuli containing parallel green double arcs, and a
single punctured torus component.

By Lemma \ref{lemma:no strips unless parallel}, Case 4, there is a
green essential curve and a blue essential curve with the same slope,
and with with geometric intersection number two. Such a pair of curves
divide the red torus into three regions, two of which must contain an
innermost bigon. Therefore one of the innermost bigons inside the
innermost blue inessential curve must lie in one of these two regions.
This means that all of the double curves with triple points in the
configuration must be parallel to the curves shown below in Figure
\ref{picture161}.  There is an even number of green curves in total,
and an even number of essential blue curves, so there is an even
number of inessential blue curves. In particular, this means that
there are at least four double curves of each colour, and the red
torus contains red strips with both green and blue boundary. We will
use the labels in Figure \ref{picture161} to refer to how many curves
there are parallel to the one shown.

\begin{figure}[H]
\begin{center}
\epsfig{file=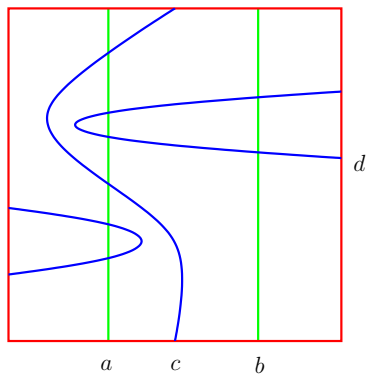, height=140pt}
\end{center}
\caption{Inessential curves of one colour.}
\label{picture161}
\end{figure}

We now show that there is a green strip with blue boundary components, such that at least one of these blue boundary components is inessential in the red torus. The red torus divides the green torus into a number of subsurfaces along the blue double curves. The green surfaces consist of strips, and a pair of 1bigons. If none of the blue boundary components of the green strips are inessential in the red torus then the inessential curves must all be part of the boundary of the two green 1bigons. In this case, the two green 1bigons share a common pair of blue boundary components, creating a green torus with no strips with blue boundary, a contradiction, as the green torus is an image of the red torus under $g$, and the red torus contains strips with both colour boundary components.

Figure \ref{picture161} shows that an inessential blue curve meets all the essential green curves, and by Lemma \ref{lemma:no strips unless parallel}, a green strip with an inessential blue boundary component may not intersect any blue strips, therefore there can't be blue strips on both sides of the red torus. This means the total number of blue double curves is at most four, so must be exactly four, with two essential blue curves, and two inessential blue curves, so $c=d=2$. The green curves parallel to $b$ have four triple points, so there are blue curves with four triple points. The blue curves parallel to $d$ have more than four triple points, so the curves parallel to $c$ have four triple points, so $a=2$, and hence $b=2$. But this implies there are $32$ triple points in total, which is not divisible by six, a contradiction.

So we may now assume that there are inessential double curves of both colours, for example as illustrated in Figure \ref{picture15} below.

\begin{figure}[H]
\begin{center}
\epsfig{file=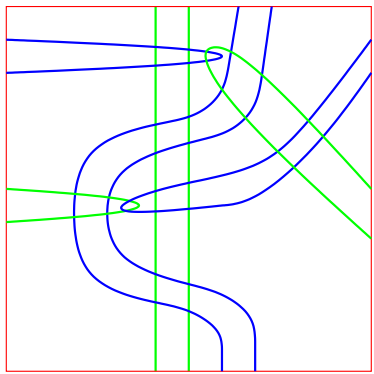, height=140pt}
\end{center}
\caption{Inessential curves of both colours.}
\label{picture15}
\end{figure}

The number of triple points on the inessential curves is strictly greater than the number of triple points on the essential curves, as each essential curve intersects essential curves of the other colour twice, but inessential curves four times. This means that essential curves are mapped to essential curves, and inessential curves are mapped to inessential curves. In particular this implies that the red disc with boundary consisting of the innermost green inessential curve gets mapped by $g$ to a green disc with a blue boundary curve which is inessential in the red torus. Similarly the red disc with boundary consisting of the innermost blue inessential curve gets mapped by $g^2$ to a blue disc with green inessential boundary.

There are an even number of essential curves of the same colour,
creating a strip with essential boundary components. Suppose there are
three or more inessential curves of the same colour. This creates at
least two strips in the red torus with inessential boundary, so there
are intersecting strips on at least one side of the red torus, but one
of these strips has inessential boundary components, contradicting
Lemma \ref{lemma:no strips unless parallel}. So there are at most two
inessential curves. We now consider the cases when there are one or
two inessential curves of each colour.

\begin{subcase}{Two inessential curves of each colour.}
\end{subcase}

\begin{notation}
The green curves divide the red torus up into subsurfaces, which lie on different sides of the blue surfaces. We will write $\{ S_1, S_2, \ldots | T_1, T_2, \ldots \}$ to mean that the curves divide the red torus up into surfaces $S_i$ and $T_i$, with the $S_i$ lying on one side, and the $T_i$ lying on the other side.
In practice it will be convenient to label the surface by its
homeomorphism type, i.e. disc, strip, etc. So $\{$ disc, strip $|$
pants $\}$ means there is a disc and a strip on one side of the blue torus, and
a pair of pants on the other side.
\end{notation}

Two inessential curves create a strip with inessential boundary. If there are four or more essential curves, then there are strips on both sides, one of which intersects a strip with an inessential boundary component, contradicting Lemma \ref{lemma:no strips unless parallel}, so there are two essential curves also. So both the green and blue double curves divide the torus into $\{$ disc, pants $|$ strip, strip $\}$.

A pair of green and blue inessential double curves must have at least four points of intersection, however, they may have arbitrarily many. Blue and green inessential double curves with minimal intersections are illustrated in Figures \ref{picture30} and \ref{picture77}. Extra points of intersection may occur, for example by doing a Dehn twist to the blue double curves along a curve parallel to the essential curves. Figure \ref{picture15} is an example of inessential blue and green curves having non-minimal intersection.

If the green and blue discs are on the same side of the red torus, then they may have at most two arcs of intersection, by Lemma \ref{lemma:intersecting discs}, implying the inessential curves have minimal intersection. If the discs are on opposite sides, then the two strips intersect an inessential disc and a pants. However, an inessential disc may intersect a strip in at most four arcs, by Lemma \ref{lemma:disc intersect strip} which again implies that the inessential curves have minimal intersections. If the inessential curves have minimal intersection, then there are $56$ triple points, as shown in Figure \ref{picture30}. This is a contradiction, as $56$ is not a multiple of six.

\begin{figure}[H]
\begin{center}
\epsfig{file=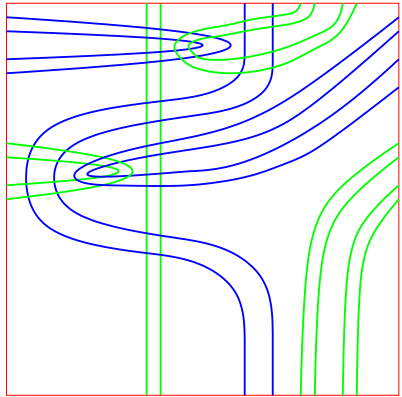, height=150pt}
\end{center}
\caption{Two inessential curves with minimal intersection gives $56$ triple points.}
\label{picture30}
\end{figure}

\begin{subcase}{One inessential curve of each colour.}
\end{subcase}

Suppose there is one inessential curve, with some even number of essential curves. Both the blue and the green curves divide the torus into $\{$ disc, strips $|$ pants, strips $\}$. If the discs are on the same side, so there are at most two arcs of intersection, by Lemma \ref{lemma:intersecting discs}, so the intersection of the inessential curves is minimal, as shown in Figure \ref{picture77} below.

\begin{figure}[H]
\begin{center}
\epsfig{file=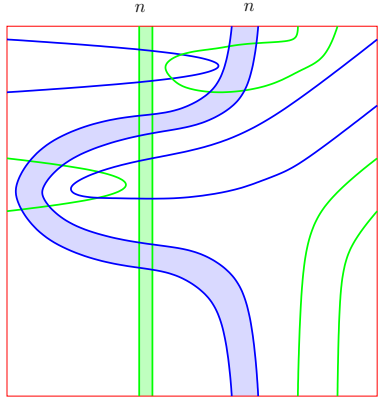, height=150pt}
\end{center}
\caption{One inessential curve.} 
\label{picture77}
\end{figure}

The inessential green curve contains $4+4n$ triple points, and each essential green curve contains $4+2n$ triple points, so in total there are $4 + 8n + 2n^2$ triple points, which it is easy to check is not divisible by six.

Finally, if the discs are on different sides, then there is a green
disc, and at least one green strip, on one side, which both intersect
a blue pair of pants. These green and blue surfaces must intersect as
their boundaries consist of blue and green double curves in the red
torus, and every blue double curve hits every green double curve, as
shown in Figure \ref{picture77}. The inessential curves have more
triple points than the essential curves, so the inessential curves are
images of each other, and the essential curves are images of each
other. This means that the strip has both boundary components
essential. Also, the blue pants does not contain any bigons, as one of
its boundary components is inessential in the blue torus, so by Lemma
\ref{lemma:at most octagon}, the pants contains only squares, together
with either a pair of hexagons or a single octagon. Figure
\ref{picture104} below illustrates the only way of embedding squares,
hexagons and octagons in the inner solid torus on one side of the
green strip. If there are blue squares in the inner solid torus, then
there may not also be any hexagons or octagons in the inner solid
torus, as this would create at least three bigons.

\begin{figure}[H]
\begin{center}
\epsfig{file=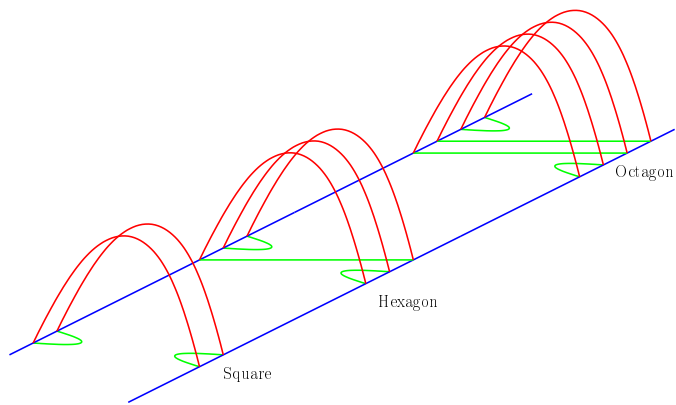, height=150pt}
\end{center}
\caption{Blue faces in the inner solid torus.}
\label{picture104}
\end{figure}

Whichever faces occur, there are two bigons in the red annulus in the boundary of the inner solid torus, so the red annulus forming part of the boundary of the outside solid torus is a strip, so the boundary of the outer solid torus consists of the union of two strips. This means that if a blue disc face which is a $2n$-gon is contained in the outer solid torus, then each blue double curve in the boundary of the two strips intersects the blue disc $n$ times, as all of the arcs in the boundary of the $2n$-gon are essential in each strip. This means that the blue boundary curves of the strips may not be meridional curves for the outer solid torus, and Figure \ref{picture107} below illustrates the possible squares, hexagons and octagons that may occur in the outer solid torus.

\begin{figure}[H]
\begin{center}
\epsfig{file=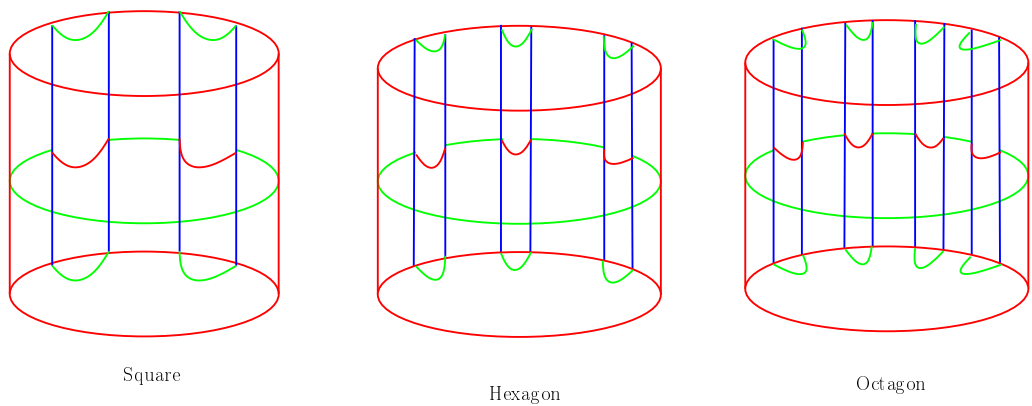, height=150pt}
\end{center}
\caption{Blue faces in the outer solid torus.}
\label{picture107}
\end{figure}

In Figure \ref{picture107} above, the top of the solid torus should be glued to the bottom with an appropriate twist to glue the green surfaces up into a single annulus. The three possibilities shown in Figure \ref{picture107} are mutually exclusive, i.e. there may only be faces of one sort (squares, hexagons or octagons) in the outer solid torus.

If there is a blue strip on the same side of the red torus as the blue pair of pants, then this blue strip gives rise to blue squares on both sides of the green strip, which means there are no hexagons or octagons in either the inner or outer solid tori, a contradiction.

If there are no other blue strips on the same side of the red torus as the blue pair of pants, then there are precisely two essential double curves of each colour, giving a total of three double curves of each colour. The red torus then bounds a solid torus containing a green strip, a green disc and a blue pair of pants, and no other surfaces. There are blue hexagons or octagons on at least one side of the green strip, and that side contains either a single octagon, or at most two hexagons. This means that there may be at most six red arcs in the green strip. However, Figure \ref{picture77} above shows the smallest number of arcs that may occur in the strips, namely eight arcs in the strips, giving a contradiction.
\end{case}

This completes the proof of Lemma \ref{lemma:no strips parallel}, which in turn completes the proof of Lemma \ref{lemma:no strips}.
\end{Proof}

\subsection{There are at least three bigons.} \label{section:three bigons}

In this section, we prove the following lemma.

\begin{lemma} \label{lemma:three bigons}
A configuration with triple points contains at least three bigons.
\end{lemma}

Consider the green double curves which contain triple points. These divide the red torus up into subsurfaces, possibly containing simple closed curves without triple points. We may assume there is at most one disc, as each disc contains at least two bigons. Similarly, the blue double curves with triple points also divide the red torus into subsurfaces. Therefore, the following possible cases may arise: the collection of surfaces may contain a punctured torus component, it may contain a pair of pants, or it may consist of annuli only. We need to consider all pairs of these possibilities, as the green curves and blue curves may divide the red torus in different ways. We deal with the different cases in the order listed in Table \ref{table} below. It doesn't matter which set of curves is green or blue, so we will make an arbitrary choice.

\begin{table}[H]
\begin{center}
\begin{tabular}{llll}
 & Punctured torus & Pants & Annuli \\
Punctured torus & $1$ & $2$ & $3$ \\
Pants && $4$ & $5$ \\
Annuli &&& $6$ \\
\end{tabular}
\end{center}
\caption{Cases} \label{table}
\end{table}

\begin{notation}
We will write $\geom{a}{b}$ to denote the geometric intersection number of two simple closed curves $a$ and $b$. 
\end{notation}

\setcounter{case}{0}

\begin{case}{The green and the blue double curves with triple points divide the red torus into subsurfaces which both have punctured torus components.}

An innermost green double curve bounds a red disc containing at least two bigons, so all green double curves with triple points are parallel, and bound strips in the red torus. Similarly, all blue double curves with triple points are parallel, and bound strips in the red torus. The total number of triple points is therefore $n^2 \geom{a}{b}$, where $n$ is the number of double curves of each colour with triple points, and $a$ is a green curve with triple points and $b$ is a blue curve with triple points.

If there are three or more curves of each colour with triple points, then there are strips of both colours on both sides of the red torus, so there are intersecting strips, contradicting Lemma \ref{lemma:no strips}, so there may be at most two curves of each colour with triple points.

\begin{subcase}{Two curves of each colour with triple points.}
\end{subcase}

If there are two curves of each colour with triple points, then the red torus is divided up into $\{$ disc, punctured torus $|$ strip $\}$ by each of the green and blue curves with triple points. The total number of triple points is four times the number of intersections between a single green curve and a single blue curve.

So the green inessential disc on one side of the red torus intersects either a blue inessential disc, or a blue strip. If the green inessential disc intersects the blue inessential disc, then by Lemma \ref{lemma:intersecting discs} it does so in one or two arcs, so the boundary curves of the discs intersect in either two or four points, giving eight or $16$ triple points in total, a contradiction. If the green inessential disc intersects the blue strip, then there are green faces coming from the green disc lying inside the inner solid torus bounded by the blue strip. As the boundary components of the blue strip are inessential in the red torus, any green square contained in the inner solid torus bounded by the blue strip creates bigons outside the innermost disc, so there may be no green square faces in the inner solid torus. If there are at most two bigons, then the green disc contains two bigon faces, and all other faces in the green disc are squares. This means there may be at most two green bigons from the green disc in the inner solid torus, which again gives either eight or $16$ triple points, a contradiction.

\begin{subcase}{One curve of each colour with triple points.}
\end{subcase}

If there is only one curve of each colour with triple points, then the
red torus is split into $\{$ disc $|$ punctured torus~$\}$ by each of
the green and blue curves with triple points. If the green disc and
the blue disc lie on the same side of the red torus, then they may
intersect in at most two arcs, by Lemma \ref{lemma:intersecting
  discs}, giving a total of two or four triple points, a
contradiction. So on one side of the red torus there is a green disc
and a blue punctured torus. If there are two bigons, then the
punctured torus may not contain a face with more than eight sides, by
Lemma \ref{lemma:at most octagon}. The green disc has inessential boundary, and divides the solid torus bounded by the red torus into a $3$-ball, which we shall call the inner $3$-ball, and a solid torus, which we shall call the outer solid torus. 

\begin{claim} \label{claim:inner ball}
The blue faces in the inner $3$-ball are either squares or hexagons.
\end{claim}

\begin{Proof}
The rotation lemma, Lemma \ref{lemma:reflection}, implies all boundary components of blue faces inside the inner $3$-ball have the same length, so they are all length at most eight. If they are length two (coming from an annulus consisting of two bigons tubed together), then there are only two or four triple points. There may be only one octagon, which gives eight triple points, which is not divisible by six, so the blue faces in the inner $3$-ball are either squares or hexagons. \end{Proof}

Now consider the blue faces in the outer solid torus. The boundary of the outer solid torus consists of a green disc and a red punctured torus. The boundary of a blue face consists of alternating red and green arcs. The red arcs are all parallel in the green disc, and the green arcs are all essential in the red punctured torus. 

The boundary of a blue square in the outer solid torus may be essential or inessential. If the boundary is inessential, then the disc it bounds in the boundary of the outer solid torus contains either a red square or a green square. If it contains a green square, then the green bigons lie outside the disc, so there are green inessential arcs in the red punctured torus, so it must contain the red square, and hence the two green bigons. These two cases are illustrated below in Figure \ref{picture163}. There is no natural choice of longitude for the outer solid torus, so we have chosen to draw the simplest one.

\begin{figure}[H]
\begin{center}
\epsfig{file=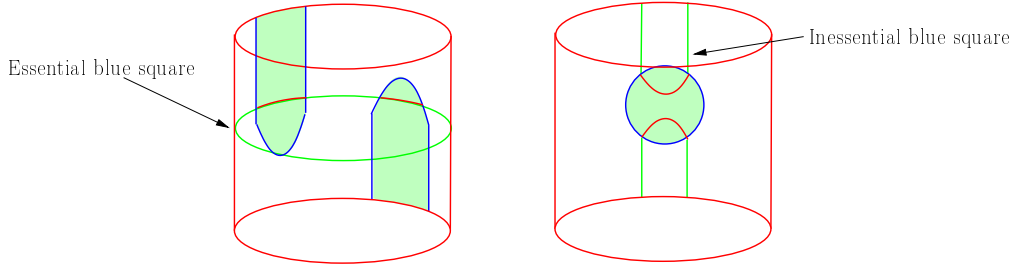, height=100pt}
\end{center}
\caption{Squares in the outer solid torus.} 
\label{picture163}
\end{figure}

In the diagram above, we have only drawn the boundary of the blue square, and we have labelled the blue squares essential or inessential, depending on whether their boundaries are essential or inessential. The top of the red cylinders should be identified, with an appropriate twist in the left hand diagram. We have also shaded the interior of the green disc with blue boundary in the boundary of the outer solid torus, but we have not shaded the complement, which is a red punctured torus with blue boundary.

We next show that all faces in the punctured torus are discs. 

\begin{claim}
All faces in the punctured torus are discs.
\end{claim}

\begin{Proof}
All the arcs in the punctured torus are essential, so if there is a
face which is not a disc, then there is a single annulus with two
bigon boundaries, which implies that all the arcs in the punctured
torus must be parallel. In particular this means that there are no
hexagonal faces, so by Claim \ref{claim:inner ball} all of the faces
in the inner $3$-ball are squares, and so there are an even number of
red arcs in the green disc in the boundary of the inner $3$-ball.
However, all arcs in the green disc must also be parallel, which
implies that all boundary curves of blue faces in the outer solid torus
have the same number of triple points. The collection of blue
faces in the blue punctured torus includes an annulus with bigon boundaries, which does not lie in
the inner $3$-ball, and so lies in the
outer solid torus.  Therefore the boundary curves of the blue faces in the
outer solid torus must all have two triple points. The only faces with
bigon boundaries are the annulus, and at most two bigons, so as there
are an even number of red arcs, there are either four or eight triple
points in total, a contradiction.
\end{Proof}

We have shown that all faces in the punctured torus are discs, so in fact there may be no essential simple closed curves without triple points. We now show that no face can be an octagon.

\begin{claim}
There are no octagons in the punctured torus.
\end{claim}

\begin{Proof} 
If the blue punctured torus contains an octagon, then it must lie in the outer solid torus, by Claim \ref{claim:inner ball}. 

Suppose there is an essential square in the outer solid torus. The boundary of an essential square contains two distinct isotopy classes of green arcs in the punctured torus, so if the punctured torus contains an octagon, these must be the only two isotopy classes of green arc, and so all other green arcs are parallel to these. However, as all the red arcs in the green disc are also parallel, this means that all curves in the boundary of the outer solid torus contain exactly four triple points, so if there is an octagon in the outer solid torus, then there may be no essential squares in the outer solid torus.

Suppose there is an inessential square in the outer solid torus. Consider an innermost inessential square, which bounds a $3$-ball in the outer solid torus, not containing any other blue face. If there is an octagon, then there is exactly one octagon, and no hexagons, so the inner $3$-ball contains squares, by Claim \ref{claim:inner ball}. The innermost inessential face has two red arcs in its boundary, which are innermost in the green disc, and therefore lie in the boundary of a common square. However, the union of the two squares is then an annulus with two green double curves in its boundary, contradicting the fact there is a single double curve of each colour.

We have shown that there may be no squares in the outer solid torus, so it contains a single octagon, so there are eight triple points in total, a contradiction.
\end{Proof}

So the remaining cases are when there are hexagonal faces in the punctured torus. 

If the boundary of a hexagon in the outer solid torus is inessential, then this creates a bigon in the boundary of the outer solid torus whose boundary consists of a green arc and a blue arc, as illustrated in Figure \ref{picture176} below. However this bigon lies in the red punctured torus, creating at least three bigons. So any hexagon in the outer solid torus must have essential boundary.

\begin{figure}[H]
\begin{center}
\epsfig{file=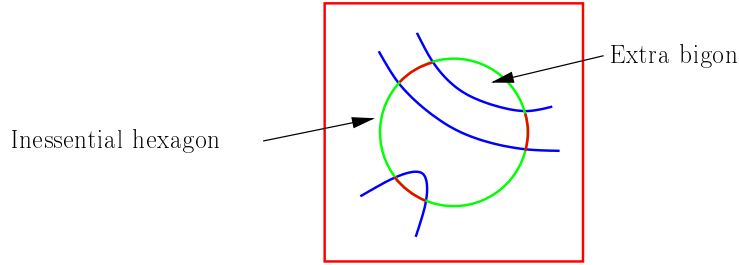, height=100pt}
\end{center}
\caption{The boundary of an inessential hexagon in the boundary of the outer solid torus.} 
\label{picture176}
\end{figure}

If there are hexagons on both sides of the green disc, then the inner
$3$-ball contains a single hexagon, so the outer solid torus also
contains a single hexagon and no squares, as there are exactly three
red arcs in the green disc in the boundary of the outer solid torus.
The single hexagon is a meridional disc which intersects a core curve
for the outer solid torus precisely once, a contradiction, as the blue
surfaces are null homologous. So both of the hexagons lie on the same
side of the green disc.

\begin{sub2case}{Both hexagons are contained in the inner $3$-ball.}
\end{sub2case}

The blue squares in the outer solid torus may have boundaries which are either essential or inessential. If there are hexagons in the punctured torus, then there are three distinct isotopy classes of green arcs in the red punctured torus, so there are both essential and inessential squares. The blue torus is a zero homology class, and in particular intersects the core curve of the outer solid torus an even number of times, so there are an even number of essential squares. Two hexagons in the inner $3$-ball give rise to six red arcs in the green disc in the boundary of the outer solid torus, so there are three squares in the outer solid torus, two of which must be essential, and one is inessential. This is illustrated on the left hand side of Figure \ref{picture177} below.

\begin{figure}[H]
\begin{center}
\epsfig{file=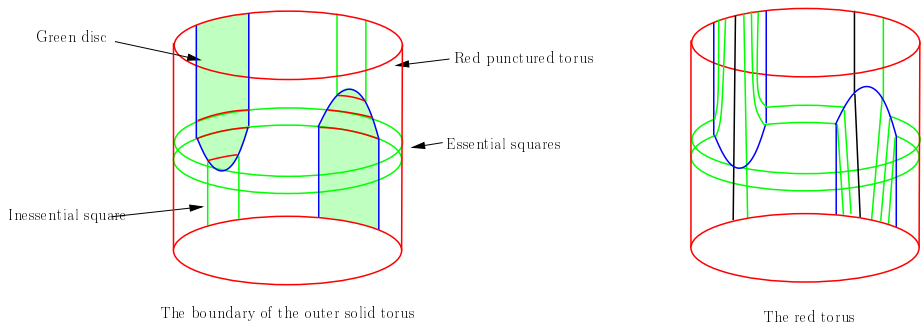, height=120pt}
\end{center}
\caption{Both hexagons contained in the inner $3$-ball.} 
\label{picture177}
\end{figure}

The right hand side of Figure \ref{picture177} shows the corresponding pattern of green arcs in the red torus. However, this means there are four green double curves, a contradiction. In Figure \ref{picture177} we have illustrated this by colouring one of the simple closed curves black instead of green.

\begin{sub2case}{Both hexagons are contained in the outer solid torus.}
\end{sub2case}

Blue hexagons in the outer solid torus must have essential boundary,
as illustrated below in Figure \ref{picture194}, which also shows the
possible essential and inessential squares in the outer solid torus.
In the right hand part of Figure \ref{picture194} the top of the red
square should be identified with the bottom with a one-third shift, so
that the blue arcs are identified to form a simple closed curve. The
outer solid torus may not contain both essential and inessential
squares.

\begin{figure}[H]
\begin{center}
\epsfig{file=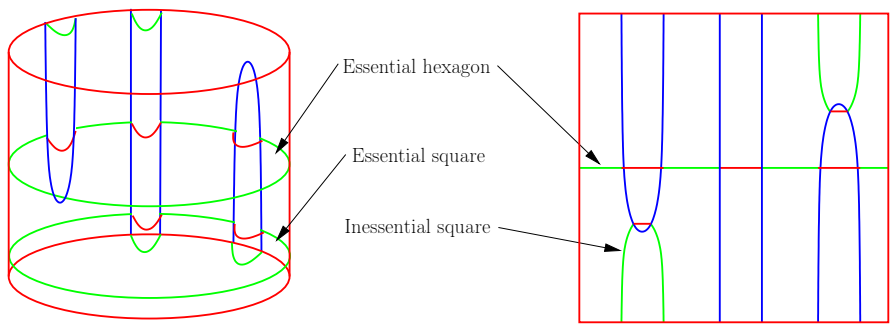, height=120pt}
\end{center}
\caption{An essential hexagon in the outer solid torus.} 
\label{picture194}
\end{figure}

We now show that there may be no essential blue squares in the outer
solid torus. The two essential blue hexagons must be parallel. Suppose
there are $n$ essential blue squares in the outer solid torus. The red
arcs in the green disc in the boundary of the outer solid torus are
all parallel, and starting from an outermost red arc, they consist of
two red arcs contained in the boundary of the two blue hexagons, then
$n$ red arcs contained in the boundary of the blue squares, then
another two red arcs contained in the boundary of the blue hexagons,
then another $n$ red arcs contained in the boundary of the blue
squares, and finally another two red arcs from the blue hexagons.  The
faces inside the inner $3$-ball are all blue squares, and so the
outermost red arcs live in the same blue square, as do the second
outermost pair of arcs, and so on. This means that the six red arcs in
the boundary of the two blue hexagons are equal to the six red arcs
contained in the boundary of three blue squares in the inner $3$-ball.
So the union of the six green arcs in the boundary of the blue
hexagons in the outer solid torus and the six green arcs in the
boundary of the three blue squares in the inner $3$-ball form a single
simple closed green double curve in the red torus. As there is at most
one simple closed curve this means that $n$ must be zero, i.e. there
are no essential squares in the outer solid torus.

If there are inessential blue squares, then each inessential
square is connected to itself by a square in the inner $3$-ball,
forming annuli, so there may be no squares in the outer solid torus.
This means there are exactly 12 triple points, and the inner $3$-ball
contains three squares. The boundary of the inner $3$-ball contains a
simple closed blue curve. On one side the blue curve bounds a red disc
containing parallel green arcs, and on the other side a green disc
containing parallel red arcs. The red and green arcs form the boundary
of three blue squares. This is illustrated below in Figure
\ref{picture178}.

\begin{figure}[H]
\begin{center}
\epsfig{file=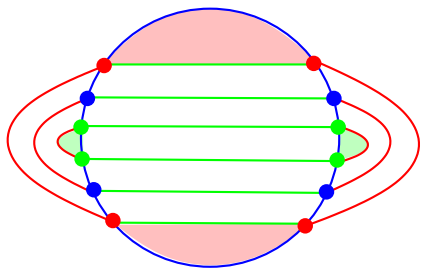, height=100pt}
\end{center}
\caption{The boundary of the inner $3$-ball.} 
\label{picture178}
\end{figure}

The red dots indicate those triple points which are vertices of the red bigons, and the green dots indicate the triple points which are vertices of the green bigons. These two sets are disjoint from each other, so the third images of these points, the boundaries of the blue bigons, are disjoint from the first two, in fact the remaining four triple points, which we have labelled with blue dots. But then the triple points with blue dots all lie on a single red double curve with four triple points, a contradiction.
\end{case}

\begin{case}{The green double curves with triple points divide the red torus into components which include a punctured torus, and the blue double curves with triple points divide the red torus into components which include a pair of pants.}

We may assume that the green curves with triple points are all inessential, and the blue curves with triple points contain an essential curve. The blue inessential curve with triple points bounds a red disc which must contain at least two bigons, so there may be no bigons outside the disc. There must also be an even number of blue essential double curves, but if there is only one essential blue double curve with triple points, then there is an annulus whose boundary consists of two blue curves, one with triple points and the other with no triple points, which creates bigons outside the disc bounded by the blue inessential curve with triple points. So there must be at least two essential blue curves with triple points, so there are at least three curves of each colour with triple points. Furthermore, the annulus bounded by two essential blue double curves with triple points may not contain bigons, so must be a strip. The green double curves with triple points dividing the torus into a punctured torus create strips on both sides of one of the tori. So at least one side of the red torus has both green and blue strips, which must intersect, which implies there are at least three bigons, by Lemma \ref{lemma:no strips}.
\end{case}

\begin{case}{The green double curves with triple points divide the red torus into components which include a punctured torus, and the blue double curves with triple points divide the red torus into components which consist of annuli only.}

Suppose there are three or more curves of each colour with triple points. The innermost green double curve with triple points bounds a red disc containing at least two bigons, so there are no bigons outside the disc, and the green curves with triple points bound a pair of adjacent strips. As there are also three blue double curves with triple points, there must be at least one annulus which is a strip, so there are intersecting strips on at least one side of the red torus, contradicting Lemma \ref{lemma:no strips}.

Suppose there is only one double curve of each colour with triple
points. All the blue arcs inside the red disc bounded by the green
double curve with triple points are parallel. There is an even number
of essential blue double curves, so there must be an essential blue
double curve with no triple points. This means that the green and blue
double curves with triple points are contained in an annulus. There
may be no bigons outside the red disc bounded by the green double
curves, so all the blue arcs outside this red disc must all be
parallel and essential.  But then if there are more than two blue arcs
in the red disc, there are more than two blue double curves with
triple points, so in fact there are at most two blue arcs in the red
disc, so there are at most four triple points in total, a
contradiction.

Finally, suppose there are two double curves of each colour with triple points.

\begin{figure}[H]
\begin{center}
\epsfig{file=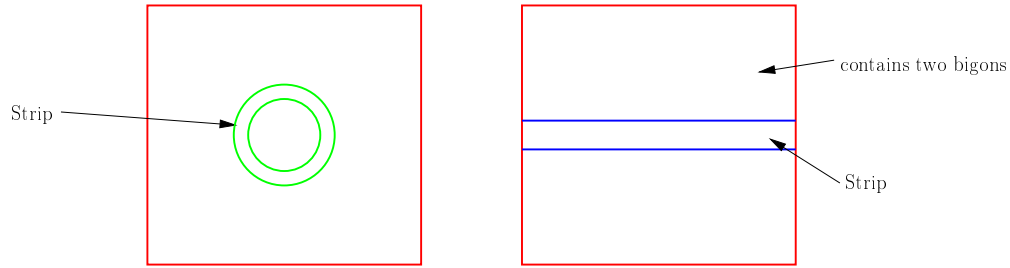, height=100pt}
\end{center}
\caption{The green double curves create a punctured torus, while the blue double curves create annuli only.}
\end{figure}

The two inessential green curves bound a strip, so they both have the
same number of triple points, so the total number of triple points is
four times the number of blue arcs inside the innermost disc. This
also means the blue curves have the same number of triple points, so
the red annuli with blue boundary cannot contain one bigon each, and
so must be a strip and an annulus which contains two bigons. The blue
strip and the green strip lie on opposite sides of the red torus, so
there is a side in which a strip intersects a disc and a punctured
torus. The disc must be essential, so may have at most two arcs of
intersection with the strip, by Lemma \ref{lemma:disc intersect
  strip}, but this implies that there are either two or four triple
points on the inner green curve, which implies there are four or eight
triple points in total, a contradiction.
\end{case}

\begin{case}{Both sets of double curves with triple points divide the red torus into components which contain pairs of pants.}

If there are five or more curves of each colour with triple points, then there must be adjacent red strips with both green and blue boundaries, so at least one side of the red torus has intersecting strips, contradicting Lemma \ref{lemma:no strips}. If there are only two curves of each colour with triple points, then one of the curves must be essential, and the other curve must be inessential. However there must be an even number of essential curves, so there is an essential  double curve with no triple points. This creates an annulus disjoint from the inessential double curve with triple points, which has triple points on one boundary component, but not the other, creating bigons outside the disc bounded by the inessential double curve with triple points, a contradiction. So the cases we need to consider are when there are three or four curves of each colour with triple points.

\begin{subcase}{The configuration contains four double curves of each colour with triple points.}
\end{subcase}

There is only one arrangement of four green curves with triple points which does not create adjacent strips, i.e. there must be two essential curves and two inessential curves. This illustrated below in Figure \ref{picture80}.

\begin{figure}[H]
\begin{center}
\epsfig{file=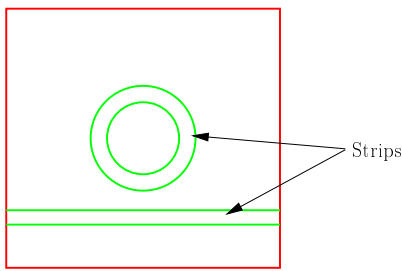, height=100pt}
\end{center}
\caption{The red strips with green boundary both lie on the same side of the blue torus.}
\label{picture80}
\end{figure}

If there are only two bigons, and there are essential curves of both colours, then all faces must be discs. As we will make use of this fact in future cases, we state it as the following lemma.

\begin{lemma} \label{lemma:disc faces}
If there are only two bigons, and the red torus contains essential double curves of both colours, then all faces must be discs.
\end{lemma}

\begin{Proof}
If there is a face which is not a disc, then the red torus contains a simple closed curve (not necessarily a double curve) which is disjoint from all of the double curves, and which either bounds a disc containing triple points in its interior, which creates at least three bigons, or else implies that the essential curves of both colours must have the same slope, contradicting Lemma \ref{lemma:no strips parallel}
\end{Proof}

Therefore, by the argument preceding Case 4.1, there may not be any
essential double curves without triple points, or any inessential
double curves without triple points which contain triple points in the
discs they bound. This means that the two red strips with green
boundary lie on the same side of the blue torus, so the images of
these strips under $g$, which are green strips with blue boundary,
must also lie on the same side of the red torus.

The same argument may be applied to the blue double curves, so there is a pair of red strips with blue boundary which lie on the same side of the green torus, and the images of these under $g^2$, which are blue strips with green boundary, must lie on the same side of the red torus. However, by Lemma \ref{lemma:no strips}, the strips of different colours may not intersect, so the pair of green strips with blue boundary must lie on the other side of the red torus from the pair of blue strips with green boundary. So on one side of the red torus we have a pair of green strips intersecting a blue disc and a blue pants. 

If both green strips have exactly one inessential boundary component, then by Lemma \ref{lemma:disc intersect strip}, the blue disc is disjoint from the green strips, so doesn't have any triple points, a contradiction. So we may assume that one strip has both boundary components essential, and the other strip has both boundary components inessential. We now show that as the green strip with inessential boundary bounds a red strip in the the red torus, the union of these two strips is in fact an unknotted solid torus. As we use this fact later on, we will state it as the following lemma.

\begin{lemma} \label{lemma:torus bigon}
Suppose a green strip and a red strip have two blue double curves in common as their boundary. Suppose further that the blue double curves are either both inessential, or both meridional on the side of the red torus containing the green strip, and that all faces are discs. Then the torus formed by the union of the two strips bounds an unknotted solid torus containing blue bigon faces only.  
\end{lemma}

\begin{Proof}
The torus formed by the union of the green strip and the red strip has a compressing disc with boundary parallel to the blue double curve boundary components of the strip. We will call the side of the torus containing the compressing disc the outside, and the other side the inside. All the blue faces inside the torus formed by the union of the two strips have parallel boundary, with the same number of red and green edges. Furthermore, the boundaries of the blue faces cross the blue double curves with the same orientation at each crossing, and are therefore compressing discs for the inside region. As compressing a torus creates a sphere, the inside region must be a solid torus, as the torus is contained in the solid torus bounded by the red torus, which is irreducible. This solid torus is unknotted, as there is a compressing disc on the outside. Therefore the compressing disc on the outside of the torus intersects each blue face inside the torus with the same orientation, so if there is a face with more than two sides, this creates a lens space connect summand inside the solid torus, contradicting irreducibility, so all faces inside the torus must be bigons. 
\end{Proof}

We have shown that all faces inside the solid torus bounded by the strips with inessential boundary are bigons, so as the blue pants contains no bigons, this means that the blue pants must be disjoint from the inessential blue curves. However, the blue pants contains three of the green double curves in its boundary, which means that only one green double curve may intersect the blue inessential double curves, which is a contradiction, as the double curves come in two parallel pairs.

\begin{subcase}{The configuration contains three double curves of each colour with triple points.}
\end{subcase}

All faces are discs, by Lemma \ref{lemma:disc faces}, so there are no essential double curves without triple points. There must be an even number of essential double curves, so there must be two essential double curves with triple points of each colour, and a single inessential double curve with triple points of each colour, as illustrated for the green double curves with triple points in Figure \ref{picture92} below.

\begin{figure}[H]
\begin{center}
\epsfig{file=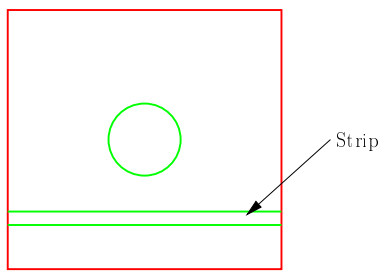, height=100pt}
\end{center}
\caption{Three double curves.} \label{picture92}
\end{figure}

The curves divide the red torus into $\{$ disc, strip $|$ pants $\}$, and the strips of different colours must lie on opposite sides of the red torus, by Lemma \ref{lemma:no strips}. We now show that the green disc must be inessential in the solid torus bounded by the red torus.

\begin{claim}{The green disc is inessential.}
\end{claim}

\begin{Proof}
Suppose the green disc is essential. Then in this case, the green strip has exactly one inessential boundary component. This is shown in Figure \ref{picture82} below.

\begin{figure}[H]
\begin{center}
\epsfig{file=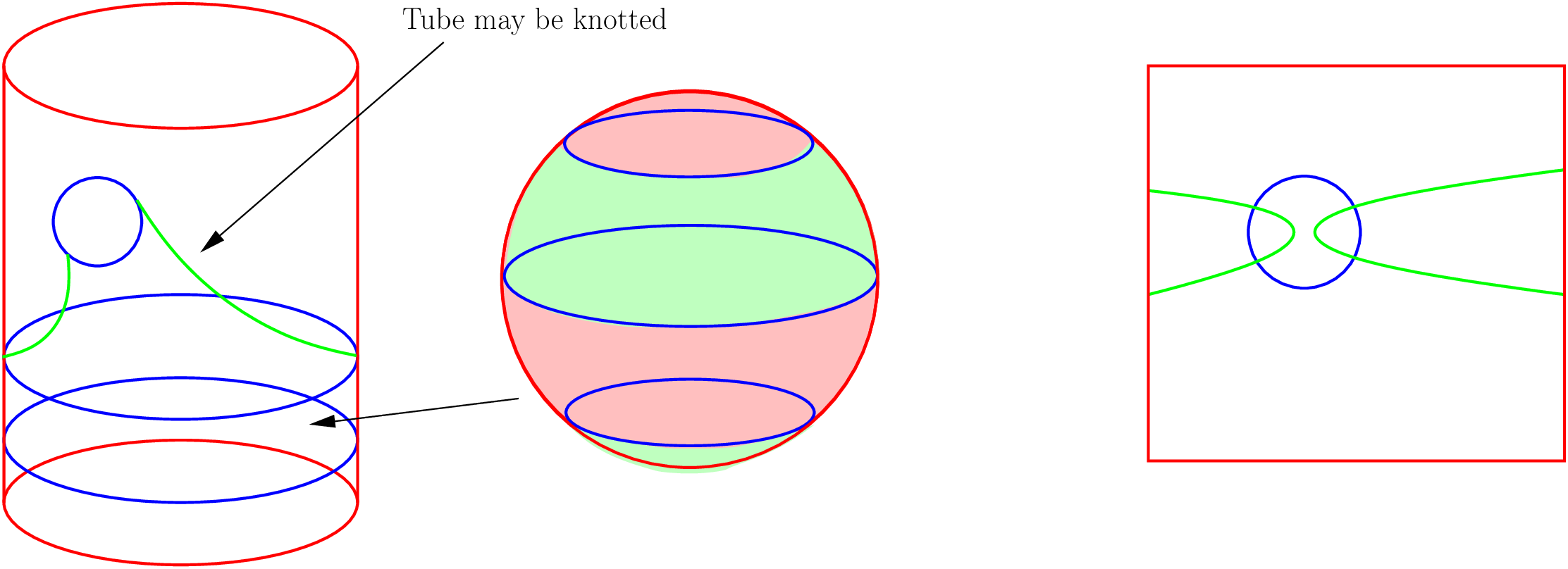, height=100pt}
\end{center}
\caption{The green disc is essential.}
\label{picture82}
\end{figure}

One of the complementary regions is a $3$-ball, whose boundary is divided up into a disc and an annulus of each colour. All arcs are parallel in the discs and the annuli, so this $3$-ball can contain hexagons only. The inessential blue double curve bounds a disc in the red torus containing two bigons, and these two bigons must also lie in the disc bounded by the inessential green curve, so the blue inessential double curve must have at least four triple points. This means that there must be two hexagons in the $3$-ball, and so in fact every blue double curve has precisely four triple points. This implies that each of the green double curves also has precisely four triple points, so the green and blue inessential curves may only intersect each other, as shown on the right hand side of Figure \ref{picture82}. However, this implies that the green and blue essential double curves have the same slope, contradicting Lemma \ref{lemma:no strips parallel}.
\end{Proof}

So we may now assume that the green disc is inessential.

The inner solid torus bounded by the green strip contains either the
red strip, or the inessential green disc, in its boundary. If it
contains the red strip, then, as an inner solid torus has meridional
discs that hit each boundary curve of the strips once, all the blue
disc faces inside the inner solid torus are bigons. By Lemma
\ref{lemma:no strips}, the green and blue strips lie on opposite sides
of the red torus, so the blue faces in the inner solid torus must lie
in the blue pants. However the blue pants contains no bigons, so the
inner solid torus contains the inessential green disc. The inessential
green disc divides the inner solid torus into two pieces. We will call
the piece which is a $3$-ball the inner $3$-ball, and the other piece,
which is a solid torus, the middle solid torus.  The middle solid
torus separates the inner $3$-ball from the outer solid torus.

We now describe all possible blue squares in the middle solid torus.
If a blue square has both red boundary arcs in the strip, then both
green arcs must connect the two essential blue double curves. However
this gives a disc with intersection number two with one of the blue
double curves, a contradiction. If a blue square has one red boundary
arc in the green strip, and one in the green disc, then the two green
arcs must each have one end point on the inessential blue double
curve, and one endpoint on an essential blue double curve. In this
case the blue square is a meridional disc for the middle solid torus,
and we will refer to squares of the sort as essential squares.
Finally, if a blue square has both red boundary arcs contained in the
green disc, then the two green boundary arcs must be essential arcs
which connect the inessential blue double curve to itself, and we will
call squares of this sort inessential squares. These two possibilities
are illustrated below in Figure \ref{picture17}.

\begin{figure}[H]
\begin{center}
\epsfig{file=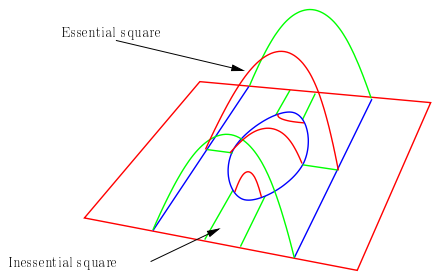, height=100pt}
\end{center}
\caption{Squares in the middle solid torus.}
\label{picture17}
\end{figure}

The pair of pants contains either a single octagon or two hexagons.
Recall that by Lemma \ref{lemma:no strips}, the green and blue strips
may not lie on the same side of the red torus, so the blue pants lies
on the same side of the red torus as the green strip. We now consider
each case in turn.

\begin{sub2case}{The pair of pants contains an octagon.}
\end{sub2case}

Suppose the octagon is contained in the inner $3$-ball. This creates exactly eight triple points on the inessential blue double curve. As the pair of pants contains an octagon, there are exactly two isotopy classes of essential arc in the pair of pants, so the middle solid torus contains essential squares only. This means that there are precisely four triple points on each essential blue double curve, giving a total of 16 triple points, which is not divisible by six, a contradiction.

Suppose the octagon is contained in the outer solid torus, then there may be no other squares in the outer solid torus, so each blue essential double curve contains precisely four triple points. Again, there may only be essential squares in the middle solid torus, so the inessential blue double curve contains precisely eight triple points, giving a total of 16 triple points, which is not divisible by six, a contradiction.

So the octagon is contained in the middle solid torus. If there are no
other faces in the middle solid torus, then there are eight triple
points in total, which is not divisible by six, a contradiction. The
only other faces in the middle torus may be squares, and there may be
no inessential squares, as they connect a boundary curve of the red
pants to itself, dividing the red pair of pants into hexagons, not
octagons. However, the red pair of pants is an image of the blue pair
of pants under $g$, which contains an octagon, so this cannot happen.
Therefore there may only be essential squares. As the pants contains
an octagon, there are only two isotopy classes of green double arcs in
the red pants, so the green boundary arcs of the octagon in the red
pants must be parallel to the green boundary arcs of the essential
square. As the blue essential double curves bound a green strip in
which each red arc runs from one boundary component to the other, the
green boundary arcs of the octagon must have the same number of
endpoints in each essential blue curve, so the octagon must have two
green boundary arcs parallel to each green boundary arc of an
essential blue square. This is illustrated below in Figure
\ref{picture187}, which shows the arrangement of green and red double
arcs in the boundary of the middle solid torus, which is a torus
consisting of a green strip, a red pair of pants and a green disc.

\begin{figure}[H]
\begin{center}
\epsfig{file=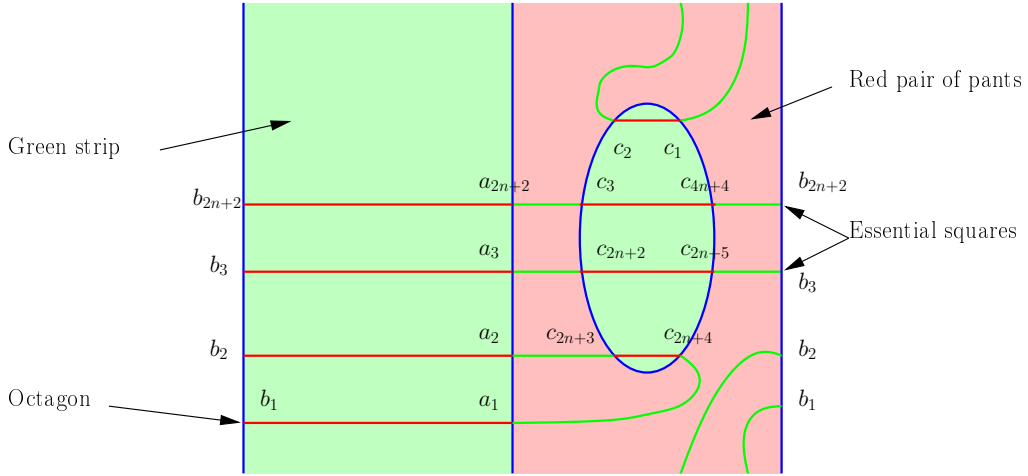, height=180pt}
\end{center}
\caption{The boundary of the middle solid torus.}
\label{picture187}
\end{figure}

In Figure \ref{picture187} above we have only drawn two essential
squares, there may be arbitrarily many parallel essential squares.

The inner $3$-ball contains squares only, so there must be an even
number of red arcs in the disc, hence an even number of essential
squares, $2n$ say, in the middle solid torus. We will label the triple
points on the two essential blue curves $a_i$ and $b_i$, respectively,
and we will label the triple points on the inessential blue curve
$c_i$. If there are $2n$ squares in the middle solid torus, then we
will label the points on one double curve $a_1, \ldots, a_{2n+2}$ in
order going around the curve, starting with $a_1$ and $a_2$ being
adjacent triple points which lie in the boundary of the octagon.
Similarly, we will label the triple points on the other blue double
curve $b_1, \ldots, b_{2n+2}$, so that $b_1$ and $b_2$ are adjacent
triple points on the blue double curve, and they both lie in the
boundary of the octagon, and the induced orientation on the blue
double curves given by going from either $a_1$ to $a_2$ or $b_1$ to
$b_2$ is the same. We will label the triple points on the inessential
blue double curve $c_1, \ldots c_{4n+4}$, in order along the double
curve, so that the triple points lying in the outermost red double
arcs are labelled $c_1, c_2, c_{2n+3}$ and $c_{2n+4}$, and
furthermore, so that the labelling of the triple points in the
boundary of the octagon is $\{c_1, c_2, b_2, a_2, c_{2n+3}, c_{2n+4},
a_1, b_1\}$, in that order.

As the inner $3$-ball and the outer solid torus contain blue squares
only, this determines the pattern of green arcs in the red torus. In
the inner $3$-ball, the red arcs connect $c_{k}$ to $c_{4n+7-k}$,
taking the indices modulo $4n+4$, while in the outer solid torus the
green arcs connect $a_{k}$ to $b_{k+n+1}$, taking the indices modulo
$2n+2$. This is illustrated below in Figure \ref{picture188}. In the
figure below, a shaded green region between two green arcs indicates
some number of parallel arcs. If the region is labelled $n$, then
there are $n$ parallel arcs in the region, including the two green
boundary arcs which demarcate the region.

\begin{figure}[H]
\begin{center}
\epsfig{file=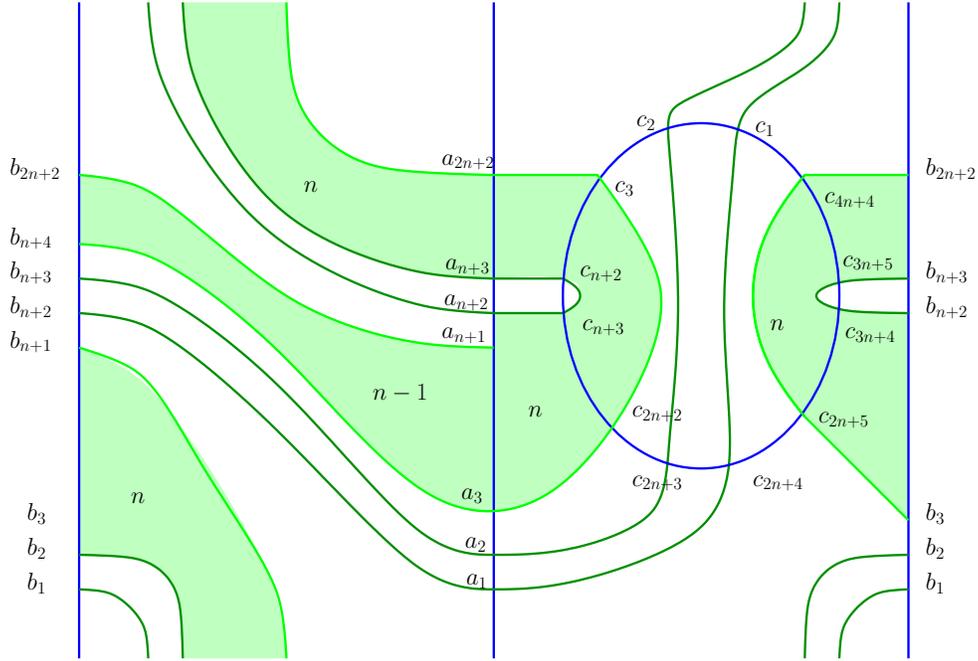, height=250pt}
\end{center}
\caption{The green and blue double curves in the red torus.}
\label{picture188}
\end{figure}

The outermost pair of red arcs in the disc bounded by the inessential
blue curve lie in the boundary of a common square, so there are green
arcs connecting $c_1$ and $c_{2n+4}$, and $c_2$ and $c_{2n+3}$. As all
the blue squares are parallel inside the inner $3$-ball, this means
that $c_k$ is connected to $c_{2n+5-k}$ by a green arc, taking the
indices modulo $4n+4$. As the blue squares inside the outer solid
torus are all parallel essential discs, the green arcs inside the
annulus connect $a_k$ to $b_{k+n+1}$, modulo $2n+2$.

These identifications mean that in fact the triple point $a_1$ lies on a simple closed green double curve with exactly $16$ triple points. To be precise, the identifications described in the preceding paragraph show that this double curve contains the triple points $\{ a_1$, $b_{n+2}$, $c_{3n+4}$, $c_{3n+5}$, $b_{n+3}$, $a_2$, $c_{2n+3}$, $c_2$, $b_2$, $a_{n+3}$, $c_{n+2}$, $c_{n+3}$, $a_{n+2}$, $b_1$, $c_1$, $c_{2n+4} \}$. This curve is drawn using a darker green line in Figure \ref{picture188} above.

The blue double curves have either $2n+2$ or $4n+4$ triple points, so $n$ may be three or seven, but the total number of triple points is $8n+8$, and this gives either $32$ or $64$ triple points, neither of which is divisible by six.

\begin{sub2case}{The pair of pants contains hexagons.}
\end{sub2case}

First we show there are no hexagons in the inner $3$-ball. If there is a single hexagon in the inner $3$-ball, then all faces in the inner $3$-ball are hexagons. If there is only one hexagon, then the inessential blue double curve has six triple points. As the pair of pants contains hexagons, it contains three isotopy classes of essential arcs, so there must be both essential and inessential squares in the middle solid torus. So four of the triple points on the inessential double curve must lie in the boundary of an inessential square, and two lie in the boundary of an essential square, which has two other triple points, one on each of the other blue double curves. However this means that there is a total of eight triple points, which is not divisible by six, a contradiction. Suppose both hexagons must lie in the inner $3$-ball. As before, both essential and inessential squares must occur in the middle solid torus. If there is exactly one inessential square, then the arcs in the red pants coming from the inessential square connect four triple points together in the inessential blue double curve, so there are four triple points on each of the essential curves, giving $20$ triple points in total. If there are two essential squares, then there are two triple points on each essential curve, giving $16$ triple points in total. In either case, the total number of triple points is not divisible by six, so there may be no hexagons in the inner $3$-ball.

We now show that there are no inessential squares in the middle solid
torus. As there are no hexagons in the inner $3$-ball, the inner
$3$-ball contains squares only. If there are inessential squares in
the middle solid torus, then there is an innermost one which bounds a
$3$-ball containing no other faces. The two red arcs in the boundary
of the innermost inessential square are innermost on the green disc in
the boundary of the inner $3$-ball, and so lie in the boundary of a
common blue square in the inner $3$-ball, and so the green arcs in the
boundary of these two squares form two green simple closed curves with
exactly two triple points each. Furthermore, this means that there may
be no other green arcs in the red pants parallel to the green arcs in
the boundary of the inessential square, as this would create at least
one more green double curve with two triple points, but this would
imply that all three green double curves have exactly two triple
points, and therefore all three blue double curves have exactly two
triple points, a contradiction, as the blue inessential double curve
has at least four. The image of these two green double curves with two
triple points under $g$ consists of two blue double curves with triple
points. The inessential blue double curve has at least four triple
points, so the two essential blue must each have two triple points. In
particular, this means that the outer solid torus contains exactly
four triple points in its boundary, and so contains a single blue
square. As all the green arcs in the red pants with blue boundary are
essential, this means that for each essential blue double curve there
is a single pair of greens arcs in the red pants connecting it to the
inessential blue curve. This gives exactly six green arcs in the red
pants with blue boundary, two of which lie in the boundary of a blue
square. However, this leaves exactly four green arcs in the boundary
of the middle solid torus, which are too few to form the boundary of
the two blue hexagons which must be contained in the middle solid
torus, a contradiction. So there are no inessential squares in the
middle solid torus.

The green arcs in the boundary of an essential square in the middle
solid torus give rise to two isotopy classes of green arc in the red
pants. However, as the red pants contains hexagons, there must be
three isotopy classes for green arcs in the red pair of pants, so
there must be at least one green arc in the red pants not contained in
a blue square in the middle solid torus. This means at least one
hexagon must be contained in the middle solid torus, containing a
green boundary arc that either connects the two essential blue double
curves, or connects the inessential blue double curve to itself.

If the hexagon contains a green arc that connects the two essential double curves, then the hexagon contains at least two red arcs lying in the green strip. If all three of the hexagon's red arcs lie in the green strip, then all three green arcs must connect the essential blue double curves, and so the boundary of the hexagon has intersection number three with the blue double curves in the middle solid torus, a contradiction. If the hexagon has exactly two red arcs in the green strip, then they each have an endpoint on the green arc which connects the two essential blue double curves, so the hexagon has intersection number two with a blue essential double curve in the middle solid torus, which is again a contradiction. So the hexagon has at most one red arc in the strip, and so must contain a green arc with both endpoints in the inessential blue double curve. All three red arcs may not be contained in the green disc, so the hexagon must have one red arc in the strip, two red arcs in the disc, and at least one green arc with both endpoints in the blue inessential curve. In this case, the remaining two green arcs must each have one endpoint in the inessential blue double curve, and one endpoint in one of the essential blue double curves, so up to isotopy there is only one possible hexagon in the middle solid torus, illustrated below in Figure \ref{picture189}.

\begin{figure}[H]
\begin{center}
\epsfig{file=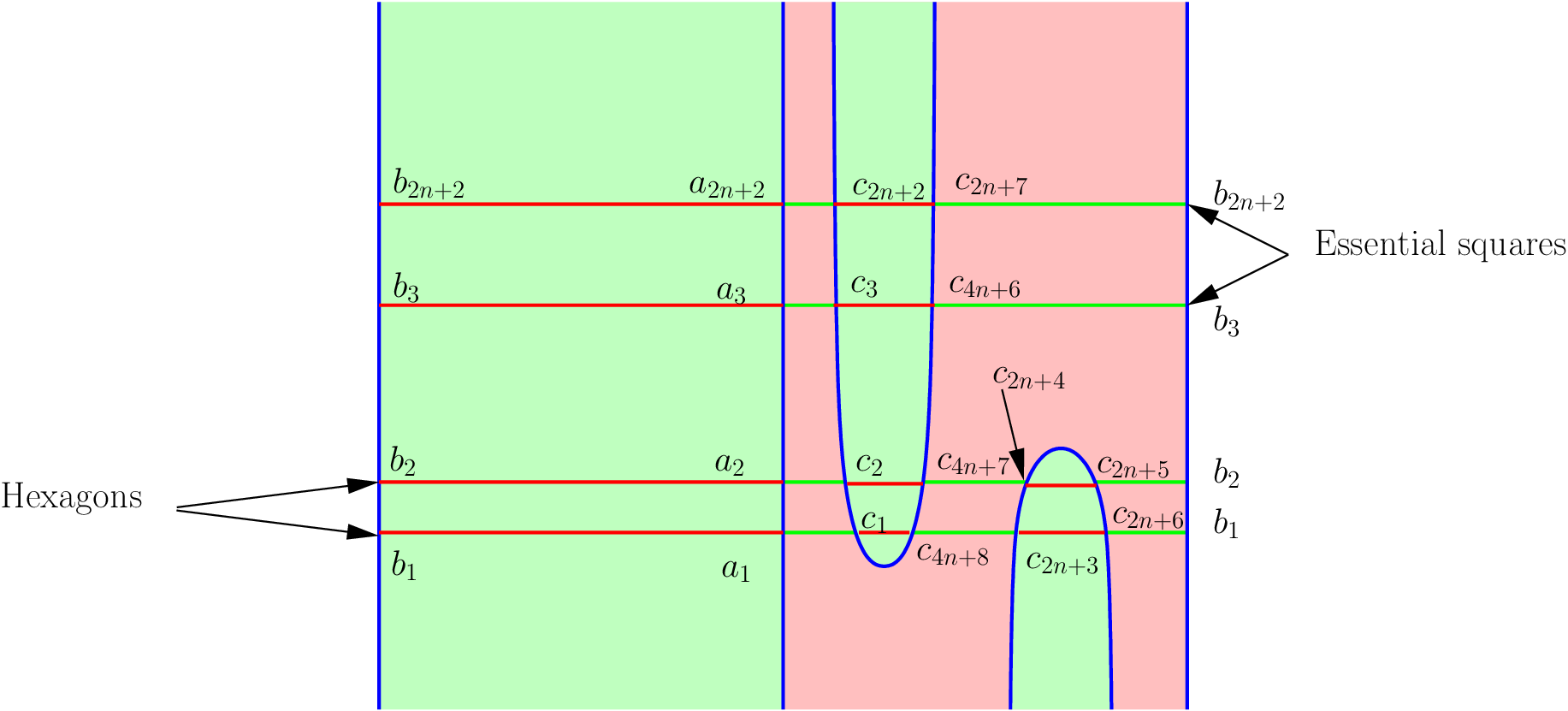, height=180pt}
\end{center}
\caption{The boundary of the middle solid torus.}
\label{picture189}
\end{figure}

The inner $3$-ball contains squares, so there is an even number of red arcs in the green disc. Each blue square in the middle solid torus contains a single red boundary arc in the green disc, and has a single triple point on each blue essential double curve. Each blue hexagon in the middle solid torus contains two red boundary arcs in the green disc, and has a single triple point on each blue essential double curve. There must be an even number of triple points on each double curve, so in fact there must be two parallel blue hexagons in the middle solid torus. We will assume there are $2n$ essential squares, so the essential blue double curves each have $2n+2$ triple points, and the inessential blue double curve has $4n+8$ triple points.

We will label the triple points as follows. Start with one of the essential blue double curves, and label the adjacent triple points lying in the two blue hexagons $a_1$ and $a_2$. Now label the remaining triple points in the essential blue double curve $a_3, \ldots a_{2n+2}$, following the circular ordering. There are two possible choices of labelling consistent with this, choose the one so that the opposite endpoints of the green arcs meeting the adjacent triple points $a_1$, $a_2$ and $a_3$ are also adjacent on the inessential blue double curve. Now label the triple points in the other essential blue double curve $b_1, \ldots b_{2n+2}$, so that $b_i$ and $a_i$ lie in a common blue face in the middle solid torus for all $i$. Finally, label the triple points in the inessential blue double curve $c_1, \ldots c_{4n+8}$, following the circular ordering, so that $a_i$ and $c_i$ are the end points of a green arc in the red pair of pants, for $1 \leqslant i \leqslant 2n+2$. This labelling is illustrated above in Figure \ref{picture189}.

The squares in the inner $3$-ball and the outer solid torus determine
the pattern of green arcs in the red torus. The green arcs in the
boundary of the inner $3$-ball connect $c_{k}$ to $c_{2n+5-k}$, taking
indices modulo $4n+8$, while the green arcs in the boundary of the
outer solid torus connect $a_k$ to $b_{k+n+1}$, taking indices modulo
$2n+2$. This is illustrated in Figure \ref{picture190} below.

\begin{figure}[H]
\begin{center}
\epsfig{file=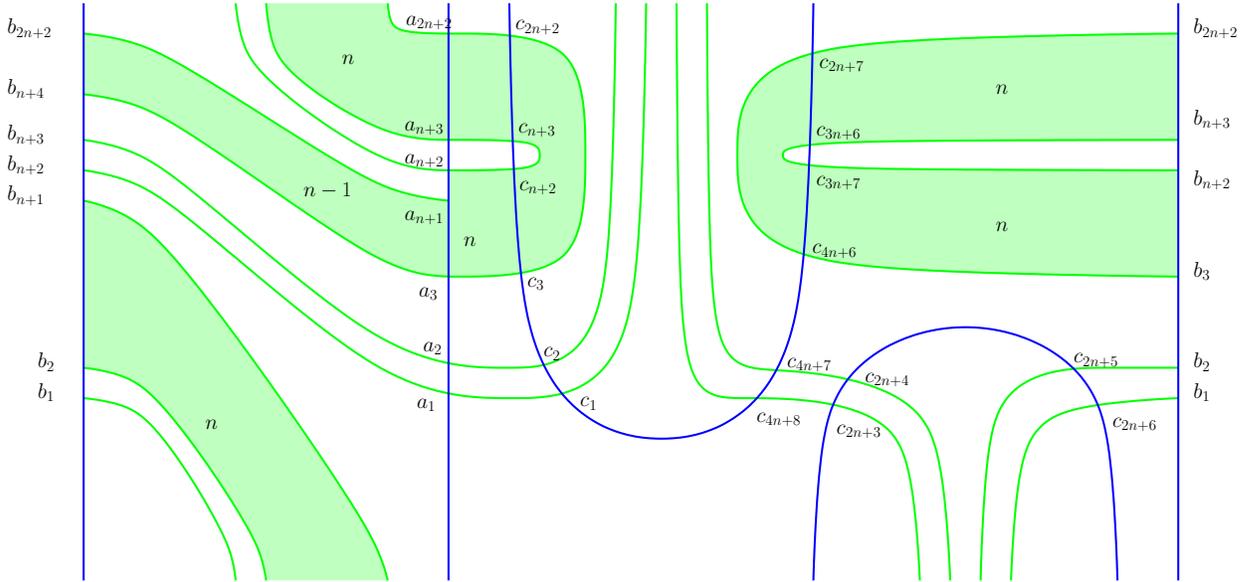, height=220pt}
\end{center}
\caption{The green double curves in the red torus.}
\label{picture190}
\end{figure}

The following eight points all lie on a single green double curve,
which contains no other triple points: $\{ a_3, c_3, c_{2n+2},
a_{2n+2}, b_{n+1}, c_{3n+8}, c_{3n+5}, b_{n+4}\}$. The number of
triple points on a simple closed curve is either $4n+8$ or $2n+2$,
which means $n$ must be either zero or three. As $n$ is not zero, this
implies that $n=3$, and that the two essential blue double curves have
eight triple points, and the inessential blue double curve has 20
triple points. We will call the green double curve with 20 triple
points $a$.  Although we now know that $n=3$, we will continue to
refer to triple points as $c_{n+2}$ rather than $c_5$ in order to aid
references to Figure \ref{picture190}.

The green double curve $a$ contains the following triple points: $\{
a_1$, $c_1$, $c_{2n+4}$, $c_{4n+7}$, $c_{2n+6}$, $b_1$, $a_{n+2}$,
$c_{n+2}$, $c_{n+3}$, $a_{n+3}$, $b_2$, $c_{2n+5}$, $c_{4n+8}$,
$c_{2n+3}$, $c_2$, $a_2$, $b_{n+3}$, $c_{3n+6}$, $c_{3n+7}$, $b_{n+2}
\}$.  The blue inessential double curve also contains 20 triple
points, so this curve is $ga$.

The four points $\{c_{n+2}, c_{n+3}, c_{3n+6}, c_{3n+7} \}$ are the vertices of the red bigons, while the four points $\{ c_1, c_{4n+8}, c_{2n+4}, c_{2n+5} \}$ are the vertices of the green bigons, contained in the green disc. All of these points lie in $a \cap ga$, so the vertices of the blue bigons must also lie in $a \cap ga$, and as the red and green bigons have distinct vertices, the vertices of the blue bigons are also distinct. As $a \cap ga$ consists of exactly twelve triple points, this means that the vertices of the blue bigons are $\{ c_2, c_{2n+3}, c_{2n+6}, c_{4n+7} \}$. Note that $c_{2n+4}$ and $c_{4n+7}$ are vertices of green and blue bigons respectively, and are connected by a green double arc. So the image of the this green double arc under $g$ is a blue double arc connecting the blue and red bigons, but this gives a contradiction, as the vertices of the red bigons are not adjacent to the vertices of the blue bigons along the blue double curve $ga$.
\end{case}

\begin{case}{The green double curves with triple points divide the red torus into components which include a pair of pants, and the blue double curves with triple points divide the red torus into components all of which are annuli.}

\begin{figure}[H]
\begin{center}
\epsfig{file=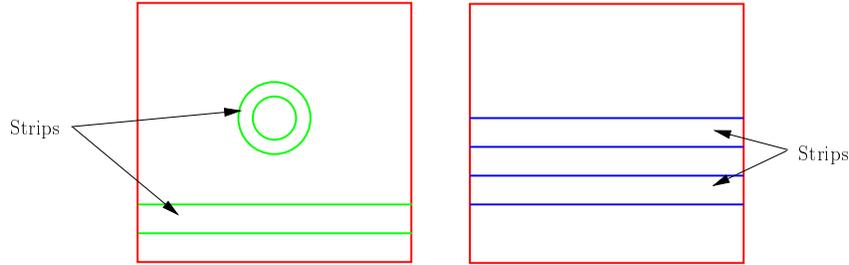, height=100pt}
\end{center}
\caption{Pants and annuli.} \label{picture94}
\end{figure}

First we show that there must be exactly four double curves of each colour.
If there are five or more double curves with triple points, then there is a triple of parallel green curves bounding strips whose images under $g$ give rise to green strips on both sides of the red torus. The blue double curves bound at least one strip, so there are intersecting strips on at least one side of the red torus, contradicting Lemma \ref{lemma:no strips}, so there may be at most four double curves with triple points. By Lemma \ref{lemma:no strips parallel}, the essential double curves have different slopes, so every essential double curve contains triple points. In particular this means that there are an even number of essential curves with triple points of each colour. This means there must be at least three green curves with triple points, so there must be exactly four double curves with triple points of each colour, as illustrated in Figure \ref{picture94} above.

We now obtain a contradiction by showing that all double curves must
have the same number of triple points. By Lemma \ref{lemma:no strips},
the green and blue strips lie on opposite sides of the red torus.
Consider the side of the red torus in which two blue strips with green
boundary intersect a green disc and a green pants with blue boundary.
Suppose one of the blue strips with green boundary has both curves
inessential in the red torus. Then this only contains bigons
boundaries by Lemma \ref{lemma:torus bigon}, so the inessential curves
have two triple points each. But the blue double curves come in pairs
that bound strips, so the inessential green curves must have at least
four triple points, a contradiction.  So each blue strip must have one
boundary curve essential, and one boundary curve inessential.  This
means that the union of the strips forms a torus, so all double curves
must have the same number of triple points.  This is a contradiction,
as some of the annuli are 1bigons, which have different number of
triple points on each boundary component.
\end{case}

\begin{case}{Both sets of double curves with triple points divide the red torus into annuli only.}

First we show there are exactly two double curves of each colour.

\begin{claim}
There are exactly two double curves of each colour.
\end{claim}

\begin{Proof}
If there are five or more double curves with triple points, then there is a pair of adjacent strips, so there are intersecting strips on both sides, contradicting Lemma \ref{lemma:no strips}, so there may be at most four double curves of each colour with triple points. The blue and green essential curves may not have the same slope, by Lemma \ref{lemma:no strips parallel}, so all faces are discs, and there are an even number of double curves of each colour with triple points. So there are either two or four double curves of each colour with triple points. 

We now show that there are exactly two double curves with triple
points. Suppose there are four double curves with triple points, then
in order to avoid intersecting strips, the double curves of each
colour divide the red torus into the following collection of surfaces:
$\{$ 1bigon, 1bigon $|$ strip, strip $\}$. Consider an innermost green
strip with blue boundary, which divides the solid torus on one side of
the red torus into inner and outer solid tori. The inner solid torus
has a red annulus in its boundary, which has the same number of triple
points on each boundary component, so the red annulus is a strip as
well. Then as all faces are discs, Lemma \ref{lemma:meridional bigons}
implies that each blue face in the inner solid torus is a bigon, so
there are at most two essential arcs in each strip, but as the green
and blue double curves have different slopes, each curve intersects
the others, so there are at least four arcs of intersection, a
contradiction.
\end{Proof}

The two essential double curves divide the red torus into two annuli. The two annuli are either both 1bigons, or a strip and a 2bigon. Every blue curve is the image of a green curve, so if both the green curves have the same number of triple points, then so do both the blue curves, so if the green double curves divide the red torus into a strip and a 2bigon, then the blue double curves also divide the red torus into a strip and a 2bigon. Similarly, if the green double curves divide the red torus into a pair of 1bigons, then the blue double curves also divide the red torus into a pair of 1bigons.

We now show that both annuli must be 1bigons. 

\begin{claim}
Both annuli are 1bigons.
\end{claim}

\begin{Proof}
Suppose the annuli consist of a strip and a 2bigon. The green and blue
strips cannot lie on the same side of the red torus by Lemma
\ref{lemma:no strips}, so one side of the red torus contains a green
strip and a blue 2bigon, which intersect in red double arcs. As one of
the annuli is a strip, both double curves have the same number of
triple points, so the 2bigon contains an innermost bigon adjacent to
each of the boundary components.

Consider the solid torus on the side of the red torus which contains
the green strip with blue boundary. The green strip divides this solid
torus into inner and outer solid tori. If the boundary of the inner
solid torus is a union of two strips, then all blue faces in the inner
solid torus are bigons, by Lemma \ref{lemma:meridional bigons}. As
there are two bigons, this gives at most four triple points in total,
a contradiction. So the boundary of the inner solid torus consists of
a green strip and a red 2bigon, and the boundary of the outer solid
torus consists of a green strip and a red strip. This implies that all
blue faces in the outer solid torus have the same number of triple
points. If the faces in the outer solid torus are bigons, then as
there are two bigons, this again implies there are at most four triple
points, a contradiction, so in fact both bigons are contained in the
inner solid torus.

As all faces are discs, the blue 2bigon contains either two hexagons,
or an octagon. We now show that in either case, these faces are
contained in the inner solid torus. All blue faces in the outer solid
torus have the same number of triple points in their boundaries, and
the possible squares and hexagons are illustrated in Figure
\ref{picture107}. Therefore, the outer solid torus may not contain an
octagon, as this implies there are exactly eight triple points. There
may not be a single hexagon in the outer solid torus, as then the
other hexagon and the two bigons would lie in the inner solid torus,
giving different numbers of red arcs in the green strip. If both
hexagons lie in the outer solid torus, then the inner solid torus
contains two blue bigons and two blue squares. This determines the
arrangement of green double curves in the red torus, illustrated below
in Figure \ref{picture205}.  The labels refer to which blue faces
contain the green double arcs in the 2bigon. All the green double arcs
in the strip are contained in the hexagons.

\begin{figure}[H]
\begin{center}
\epsfig{file=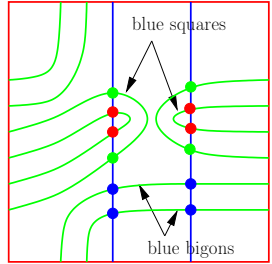, height=100pt}
\end{center}
\caption{Two hexagons in the outer solid torus.}
\label{picture205}
\end{figure}

In Figure \ref{picture205} above, we have labelled the vertices of the
red bigons with red dots, and we have labelled the vertices of the
blue bigons with blue dots. For a given bigon, the orbits of its
vertices are distinct, and so consist of exactly six points. As all the
red and blue bigons are disjoint from each other, this implies that
the green bigons are disjoint from red and blue bigons as well, and so
the remaining four triple points, marked with green dots above, must
be the vertices of the two green bigons. However, the two vertices of
a green bigon must be adjacent along a blue double curve, and none of
the triple points marked with green dots in Figure \ref{picture205}
above have this property, so the outer solid torus may not contain a
pair of hexagons.

We have shown that the blue hexagons or octagons are contained in the
inner solid torus. These faces are disjoint from the blue bigons, so
their boundaries are inessential curves in the boundary of the inner
solid torus, and we now deal with each case in turn.

The boundary of a hexagon consists of three essential red arcs in the
green strip, together with three green arcs in the red 2bigon. As the
boundary of the hexagon is inessential, there is at least one
inessential green arc connecting a pair of endpoints of red arcs on
the same blue double curve. The remaining triple point on that blue
double curve must be connected to a triple point on the other blue
double by an essential green arc, and so the remaining two triple
points on the second blue double curve must be connected by an
inessential green arc. This corresponds to the hexagon illustrated in
Figure \ref{picture104}. However, there must be two hexagons in the
inner solid torus, and this cannot happen if all the inessential green
arcs with endpoints in the same blue double curve are parallel.

The boundary of an octagon consists of four essential red arcs in the
green strip. The boundary of the octagon is inessential in the inner
solid torus, so contains at least one green arc which is inessential
in the red 2bigon. If the boundary of the octagon contains two
inessential green arcs in the red 2bigon with endpoints in the same
blue double curve, then the other two green arcs in the red 2bigon are
also inessential, with endpoints in the other blue double curve.
However, as the inessential green arcs must be parallel in the red
2bigon, this creates two squares instead of an octagon. So the
boundary of the octagon contains one green arc with both endpoints in
one of the double curves, and the other two triple points in this
double curve are endpoints of essential green arcs in the red 2bigon.
This implies that the boundary of the octagon contains two essential
green arcs, and two inessential green arcs, and the two inessential
green arcs have endpoints in the two different blue double curves.
This gives rise to the octagon illustrated in Figure \ref{picture104}.
In this case, there may be no squares in the inner solid torus, so
there are two bigons and an octagon in the inner solid torus, and
exactly three squares faces in the outer solid torus. The squares are
meridional discs for the outer solid torus, as illustrated in Figure
\ref{picture107}. However, the blue surface is null-homologous, so
there must be an even number of blue squares, a contradiction.
\end{Proof}

So we have shown that the annuli consist of a pair of 1bigons. In
particular, this means that the faces in the red torus consist of a
pair of bigons, a pair of hexagons, and all other faces are squares.
Furthermore, note that the red bigons each lie on different sides of
both the green and blue double curves, so the blue bigons lie on
different sides of the red torus. Similarly, the blue hexagons lie on
different sides of the red torus.

We now show there are more than six triple points. If there are only six triple points, the diagram is the one shown in Figure \ref{picture34} below.

\begin{figure}[H]
\begin{center}
\epsfig{file=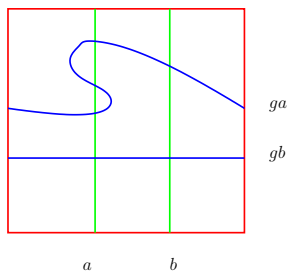, height=100pt}
\end{center}
\caption{Six triple points.}
\label{picture34}
\end{figure}

Recall from \cite[Section 2.5]{mr} that each triple point may be given
a sign, either positive or negative, which is preserved by $g$, and
which has the property that triple points which are adjacent to each
other in a double curve have opposite signs. There is only one point
in $b \cap gb$, and its image lies in $gb$. The triple point $b \cap
gb$ can't be a fixed point, so this triple point gets mapped to $a
\cap gb$, but this is adjacent to $b \cap gb$ in the blue double curve
$gb$, and so has the opposite sign, a contradiction. Therefore, there
must be at least $12$ triple points.

The green 1bigon divides the solid torus bounded by the red torus into two regions, one of which is the inner solid torus. As there are at least $12$ triple points, there must be at least two blue squares on each side of the green 1bigon. We now construct all possible squares in these two regions. 

Consider squares in the inner solid torus. If both red arcs are
essential in the green annulus, then both green arcs must be
inessential in the red 1bigon, but then the two inessential green arcs
meet different boundary components of the 1bigon, a contradiction. If
both red arcs are inessential, then as the green annulus is a 1bigon,
they must both have endpoints on a common blue double curve, and be
parallel. However, this means that both green arcs must also be
inessential and parallel. In order for the union of the arcs to form a
simple closed curve, each red arc must have endpoints in a
different green arc, but then the resulting simple closed curve is
homotopic to a core curve of the solid torus, and so does not bound a
disc inside the solid torus. The remaining case is when there is one
essential red arc, and one inessential red arc, and this is
illustrated below on the left hand side of Figure \ref{picture170}.

Consider squares in the outer solid torus. If both red arcs in the boundary of the square are inessential, then the green arcs must also be inessential, and as before, this means that the boundary of the square is isotopic to the blue double curve. If this happens, then the blue essential double curve is a meridian curve for the solid torus, and if this happens on both sides, then the manifold is $S^2 \cross S^1$, so we may assume we have chosen a side of the red torus where the blue double curve is not a meridian, and so there are no squares in the outer region with both red boundary arcs inessential. If there are one or two essential arcs, then the possible squares are illustrated below in the center and right hand diagrams in Figure \ref{picture170}. In both of these two cases the square forms a compressing disc for the outer region, so by irreducibility, the outer region is a solid torus.

\begin{figure}[H]
\begin{center}
\epsfig{file=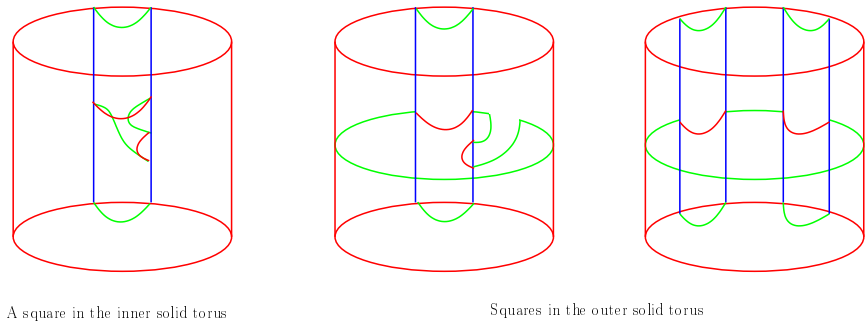, height=120pt}
\end{center}
\caption{Squares in the inner and outer solid tori.}
\label{picture170}
\end{figure}

There are therefore two types of mutually exclusive square that may occur in the outer region, and we consider each case in turn.

\begin{subcase}{The squares in the outer solid torus have essential boundary arcs only.}
\end{subcase}

Consider the case in which the boundaries of the outer squares contain essential red arcs only, as illustrated in the right-most diagram in Figure \ref{picture170}. There are at least two squares in the inner solid torus, and these have at least two inessential red arcs in their boundaries. These inessential red arcs cannot lie in squares in the outer region, and at most one essential red arc lies in the blue bigon, so the blue hexagon must be contained in the outer region. 

If the blue bigon is contained in the inner solid torus, then both boundary arcs of the bigon are essential, so if there are $n$ squares in the inner solid torus, then there are $n+1$ essential red arcs, and $n$ inessential red arcs. This means there are $2n+1$ red arcs in total, and so there are $4n+2$ triple points, which must be divisible by six, and at least $12$, so $n$ must be at least four. This means there are at least four inessential red arcs, but this is a contradiction, as all faces in the outer solid torus have essential red arcs in their boundaries, except for at most three inessential arcs in the hexagon.

If the blue bigon is contained in the outer region, then there are only squares in the inner solid torus. If there are $n$ such squares, then there are $n$ essential red arcs and $n$ inessential red arcs, giving a total of $2n$ red arcs and $4n$ triple points, so $n$ must be at least three. The blue bigon must have inessential boundary, and the hexagon has at most three red inessential arcs in its boundary, so $n=3$, and the hexagon has two inessential red arcs, and one essential red arc in its boundary. However, this means that the hexagon intersects the blue double curve with fewer triple points precisely once, so is a disc with $\pm 1$ intersection number with an essential blue double curve, which is a contradiction, as every disk in the outer solid torus has even intersection number with the essential blue double curves.

\begin{subcase}{The squares in the outer solid torus have inessential boundary arcs.}
\end{subcase}

Now suppose that the outer squares contain both essential and
inessential red arcs in their boundaries. 

We now show that the boundary of the blue bigon consists of a pair of
essential arcs, no matter which side it lies on. Suppose that the blue
bigon has both boundary arcs inessential. If it lies in the inner
solid torus, then the presence of a square in the inner solid torus
forces the red annulus in the boundary of the inner solid torus to be
a 2bigon. If it is in the outer solid torus, then the green arc in the
bigon must cross the green arcs of a square in the outer solid torus.

The squares in both the inner and outer regions all have exactly one
essential arc and one inessential arc, and the total number of
essential red arcs of faces in the inner solid torus must be the same
as the total number of essential red arcs of faces in the outer
region. So if the blue hexagon and the blue bigon lie on the same side
of the green 1bigon, then the hexagon must have two inessential red
arcs, and one essential red arc. Otherwise, if the blue hexagon and
the blue bigon lie on different sides of the green 1bigon, then the
blue hexagon must have two essential red arcs and one inessential red
arc.

\begin{claim}The hexagon has two essential red arcs.
\end{claim}

\begin{Proof}
Suppose there is a hexagon with exactly one essential red arc. The two inessential red arcs are parallel, and the innermost one has endpoints in two distinct green arcs. Removing the innermost red arc and connecting the resulting free endpoints of the two green arcs creates a blue square with one essential red arc and one inessential red arc, as illustrated in the leftmost diagram in Figure \ref{picture170}. This operation is equivalent to doing a boundary compression of the disc along a disc parallel to the green disc bounded by the innermost red arc and the corresponding innermost arc of the blue double curve with the same endpoints, and discarding the resulting component with two triple points. This means that the hexagon may be recovered by the reverse process, i.e. choose an arc in the red torus with one endpoint in a green double arc component of the boundary of the square, and the other endpoint in the innermost blue double arc bounded by the inessential red double arc in the boundary of the square, and then do a ``disc slide'' of the blue disc along this arc. Furthermore, this arc must be chosen so that the disc slide does not create any green bigons which are not parallel to the green inessential boundary arc in the blue square. There are only two possible ways of doing this, and they are illustrated in Figure \ref{picture191} below.

\begin{figure}[H]
\begin{center}
\epsfig{file=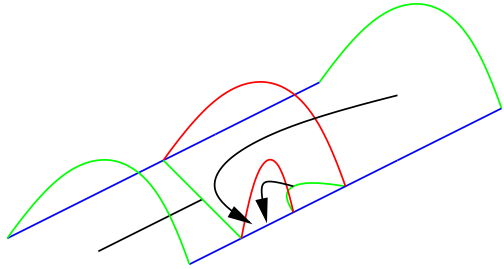, height=100pt}
\end{center}
\caption{``Disc slide'' arcs for a blue square.}
\label{picture191}
\end{figure}

Figure \ref{picture191} above shows part of the region bounded by the red and green annuli; the two ends of these surfaces should be identified to form a torus.

The resulting two hexagons with exactly one essential arc are illustrated in Figure \ref{picture171} below.

\begin{figure}[H]
\begin{center}
\epsfig{file=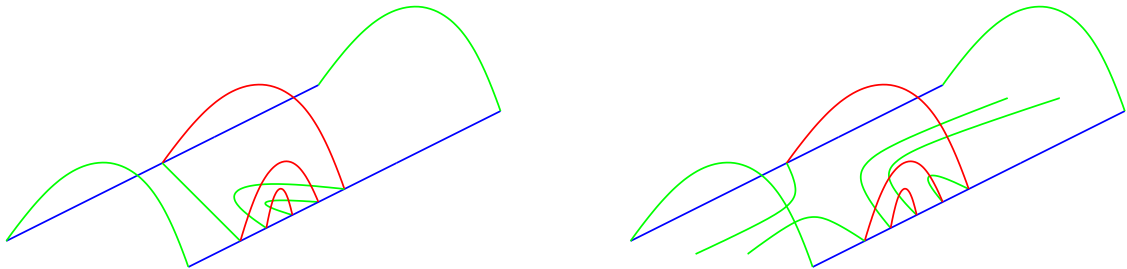, height=100pt}
\end{center}
\caption{Hexagons in the inner solid torus with one essential red arc.}
\label{picture171}
\end{figure}

In the left hand case, there may be no squares in the solid torus, so there are eight triple points in total, a contradiction, while in the right hand case, there may be no essential bigon in the solid torus, which is also a contradiction.
\end{Proof}

So we may assume that the hexagon has two essential arcs. In the left hand side of Figure \ref{picture183} below we have drawn the only possible hexagon with exactly two essential arcs, up to Dehn twists in a curve in the red 1bigon parallel to the blue double curves. As there must also be squares in the solid torus, there must in fact be exactly one Dehn twist, as illustrated on the right hand side of Figure \ref{picture183}. We have drawn the square face with black boundary arcs, in order to distinguish it from the hexagon.

\begin{figure}[H]
\begin{center}
\epsfig{file=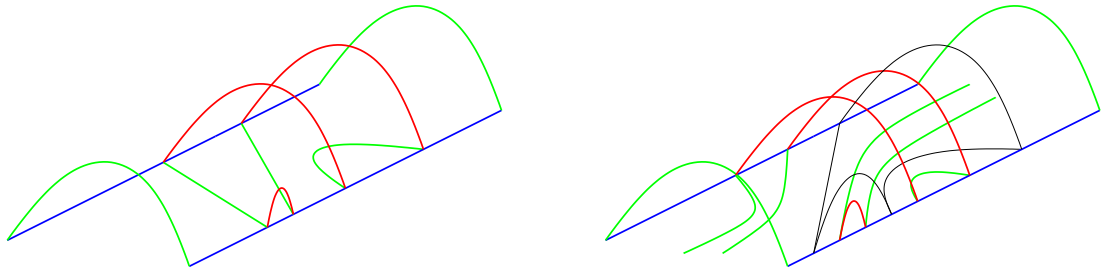, height=100pt}
\end{center}
\caption{Hexagons with two essential red arcs.}
\label{picture183}
\end{figure}

If there are $n$ parallel squares in the solid which contains the hexagon, then in the other solid torus, the blue faces consist of an essential blue bigon, and $n+1$ parallel squares. This gives the following pattern of green arcs in the red torus. The shaded regions correspond to parallel green arcs; the number of parallel green arcs is indicated by the label in the region.

\begin{figure}[H]
\begin{center}
\epsfig{file=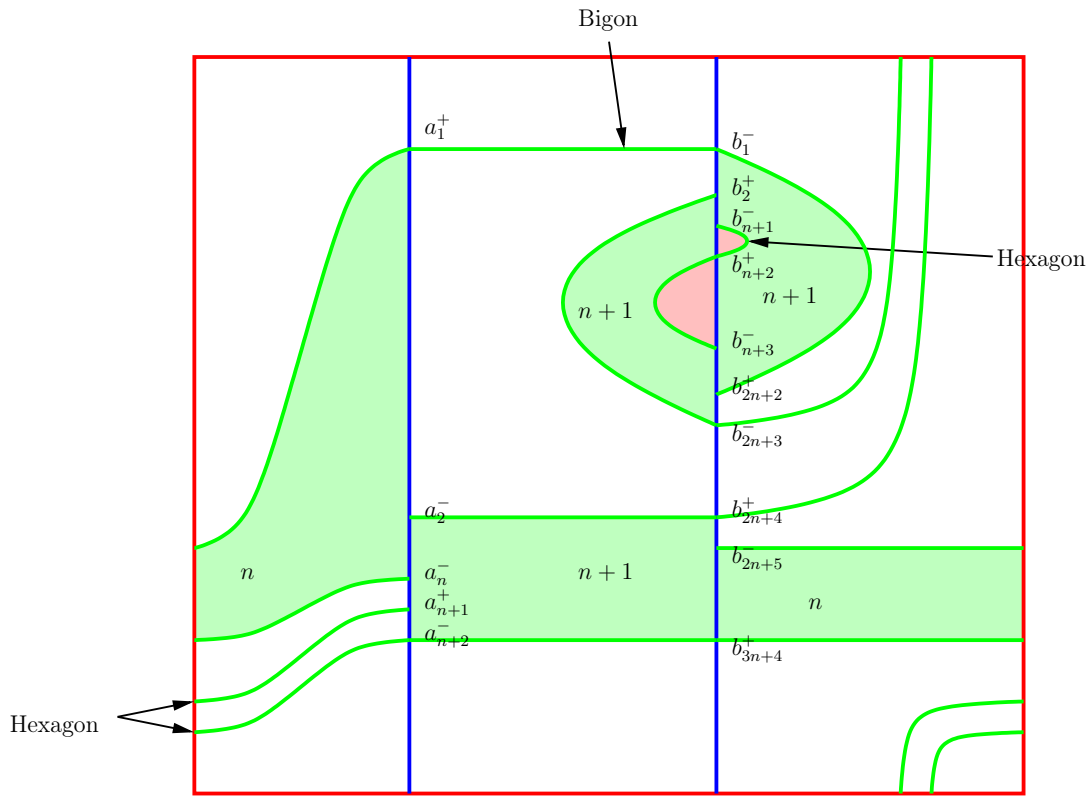, height=300pt}
\end{center}
\caption{The green arcs in the red torus.}
\label{picture184}
\end{figure}

We will label the blue double curve with fewest triple points $a$, and the other blue double curve $b$. We will give the triple points on $a$ labels $a_1, \ldots, a_{n+2}$, starting with the the triple point on $a$ which lies in the essential green arc of the bigon. Similarly we will give the triple points on the double curve $b$ the labels $b_1, \ldots b_{3n+4}$, starting with the triple point which lies in the essential green arc of the bigon. Each vertex has a sign, and we will label $a_1$ as $+$, and this gives a sign $\pm$ to the other vertices, depending on whether they are an even or odd distance from $a_1$, and we will indicate the sign with a superscript. The two red bigons share a triple point, and the triple points in their boundaries are $\{ b^-_{n+1}, b^+_{n+2}, b^-_{n+3}\}$. One of the green bigons has a red arc corresponding to the inessential red arc in the hexagon, so $\{ b^-_{2n+3}, b^+_{2n+4}\}$ lie in the boundary of the green hexagons. The third vertex is either $b^+_{2n+2}$ or $b^-_{2n+5}$, and so in fact must be $b^-_{2n+5}$, as the sign is preserved by $g$, so there are two negative vertices and one positive vertex. One of the blue essential bigons has vertices $\{ a^+_1, b^-_1\}$, so the third vertex is either $b^+_{2n+2}$ or $b^-_{2n+5}$, so as it must have negative sign, it must in fact be $b^-_{2n+5}$. However, this means that the green and blue bigons share a vertex, but the blue and red bigons do not share a vertex, a contradiction.
\end{case}

\subsection{Disjoint bigons for saddle moves.} \label{section:saddle}

In this section we show that there are disjoint bigons for saddle
moves. The only result we use from Sections \ref{section:intersection}
and \ref{section:three bigons} is Lemma \ref{lemma:three bigons}, that
there are at least three bigons.

\begin{lemma} \label{lemma:torus_saddle}
There are disjoint bigons for saddle moves.
\end{lemma}

First we set up some notation, and prove some preliminary lemmas.  We
will always call the singular curve $a$. We will draw the curve
involved in the saddle move as either an eyeglass curve, i.e. as two
smooth circles connected by an arc, or as a theta curve, i.e. a single
smooth curve, together with an arc with both endpoints in the curve.
The connecting arc we will call the {\bf saddle arc}, and we will
refer to the one or two smooth circle components as {\bf loops}. The
blue singular curve is the image of the green singular curve under
$g$, so the green curve is an eyeglass curve if and only if the blue
curve is an eyeglass curve. Recall that we call the union of the blue
and green double curves with triple points the blue-green diagram. We
will refer to the union of the blue-green diagram and the blue and
green saddle arcs, as the {\bf saddle diagram}. We will call the union
of a blue arc and a green arc which meet exactly at their endpoints
and contain no other triple points a {\bf bigon-boundary}. A
bigon-boundary may not bound a bigon as it may not necessarily bound a
disc.  If we orient the green and blue arcs in the bigon-boundary,
then each point of intersection has opposite orientation, as if they
had the same orientation, then a regular neighbourhood of the
bigon-boundary would have a boundary component which is a simple
closed curve in the red torus that hits the union of the green curves
precisely once. However, the union of the green curves is
null-homologous in the red torus, so this cannot happen. So a regular
neighbourhood of the bigon-boundary in the red torus is an annulus
with one boundary curve disjoint from the double curves.

By Lemma \ref{lemma:three bigons}, the configuration both before and
after the saddle move contains at least three bigons. So we may assume
that each saddle arc in the saddle diagram has at least one endpoint
on a bigon. A third red bigon may also have an endpoint of either the
green or blue saddle arc meeting its boundary.  The interior of the
third bigon either intersects the interior of the saddle arc, or else
meets the saddle arc only in its endpoint. If the interior of the
saddle arc meets the interior of the bigon, then the third bigon is
adjacent to one of the first two bigons. These two different cases are
illustrated below in Figure \ref{picture147}.

\begin{figure}[H]
\begin{center}
\epsfig{file=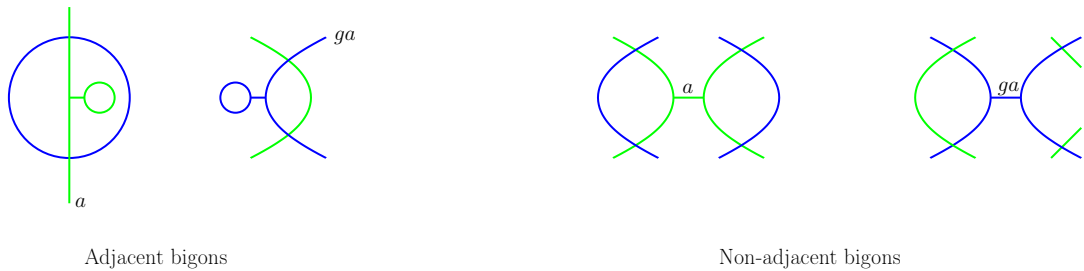, height=100pt}
\end{center}
\caption{Adjacent and non-adjacent bigons.}
\label{picture147}
\end{figure}

In the adjacent bigons case we have drawn the case in which the
singular curve is an eyeglass. There is a similar picture in which the
saddle arc and the loop with no triple points are replaced by a saddle
arc with both endpoints in the loop with triple points, forming a
theta curve. We will refer to both configurations as the case of
adjacent bigons. In the non-adjacent bigons picture the part of the
singular curve shown may lie in either a theta curve or an eyeglass.

We deal with the two cases of adjacent or non-adjacent bigons
separately. However, in each case the strategy is the same. First we
show that the orbits of the triple points closest to the saddle arcs
are all distinct, and construct the local configuration of the saddle
diagram near these points. We then show that the saddle diagram may
not have more than one connected component, and then finally
eliminate this case.

\subsubsection{Adjacent bigons.}

We start by showing that the orbits of the triple points closest to
the saddle arc are distinct.

\begin{lemma}
If the saddle arcs meet adjacent bigons then the orbits of the
triple points closest to the green saddle arc are all distinct.
\end{lemma}

\begin{Proof}
We may assume that the adjacent bigons share a common green arc, and
we will label the blue double curve made up of the blue arcs of the
adjacent bigons $gc$. The double curve $gc$ contains exactly two
triple points, which we will call $x$ and $y$. This is illustrated
below in Figure \ref{picture203}.

\begin{figure}[H]
\begin{center}
\epsfig{file=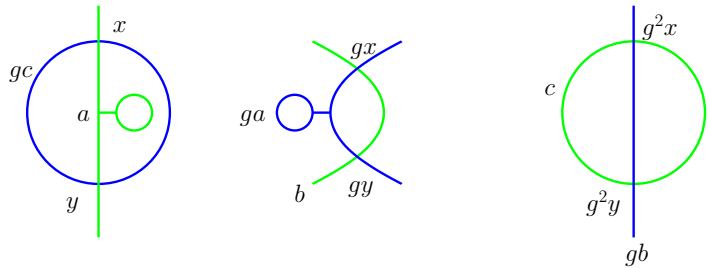, height=100pt}
\end{center}
\caption{The local configuration for adjacent bigons.}
\label{picture203}
\end{figure}

If the triple points closest to the saddle arcs, namely $x$, $y$, $gx$
and $gy,$ are not all distinct, then the pair of points $\{x,y\}$ is
the same as the pair of points $\{gx,gy\}$, as the blue double curve
$gc$ is a loop of the eyeglass curve $ga$. This means there is a pair
of points invariant under $G$, a contradiction. So we may assume that
the points in the orbits of $x$ and $y$ are all distinct. This also
implies that the double curve $c$ is not one of the loops of $a$.
Figure \ref{picture203} above shows the blue-green diagram close to
the orbits of $x$ and $y$.
\end{Proof}

We now show that the saddle diagram is connected.

\begin{lemma} \label{lemma:adjacent connected}
Suppose the saddle diagram has more than one component. Then there is
a bigon disjoint from the saddle arcs.
\end{lemma}

\begin{Proof}
We divide the proof of this lemma into two steps. We first show that
no connected component of the saddle diagram may be contained in a
disc. We then argue that the only other possibility is that the saddle
diagram consists of two components, both contained in annuli, and then
deal with that case.

\setcounter{case}{0}
\begin{case}{Suppose the saddle arcs meet adjacent bigons, and a component of the
saddle diagram is contained in a disc. Then there is a bigon disjoint
from the saddle arcs.}

A pair of intersecting curves in a disc creates at least three bigons.
A single image of a saddle arc may intersect at most two bigons, so
both the green and blue saddle arcs must be contained in the same connected
component of the saddle diagram, lying in the disc.

The green double curve $c$, together with the arc of $gb$ between
$g^2x$ and $g^2y$ with no triple points in its interior, creates two
bigon-boundaries, as illustrated in Figure \ref{picture203}. As we
have assumed that the saddle diagram is contained in a disc, at most
one of these bigon-boundaries may not bound a disc. In order to avoid
a disjoint bigon, the bigon-boundary which bounds a disc must contain
either a green or blue saddle arc in its interior, contained in a
connected component of the saddle diagram. But this means that the
different saddle arcs live in different connected components of the
saddle diagram, so there is a disjoint bigon.
\end{case}

If any connected component of the saddle diagram is contained in a
punctured torus, then any other component must be contained in a disc,
which we have just shown may not occur. Any collection of intersecting
double curves contained in an annulus creates a bigon, so if there are
three or more connected components of the saddle diagram contained in
annuli then there is a disjoint bigon, as there are only two saddle
arcs. So the final disconnected case to consider is if the saddle
diagram has two components, each contained in an annulus.

\begin{case}{Suppose the saddle arcs meet adjacent bigons, then the saddle diagram
may not contain two connected components, both contained in annuli.}

The leftmost picture in
Figure \ref{picture203} contains a pair of bigon-boundaries. If either
bounds a disc, then it is a bigon disjoint from the saddle arcs, as we
have shown that no connected component of the saddle-diagram is
contained in a disc. Therefore, neither bounds a disc, and the only
way this can happen is if they are tubed together to form an annulus
with a pair of bigon-boundaries, as illustrated below in Figure
\ref{picture148}.

\begin{figure}[H]
\begin{center}
\epsfig{file=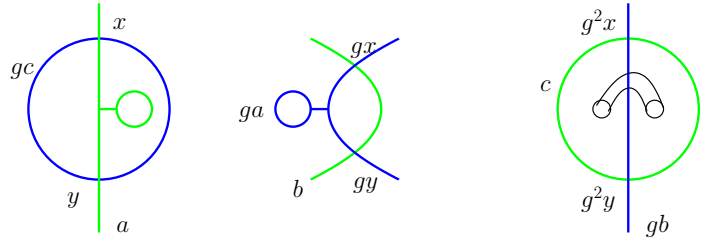, height=90pt}
\end{center}
\caption{The local configuration when the saddle diagram has two annular components.}
\label{picture148}
\end{figure}

This implies that the green double curve $c$ bounds a punctured torus
on one side, and therefore must bound a disc on the other side, so the other
component of the saddle diagram is contained in the annulus with two
bigon-boundaries.

The blue double curve $gb$ divides the disc bounded by the green curve
$c$ into two discs, each of which contains at least one bigon,
separated by the blue curve. In order for both of these bigons to
intersect a single saddle arc, the blue saddle arc must lie in this component,
and $gb$ must be the loop of $ga$ containing triple points. This is
illustrated below in the left hand side of Figure \ref{picture195}.

\begin{figure}[H]
\begin{center}
\epsfig{file=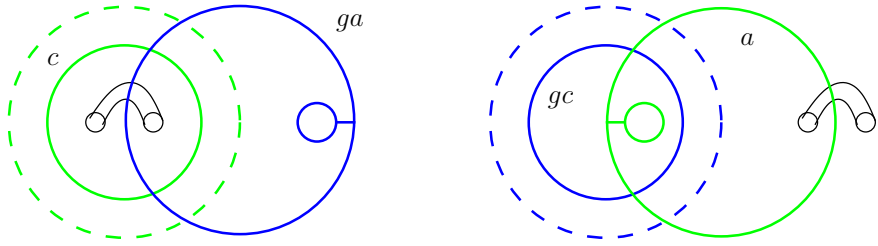, height=90pt}
\end{center}
\caption{The saddle arcs contained in two annular components of the saddle diagram.}
\label{picture195}
\end{figure}

In the left hand side of Figure \ref{picture195} above we have chosen
to draw the saddle arc on one particular side of the blue double curve
$ga$, but it makes no difference which side it lies on.

Any green arcs with endpoints in $ga$ must be parallel to the green
arcs of $c$, and so must be parallel to the green curve $c$. Any other
blue double curves create disjoint bigons, so there may be no other
blue double curves with triple points. Therefore the left hand side of
Figure \ref{picture195} shows all possible blue and green double
curves in this annular component of the saddle diagram.  The dotted
green curve represents some number, which may be zero, of parallel
green double curves.

The green saddle arc lies in the other connected component of the
saddle diagram. The blue double curve $gc$ bounds a disc containing a
pair of bigons, so is inessential. If the green double curve $a$ is
also inessential, then it contains a disjoint bigon, so it must be
essential. This is illustrated above in the right hand side of Figure
\ref{picture195}. Any blue double arcs with endpoints in $a$ must be
parallel to $gc$, and so all blue curves in this component are
parallel to $gc$. Any other green double curves contained in this
component create disjoint bigons, so there may be no other green
double curves. Therefore the left hand side of Figure \ref{picture195}
shows all possible green and blue double curves in the connected
component of the saddle diagram containing the green saddle arc. The
blue dotted curve represents some number, which may be zero, of
parallel blue double curves.

The configuration has the property that all triple points lie on
either $a$ or $ga$. As each triple point has three orbits under $G$,
this means that two of them must lie on the same double curve. By
choosing labels appropriately, we may assume that $x$ and $gx$ lie in
the same double curve. Suppose that $x$ and $gx$ lie on the green
double curve $a$. But if $x$ is in $a$, then $gx$ lies in $ga$, so the
double curves $a$ and $ga$ intersect, a contradiction. Similarly, if
$x$ and $gx$ lie in the blue double curve $ga$, then as $gx$ is in
$ga$, this implies that $x$ is in $a$, so again the double curves $a$
and $ga$ must intersect, a contradiction.
\end{case}

This completes the proof of Lemma \ref{lemma:adjacent connected}.
\end{Proof}

We have shown that there are disjoint bigons if the saddle diagram is
disconnected. So we now assume that the saddle diagram is connected,
and contains the local configurations shown in Figure \ref{picture203}.

The green double curve $c$, together with the arc of $gb$ between
$g^2x$ and $g^2y$ with no triple points in its interior, creates two
bigon-boundaries. By Lemma \ref{lemma:adjacent connected} the saddle
diagram is connected, so they cannot bound discs, so they are tubed
together, as illustrated previously in Figure \ref{picture148}
above.

The green and blue saddle arcs lie in the disc bounded by $c$. There
is a blue arc of $gb$ which divides this disc into two discs, each of
which contains a bigon. If the double curve $gb$ is not part of $ga$,
then the green and blue saddle arcs must lie on different sides of
$gb$ in the disc bounded by $c$. But then the loop of $ga$ containing
triple points is inessential, so contains two bigons, at least one of
which is disjoint from the saddle arcs. So we may assume that $b$ is a
loop of $a$.

This means that the blue loop of $ga$ with triple points is essential,
and the green loop of $a$ with triple points is disjoint from $c$, and
hence contained in the disc bounded by $c$, and so is inessential, and
bounds a disc containing the blue saddle arc in its interior. This is
illustrated in Figure \ref{picture149} below. There may be many double
curves parallel to the dotted curves. There must be an even number of
essential double curves of each colour, as the union of the double
curves of a single colour is null-homologous, so there must be at
least one blue curve without triple points. The saddle arcs may lie on
either side of the loops with triple points, this makes no difference
to the argument that follows.

\begin{figure}[H]
\begin{center}
\epsfig{file=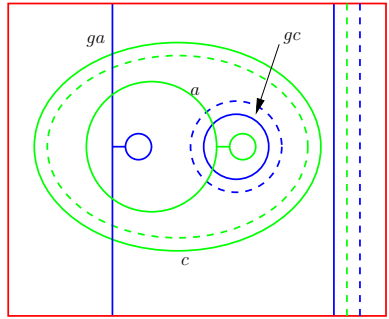, height=120pt}
\end{center}
\caption{The eyeglass $ga$ has an essential loop.}
\label{picture149}
\end{figure}

Consider the configuration after the saddle move. The green curves
divide the red torus into a disc, some number of annuli, and a
punctured torus. The red disc bounded by the green curve $a$ gets
mapped to a green disc bounded by the blue curve $ga$, so $ga$ is a
meridian on this side of the red torus. The adjacent red annulus
between $a$ and the next green curve gets mapped to a green annulus
with blue boundary, one of whose boundary curves is $ga$. This green
annulus lies on the other side of the red torus to the green disc
whose boundary is $ga$. The blue curve $ga$ does not bound a disc on
the other side of the red torus to the green disc, as $L$ is not $S^2
\cross S^1$, so this green annulus has both blue boundary components
essential and with the same slope on the red torus. But only one
essential blue curve has triple points, and the annulus has triple
points on both components, a contradiction.

This completes the proof of Lemma \ref{lemma:torus_saddle} in the case
that the saddle arcs meet adjacent bigons. 

\subsubsection{Non-adjacent bigons.}

We now consider the case in which the saddle arcs meet non-adjacent
bigons.  We start by showing that the orbits of the triple points
closest to the saddle arc are distinct.

\begin{lemma} \label{lemma: non-adjacent distinct orbits}
Suppose the saddle arcs meet non-adjacent bigons. Then the orbits of
the triple points closest to the green saddle arc are all distinct.
\end{lemma}

It will be convenient to give a name to the union of the saddle arc
and the two segments of the loops between the closest triple points to
the saddle arc.

\begin{definition}{\bf The green $H$ and the blue $H$, and the $H$-points.}

Consider the green singular curve $a$, minus its triple points. Call the
closure of the connected component of this which contains the saddle
arc the {\bf green $H$}. Define the {\bf blue $H$} similarly. We shall
refer to the endpoints of the $H$ as the \bf{$H$-points}.
\end{definition}

\begin{Proof}
We start by showing that the endpoints of each $H$ are distinct from
each other. The triple points alternate in sign going clockwise around the
endpoints of the green $H$, and $G$ preserves the sign of the triple
points. If opposite points are the same, then the green arcs are part
of a single essential curve which is saddled to itself by a saddle arc
which is essential in the annulus formed by the complement of the
green double curve in the red torus. But this means that the union of
the green curves is a non-zero homology class in the red torus, a
contradiction.

We now wish to show that in fact the orbits of all the $H$-points are
distinct. Figure \ref{picture150} below shows the configuration close
to the orbits of the $H$-points.

\begin{figure}[H]
\begin{center}
\epsfig{file=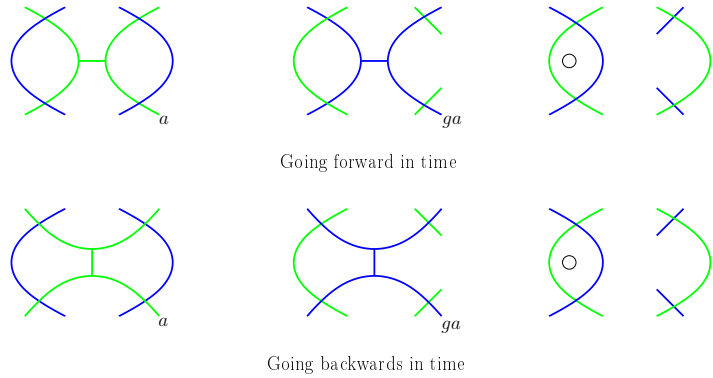, height=140pt}
\end{center}
\caption{The orbits of the $H$-points are all disjoint.}
\label{picture150}
\end{figure}

The top row of Figure \ref{picture150} shows the saddle move going
forward in time, while the bottom row shows it going back in time. The
rightmost bigon-boundary cannot bound a disc, or else there is a
disjoint bigon, so it is tubed to some other face. The bottom
configuration also contains at least three bigons, which may not be
disjoint from the saddles, or else there is a disjoint bigon. So one
of the saddle arcs has its endpoints on two bigons, and the other on
at least one. There are two different ways in which this can happen,
illustrated in Figure \ref{picture151} below.

\begin{figure}[H]
\begin{center}
\epsfig{file=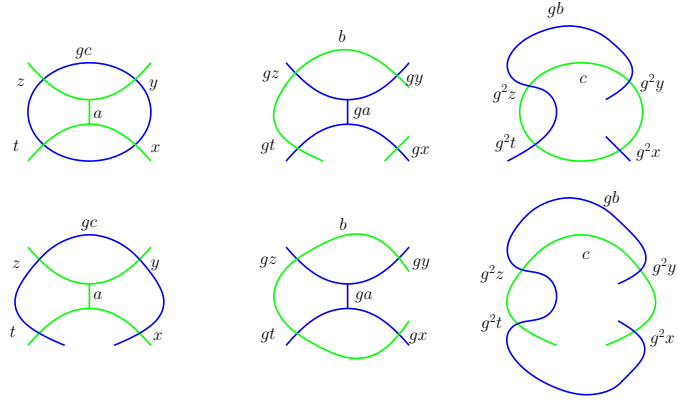, height=150pt}
\end{center}
\caption{The two different local configurations for non-adjacent bigons.}
\label{picture151}
\end{figure}

We now deal with different cases depending on how the orbits of the
endpoints of the $H$ may overlap. We label the triple points $\{
x,y,z,t \}$ as illustrated in Figure \ref{picture151}. If the
$H$-points are not all distinct, then, in particular, the sets of
points $\{ x,y,z,t\}$ and $\{ gx,gy, gz, gt \}$ must have an element
in common. Recall that each triple point has a sign, which is
preserved by $G$, and triple points which are adjacent along a double
curve have opposite sign. This means that if $x$ is an element of $\{
gx, gy, gz, gt \}$, then it may only be equal to $gx$ or $gz$. As
there are no fixed points, this means that if $x$ is a member of $\{
gx,gy, gz, gt \}$, then it must be $gz$. Analogous arguments show that
the only possible identifications between elements of $\{ x,y,z,t\}$
and $\{ gx,gy, gz, gt \}$ are $x=gz, y=gt, z=gx$ or $t=gy$.

We now deal with the upper and lower configurations in turn. In each
case, the basic strategy is to observe that if two triple points are
the same, then their adjacent triple points along a double curve are
also the same. Hence any pair of identifications between the triple
points in the diagrams gives rise to extra identifications among the
triple points which are adjacent to it along double curves.  Repeating
this procedure, we eventual find so many identifications that it gives
rise to a triple point that is its own image, i.e. a fixed point,
which is a contradiction.

\setcounter{case}{0}
\begin{case}{The upper diagram in Figure \ref{picture151}.}

First suppose that $x = gz$. In the top left diagram of Figure
\ref{picture151}, the triple point $x$ is adjacent to $y$ and $t$
along a blue double curve. In the top middle diagram of
\ref{picture151}, the triple point $gz$ is adjacent to $gy$ along a
blue double curve. This implies that $gy$ is equal to either $y$ or
$t$, but as there are no fixed points, this implies $t = gy$. Applying
$g$, this implies $gt = g^2y$. In the top middle diagram of Figure
\ref{picture151}, the triple point $gt$ is adjacent to $gz$
along a green double arc, while in the top right diagram of Figure
\ref{picture151}, $g^2y$ is adjacent to the triple points $g^2x$ and
$g^2z$ along green double arcs. As there are no fixed points, this implies that $gz =
g^2x$. But we assumed that $x = gz$, so this gives $x = g^2x$, giving
a fixed point, a contradiction.

Now suppose that $y=gt$. The triple point $y$ is adjacent to $x$ and
$z$ along blue arcs, and $gt$ is adjacent to $gx$. As there are no
fixed points, this implies that $z=gx$. Applying $g$, this implies
that $gz = g^2x$. The triple point $gz$ is adjacent to $gy$ and $gt$
along green arcs, whereas the triple point $g^2x$ is adjacent to
$g^2y$ and $g^2t$ along green arcs. This implies that the pair of
points $\{gy, gt\}$ is equal to the pair of points $\{ g^2y,g^2t\}$,
but this gives a pair of points fixed by $g$. As $g$ has order three,
this implies that $g$ has fixed points, a contradiction.

Now suppose that $z=gx$. The triple point $z$ is adjacent to $y$ and
$t$ along blue arcs, while $gx$ is adjacent to $gt$ along a blue
arc. As there are no fixed points, this implies $y=gt$, which is the
case we have just dealt with.

Finally suppose that $t=gy$. The triple point $t$ is adjacent to $x$
and $z$ along blue arcs, while $gy$ is adjacent to $gz$ along a blue
arc, so as there are no fixed points, $x=gz$, which was the first case
we considered.
\end{case}

\begin{case}{The lower diagram in Figure \ref{picture151}.}

It will be convenient to consider the cases in a different order.

First, suppose that $t=gy$. This implies that $gt=g^2y$. The triple
point $g^2y$ is adjacent to $g^2x$ and $g^2z$ along green arcs, so as
$g^2y$ is the same triple point as $gt$, this implies that the triple
points $gz$ and $gx$, which are adjacent to $gt$ along green arcs are
the same as the pair of triple points $g^2x$ and $g^2z$. However, this
implies that the pair of points $gz$ and $gx$ are $G$-invariant, but
as $g$ has order three, this means they are fixed points, a
contradiction.

Now suppose that $x = gz$. The triple point $x$ is adjacent to $t$
along a green double curve, and the triple point $gz$ is adjacent to
$gy$ and $gt$ along a green double curve. As there are no fixed
points, this means that $t = gy$. But this is the case we have just
considered.

Now suppose that $y = gt$. The triple point $y$ is adjacent to $z$
along a green arc, while $gt$ is adjacent to $gz$ and $gx$ along green
double arcs, so as there are no fixed points, this implies that
$z=gx$.  Applying $g$, this implies that the triple point $gz = g^2x$.
The triple point $gz$ is adjacent to the triple points $gy$ and $gt$
along green arcs, while the triple point $g^2x$ is adjacent to $g^2y$
along green arcs. As there are no fixed points, this implies that
$gt=g^2y$. However, as $y = gt$, this implies that $y = g^2y$, giving
a fixed point, a contradiction.

Finally suppose that $z=gx$. The triple point $gx$ is adjacent along a
blue arc to $gt$, so $gt$ is equal to either $t$ or $y$. As there are
no fixed points, this implies that $y=gt$, but this is the previous
case we have discussed.
\end{case}

This completes the proof of Lemma \ref{lemma: non-adjacent distinct orbits}.
\end{Proof}

We now show that the saddle diagram is connected. As in the case of
adjacent bigons, we start by showing that no component of the saddle
diagram may be contained in a disc. We then argue that the only
remaining disconnected case is when there are two components, each
contained in an annulus, and we then deal with that case.

\begin{lemma} \label{lemma:non-adjacent connected}
Suppose the saddle arcs meet non-adjacent bigons. Then the saddle
diagram is connected.
\end{lemma}

\begin{Proof}
We first show that no connected component of the saddle diagram may be
contained in a disc.

\setcounter{case}{0}
\begin{case}{Suppose the saddle arcs meet non-adjacent bigons, and a component of the
saddle diagram is contained in a disc. Then there is a bigon disjoint
from the saddle arcs.}

A pair of intersecting curves in a disc creates at least three bigons.
A single image of a saddle arc may intersect at most two bigons, so
both the green and blue saddle arcs must be contained in the same connected
component of the saddle diagram, lying in the disc.

In both the upper and lower cases illustrated in Figure
\ref{picture151}, the right hand diagrams contain a pair of
bigon-boundaries, formed by arcs of $c$ and $gb$, at most one of which
may not bound a disc. A bigon-boundary which bounds a disc is
either a disjoint bigon, or contains a connected component of the
saddle diagram, which either creates a disjoint bigon, or contains at
least one saddle arc. But by the argument in the previous paragraph, a
connected component of the saddle diagram contained in a disc contains
both saddle arcs, so in fact both saddle arcs are contained in the
disc bounded by the bigon-boundary, and the component of the saddle
diagram containing $c$ and $gb$ must have a disjoint bigon.
\end{case}

So we have shown that no connected component of the saddle diagram may
be contained in a disc.  Any collection of intersecting curves in an
annulus creates at least two bigons, so if there are three or more
components then there is a disjoint bigon. So there are exactly two
connected components in the saddle diagram, each of which is contained
in an annulus, and each component contains one of the blue or green
saddle arcs.

\begin{case}{Suppose the saddle arcs meet non-adjacent bigons. Then the saddle
diagram may not have two connected components, both contained in annuli.}

The two possible local configurations near the orbits of the
$H$-points is shown in Figure \ref{picture151} above. In both cases,
the right hand diagrams contain a pair of bigon-boundaries, which may
not bound discs, as this would create disjoint bigons. Therefore, they
must be tubed together to form an annular region, whose boundary
consists of a pair of bigon-boundaries, as illustrated in Figure
\ref{picture204} below.

\begin{figure}[H]
\begin{center}
\epsfig{file=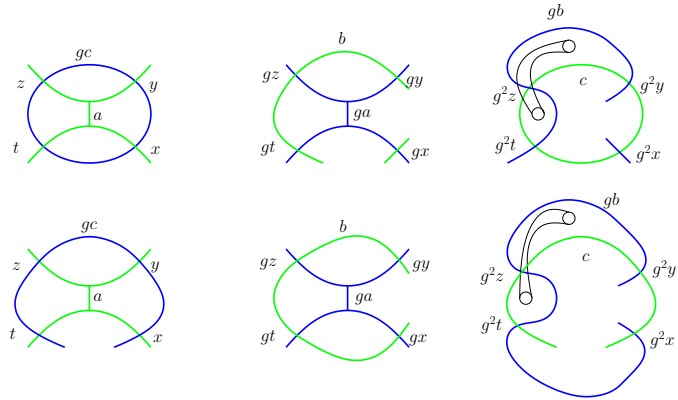, height=150pt}
\end{center}
\caption{Annuli with bigon-boundaries.}
\label{picture204}
\end{figure}

A regular neighbourhood of the annulus with boundary a pair of
bigon-boundaries is a punctured torus, so all other faces in the
annular component containing the green double curve $c$ must be
discs. However, there must therefore be bigons on both sides of the
green double curve $c$, in the complement of the annular region with two
bigon-boundaries. These bigons may meet a single saddle arc only if
the saddle arc is green. Similarly, there are bigons on both sides of
the blue double curve $gb$ in the complement of the annular region
with two bigon-boundaries, and again the only way this can happen is
if the two bigons meet the blue saddle arc. However, this implies that
both the green and blue saddle arcs are contained in a single
annular connected component of the saddle diagram, but there must be
one saddle arc in each component, a contradiction.
\end{case}

This completes the proof of Lemma \ref{lemma:non-adjacent connected}.
\end{Proof}

Finally, we consider the case in which the saddle diagram contains a
single connected component. In both cases illustrated in Figure
\ref{picture151}, the rightmost diagram contains a pair of
bigon-boundaries, neither of which may bound a disc, so they must be
tubed together to form an annulus whose boundary consists of a pair of
bigon-boundaries, as illustrated in Figure \ref{picture204}. The
annulus with bigon-boundaries may not contain any triple points in its
interior, as the saddle diagram is connected by Lemma
\ref{lemma:non-adjacent connected}. As there are no double curves with
triple points in the annulus, all curves with triple points are
inessential, except for the green and blue curves in the rightmost
diagram in Figure \ref{picture204}. We consider two cases depending on
whether or not the green double curve $c$ is a loop of the singular
green curve $a$.

\setcounter{case}{0}
\begin{case}{The green double curve $c$ is not part of the singular
    green curve $a$.}

Suppose that the green curve $c$ shown in the rightmost diagram does
not form part of the singular curve $a$. This has to happen in the top
case, as the orbits of the $H$-points are distinct, but may not
necessarily happen in the bottom case. The saddle diagram consists of
a single component contained in an annulus, and the green double curve
$c$ is essential in this annulus, and divides it into two pieces, each
of which contains at least one bigon.  As the green and blue saddle
arcs are disjoint from $c$, there must be one saddle arc on each side
of $c$. Suppose the green singular curve has two loops.  Then there is
an innermost green loop which bounds a disc containing at least two
bigons, only one of which may meet the saddle arc, so there is a
disjoint bigon. If the green singular curve contains a single loop,
then this loop bounds a disc containing at least two bigons, which
must both intersect the green saddle arc, so the green saddle arc may
not lie in the interior of the green loop. As the saddle arc lies in
the exterior of the green loop, it divides the complement of the disc
bounded by the green loop into two parts, one of which is a disc. The
blue arcs in the interior of the disc bounded by the green loop are
contained in blue double curves which may not cross the green saddle
arc, so this creates a disjoint bigon inside this disc.
\end{case}

\begin{case}{The green double curve $c$ is a loop of the singular
    curve $a$.}

Now suppose that the green double curve $c$ is one of the loops of the
singular green curve $a$. As the orbits of all of the $H$-points are
distinct, this can only happen in the lower case of Figure
\ref{picture204}. The singular curve $a$ has either one or two loops,
and we deal with each case in turn.

\begin{subcase}{The singular curve $a$ has a single loop.}
\end{subcase}

Suppose that the singular curve $a$ contains a single loop. 
Then the green saddle arc, together with a subsection of
the green loop, bounds a disc. The subsection of the green loop
between the endpoints of the saddle arc must contain triple points, so
there is a blue arc which creates a bigon inside this disc, disjoint
from the green saddle arc. This bigon must therefore meet the blue
saddle arc. But again, the blue saddle arc, together with part of the
blue loop of $ga$, bounds a disc, and there must be triple points in
between the endpoints of the blue saddle arc. Therefore the disc
bounded by the blue saddle arc and part of the blue loop contains a
disjoint bigon. This is illustrated below in Figure \ref{picture200}.

\begin{figure}[H]
\begin{center}
\epsfig{file=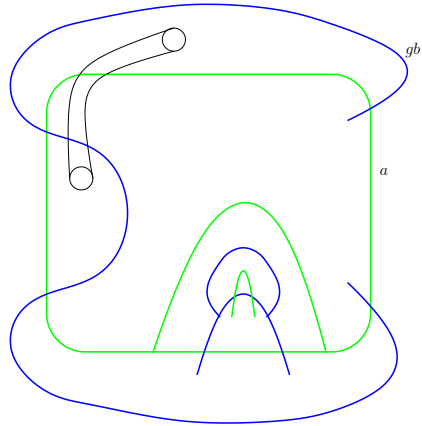, height=160pt}
\end{center}
\caption{The singular curve $a$ has a single loop.}
\label{picture200}
\end{figure}

In Figure \ref{picture200} above, we have chosen to draw the saddle
arcs as lying above the loops; this makes no difference to the
argument.

\begin{subcase}{The singular curve $a$ has two loops.}
\end{subcase}

Finally, suppose that the singular curve $a$ has two loops. There are
two ways in which this can happen, illustrated below in Figure
\ref{picture153}.

\begin{figure}[H]
\begin{center}
\epsfig{file=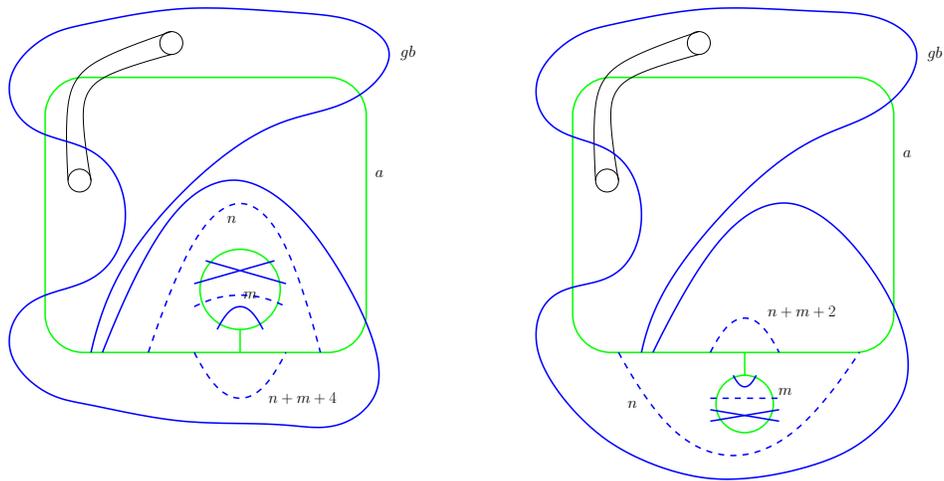, height=180pt}
\end{center}
\caption{The singular curve $a$ has two loops.}
\label{picture153}
\end{figure}

Figure \ref{picture153} above shows all possible blue arcs that may
occur. The dotted lines show that there may be many arcs parallel to
the one shown. We have drawn the saddle arcs as points, as there are
two different ways in which they may occur, either horizontal or
vertical. In both cases, the arcs of $gb$ meet $a$ on the left hand
side of the green saddle arc, and there are two collections of blue
arcs that have endpoints on either side of the green saddle arc. The
labels on the dotted arcs refer to how many curves there are parallel
to each arc, so $n$ and $m$ are non-negative integers, and the number
of dotted arcs on one side of $a$ determines the total number of
dotted arcs on the other side of $a$. In both cases, there are enough
blue arcs on the other side of the green saddle arc to create a blue
double curve (corresponding to the one labelled $gc$ in Figure
\ref{picture204}) containing precisely the four triple points closest
to the green saddle, and in particular, this blue double curve is
disjoint from the blue saddle arc. However, we have assumed that $c$
is part of the singular curve $a$, a contradiction.
\end{case}

This completes the proof of Lemma \ref{lemma:torus_saddle}.



\end{document}